\pdfoutput=1
\RequirePackage{ifpdf}
\ifpdf 
\documentclass[pdftex]{sigma}
\else
\documentclass{sigma}
\fi

\numberwithin{equation}{section}

\usepackage[all, 2cell]{xy}
\UseTwocells
\definecolor{e-mail}{rgb}{0,.40,.580}
\definecolor{reference}{rgb}{0,.40,.80}
\definecolor{citation}{rgb}{0,.40,.80}

\usepackage{mathtools}
\usepackage{tikz-cd}
\usepackage{eucal}

\newtheorem{Theorem}{Theorem}[section]
\newtheorem*{Theorem*}{Theorem}
\newtheorem{Corollary}[Theorem]{Corollary}
\newtheorem{Lemma}[Theorem]{Lemma}
\newtheorem{Proposition}[Theorem]{Proposition}
{ \theoremstyle{definition}
	\newtheorem{Definition}[Theorem]{Definition}
	
	\newtheorem{Example}[Theorem]{Example}
	\newtheorem{Remark}[Theorem]{Remark} }

\newtheorem{thmx}{Theorem}

\newcommand{\cC}{\mathcal{C}}

\newcommand{\cM}{\mathcal{M}}

\newcommand{\cO}{\mathcal{O}}
\newcommand{\calP}{\mathcal{P}}

\newcommand{\pt}{\mathrm{pt}}
\newcommand{\CG}{\mathrm{CG}}
\newcommand{\Tr}{\mathrm{Tr}}

\newcommand{\C}{\mathbf{C}}

\newcommand{\R}{\mathrm{R}}

\newcommand{\U}{\mathrm{U}_{\hbar}}
\newcommand{\Z}{\mathbf{Z}}
\newcommand{\A}{\mathbf{A}}

\newcommand{\rY}{\mathrm{Y}_{\hbar}}
\newcommand{\GL}{\mathrm{GL}}
\newcommand{\rM}{\mathrm{M}}

\renewcommand{\b}{\mathfrak{b}}
\newcommand{\g}{\mathfrak{g}}
\newcommand{\h}{\mathfrak{h}}
\newcommand{\n}{\mathfrak{n}}
\renewcommand{\sl}{\mathfrak{sl}}

\renewcommand{\u}{\mathfrak{u}}
\newcommand{\m}{\mathfrak{m}}

\newcommand{\Wh}{\mathrm{Wh}_{\hbar}}
\newcommand{\act}{\mathrm{act}}

\DeclareMathOperator{\ad}{ad}

\newcommand{\BiMod}[2]{{}_{#1}\mathrm{BiMod}_{#2}}

\newcommand{\Mod}[1]{{\mathrm{Mod}}_{#1}}

\newcommand{\dyn}{\mathrm{dyn}}
\newcommand{\End}{\mathrm{End}}
\newcommand{\ev}{\mathrm{ev}}

\newcommand{\forget}{\mathrm{forget}}
\newcommand{\Frac}{\mathrm{Frac}}

\newcommand{\free}{\mathrm{free}}

\newcommand{\gen}{\mathrm{gen}}

\newcommand{\HC}{\mathrm{HC}_{\hbar}}

\newcommand{\Hom}{\mathrm{Hom}}

\newcommand{\id}{\mathrm{id}}

\newcommand{\RMod}{\mathrm{RMod}}

\DeclareMathOperator{\Rep}{Rep}
\newcommand{\res}{\mathrm{res}}

\newcommand{\Sym}{\mathrm{Sym}}

\newcommand{\univ}{\mathrm{univ}}

\newcommand{\Vect}{\mathrm{Vect}}
\newcommand{\wt}{\mathrm{wt}}
\newcommand{\qdet}{\mathrm{qdet}}

\newcommand{\gl}{\mathfrak{gl}}

\newcommand{\adj}[2]{
	\xymatrix{
		#1 \ar@<.5ex>[r] & #2 \ar@<.5ex>[l]
	}
}

\begin{document}
\allowdisplaybreaks

\newcommand{\arXivNumber}{2310.06669}

\renewcommand{\PaperNumber}{025}

\FirstPageHeading

\ShortArticleName{Yangians, Mirabolic Subalgebras, and Whittaker Vectors}
	
\ArticleName{Yangians, Mirabolic Subalgebras, \\ and Whittaker Vectors}
	
\Author{Artem KALMYKOV~$^{\rm ab}$}
	
\AuthorNameForHeading{A.~Kalmykov}
	
\Address{$^{\rm a)}$~Department of Mathematics, Massachusetts Institute of Technology, \\
\hphantom{$^{\rm a)}$}~Cambridge MA, 02139, USA}
\EmailD{\href{mailto:artem.o.kalmykov@gmail.com}{artem.o.kalmykov@gmail.com}}
\Address{$^{\rm b)}$~Saint Petersburg University, 7/9 Universitetskaya nab., St. Petersburg, 199034, Russia}
	
\ArticleDates{Received March 26, 2024, in final form April 04, 2025; Published online April 18, 2025}
	
\Abstract{We construct an element in a completion of the universal enveloping algebra of $\gl_N$, which we call the Kirillov projector, that connects the topics of the title: on the one hand, it is defined using the evaluation homomorphism from the Yangian of $\gl_N$, on the other hand, it gives a canonical projection onto the space of Whittaker vectors for any Whittaker module over the mirabolic subalgebra. Using the Kirillov projector, we deduce some categorical properties of Whittaker modules, for instance, we prove a mirabolic analog of Kostant's theorem. We also show that it quantizes a rational version of the Cremmer--Gervais $r$-matrix. As application, we construct a universal vertex-IRF transformation from the standard dynamical $R$-matrix to this constant one in categorical terms.}
	

\Classification{17B37; 17B38}

\section{Introduction}
	
Solutions to the \emph{quantum Yang--Baxter equation} (QYBE for short)
$R_{12} R_{13} R_{23} = R_{23} R_{13} R_{12}$,
	the so-called \emph{quantum $R$-matrices}, have many applications in mathematical physics and theory of quantum groups. The most well-studied example is the \emph{standard solution}, associated to any reductive Lie algebra $\g$, that leads to Drinfeld--Jimbo quantum groups \cite{DrinfeldICM,Jimbo}.
	
The starting point for the paper is a particular non-standard solution in type A which is called the \emph{Cremmer--Gervais $R$-matrix}. Originally, it appeared in \cite{CremmerGervais} and was later studied, for instance, in \cite{BuffenoirRocheTerras,Hodges}. This solution is special among other non-standard solutions in at least two ways (which might be connected to each other).
	
First, it is related to the so-called \emph{vertex-IRF transformation}. Namely, there is another important equation in mathematical physics, which is a generalization of the QYBE, called \emph{quantum dynamical Yang--Baxter equation}
	\[
	R_{12}(\lambda) R_{13}(\lambda+h_2) R_{23} (\lambda)= R_{23}(\lambda + h_1) R_{13}(\lambda) R_{12}(\lambda + h_3),
	\]
	where $\lambda$ is a parameter on the dual space $\h^*$ to an abelian Lie algebra $\h$ and $h_i$ are weight operators, see Definition~\ref{def::all_dynamical}. As explained by Felder \cite{Felder}, it is closely related to the star-triangle relation for face-type statistical models \cite{Baxter}. Moreover, it naturally appears in the description of the Liouville and Toda field theories \cite{GervaisNeveu}. As in the non-dynamical case, it leads to the theory of \emph{dynamical quantum groups}, and for any reductive $\g$, there are canonically defined standard solutions, see \cite{EtingofICM,EtingofSchiffmann}. It turns out that in type A, it is possible to gauge the standard dynamical $R$-matrix to a constant one which is exactly the Cremmer--Gervais $R$-matrix (in fact, this is how it was discovered in \cite{CremmerGervais}). What is even more interesting is that such kind of gauge transformations (which are called \emph{vertex-IRF} transformations) is unique, that is, it exists only in type A and the only constant $R$-matrix we can get is the Cremmer--Gervais one, see \cite{BalogDabrowskiFeher}. So, it is natural to ask: what is a representation-theoretic explanation of this fact? Or, at least, what is a representation-theoretic meaning of \emph{the} vertex-IRF transformation?
	
	This brings to the second special feature of the Cremmer--Gervais $R$-matrix: its relation to W-algebras. In principle, the exchange algebra of the Toda field theory (which is related to affine W-algebras) is ``controlled'' by the standard dynamical $R$-matrix \cite{GervaisNeveu}, but, as we mentioned before, it was shown in \cite{CremmerGervais} that is possible to gauge away the dependence on the dynamical parameter. Likewise, it was observed in \cite{BuffenoirRocheTerrasAffine} that the Cremmer--Gervais $R$-matrix can be constructed using the \emph{Sevostyanov characters} \cite{SevostyanovWhittakerModules} which are used to define a $q$-analog of the principal finite W-algebra. So, one may ask: what is a direct relation between the Cremmer--Gervais $R$-matrix and W-algebras? This paper is an attempt to answer these questions in a somewhat simpler situation of the rational analog of the Cremmer--Gervais solution.
	
	Let us recall Konstant's construction of the finite principal W-algebra \cite{Kostant}. Let $G$ be a simple Lie group with the Lie algebra $\g$. Choose a character $\psi \in \n_-^*$ of the negative nilpotent subalgebra $ \n_- \subset \g $ taking nonzero values on simple root generators and denote by~\smash{$\n_-^{\psi} = \{x - \psi(x) \mid x\in \n_-\}$} the shift of $\n_-$ along $\psi$ in the universal enveloping algebra $\mathrm{U}(\g)$. Its quotient by the left action~\smash{$\n_-^{\psi} \backslash \mathrm{U}(\g)$} has a right $\mathrm{U}(\g)$-action; let~\smash{$\bigl(\n_-^{\psi} \backslash \mathrm{U}(\g)\bigr)^{N_-}$} be the space of \emph{Whittaker vectors}, i.e., the elements on which $\n_-$ acts with character $\psi$. It is an example of a~quantum Hamiltonian reduction, and as such, it has a natural product. By the result of loc.\ cit., this algebra is isomorphic to the center $\Z_{\g}$ of $\mathrm{U}(\g)$.
	
	There is a categorical analog of this construction involving the category $\mathrm{HC}(G)$ of \emph{Harish-Chandra bimodules}. Those are $\mathrm{U}(\g)$-bimodules with certain $G$-integrability condition, in particular, $\mathrm{HC}(G)$ is a monoidal category. A similar assignment \smash{$X \mapsto \bigl(\n_-^{\psi} \backslash X\bigr)^{N_-}$} for $X\in \mathrm{HC}(G)$ can be enhanced to a functor to the category of $\Z_{\g}$-bimodules that has a natural ``multiplicative'' (lax monoidal) structure analogous to the standard product construction on quantum Hamiltonian reductions
	\[
	\bigl(\n_-^{\psi} \backslash X\bigr)^{N_-} \otimes_{\Z_{g}} \bigl(\n_-^{\psi} \backslash Y\bigr)^{N_-} \rightarrow \bigl(\n_-^{\psi} \backslash X \otimes_{\mathrm{U}(\g)} Y\bigr)^{N_-}, \qquad [x] \otimes [y] \mapsto [x\otimes y].
	\]
	This is the functor $\res_{\g}^{\psi}$ of the \emph{Kostant--Whittaker reduction} studied in \cite{BezrukavnikovFinkelberg}. By the results of loc.\ cit., the lax monoidal maps above are in fact isomorphisms.
	
	There is a version for the trivial character where one can consider the space $(X/\n)^{N}$ of the invariant vectors under the left action of $\n$ in the quotient $X/\n$. In particular, for $X = \mathrm{U}(\g)$, the result is isomorphic to the universal enveloping algebra $\mathrm{U}(\h)$ of the Cartan subalgebra~$\h$ of the maximal torus $H\subset G$. According to \cite{KalmykovSafronov}, this assignment can be enhanced to a~functor~${\res_{\g} \colon \mathrm{HC}(G) \rightarrow \mathrm{HC}(H)}$ of \emph{parabolic restriction} which, upon extending the scalars to the fraction field of $\mathrm{U}(\h)$ (that we ignore in the introduction), also becomes monoidal.
	
	In fact, there is more to the latter construction. For any $G$-representation $V$, the free right $\mathrm{U}(\g)$-module $V \otimes \mathrm{U}(\g)$ can be endowed with a left action making it an object of $\mathrm{HC}(G)$. This assignment can be upgraded to a monoidal functor $\free \colon \Rep(G) \rightarrow \mathrm{HC}(G)$ of \emph{free} Harish-Chandra bimodules. According to the classical results of \cite{ZhelobenkoSalgebras} and \cite{EtingofVarchenkoExchange}, the parabolic restriction functor ``trivializes'' on such bimodules, in other words, there is a natural isomorphism~${\res_{\g}(\free(V)) \cong V \otimes \mathrm{U}(\h)}$ after extending the scalars. By loc.\ cit., the image of the braiding on $\Rep(G)$ under the monoidal isomorphisms of $\res_{\g}$ gives a solution to the QDYBE.
	
	The goal of this paper is to adapt this construction to the Whittaker setting for $\g=\gl_N$. The only missing ingredient is a trivialization isomorphism \smash{$\res_{\gl_N}^{\psi}(\free(V)) \cong V \otimes \Z_{\gl_N}$}. To define it, we follow the approach of \cite{ZhelobenkoSalgebras} that uses the so-called \emph{extremal projector} in the setting of parabolic restriction introduced in \cite{AsherovaSmirnovTolstoy}. It is an element $P$ lying in a certain completion of~$\mathrm{U}(\g)$ localized by rational functions on $\h^*$; in particular, its action is well defined only for those $\n$-integrable $\mathrm{U}(\g)$-modules whose weights are generic. For such modules, it projects to $\n$-invariants along the image of the action of $\n_-$. In particular, an isomorphism $V \otimes \mathrm{U}(\h) \rightarrow \res_{\g}(\free(V))$ can be constructed by applying $P$ to $v\otimes 1$ in $V \otimes \mathrm{U}(\g)/\n$ (again, upon extending the scalars).\looseness=1
	
	The main result of the paper is the construction of an analogous element $P_{\m_N}(\vec{u})$ for $\g=\gl_N$ depending on a vector of parameters $\vec{u}=(u_1,\dots,u_{N-1})$ that we call the \emph{Kirillov projector} as it bears some resemblance to the Kirillov model from the theory of $p$-adic groups \cite{Kirillov}, see Definition~\ref{def::kirillov_projector}. For instance, it is defined not just for $\gl_N$, but for the so-called \emph{mirabolic} subalgebra~$\m_N$ of the mirabolic subgroup $\mathrm{M}_N \subset \GL_N$ of matrices preserving the last basis vector in $\C^N$. Theorem~\ref{thm::main_thm} is the main theorem of the paper, here $E_{ij} \in \gl_N$ are the standard matrix unit generators.
	
	\begin{thmx}
		\label{intro_thm::main_thm}
		For any right Whittaker module over $\m_N$, the Kirillov projector defines a unique linear operator satisfying
		\begin{gather*}
			P_{\m_N}(\vec{u})^2 = P_{\m_N}(\vec{u}), \qquad
			P_{\m_N} (\vec{u})\cdot (E_{ij} - \delta_{i,j+1}) = 0,\qquad 1 \leq j < i \leq N, \\
			(E_{ij} + \delta_{ij} u_j)\cdot P_{\m_N} (\vec{u}) = 0,\qquad 1 \leq i \leq j \leq N-1,
		\end{gather*}
		and acting by identity on the space of Whittaker vectors.
	\end{thmx}
	
	Let us indicate main features of the construction. First, contrary to the case of the extremal projector whose action requires certain genericity assumption on the weights, the Kirillov projector does not involve any denominators, thus has a well-defined action on \emph{any} right Whittaker module. Second, it is constructed inductively using a sequence of mirabolic subalgebras $\m_1 \subset \m_2 \subset \dots \subset \m_N$ embedded via lower right corner, at each step projecting to the space of Whittaker vectors for the corresponding $\m_i$ (and this was another motivation for the name). Third, its construction involves the evaluation homomorphism to $\mathrm{U}(\gl_N)$ from the \emph{Yangian}~$\mathrm{Y}(\gl_N)$. More precisely, at each step, the Kirillov projector is defined using certain \emph{quantum minors}, and the proof of Theorem~\ref{intro_thm::main_thm} uses extensively their properties expounded in~\cite{MolevYangians}. Unfortunately, as we mentioned, the evaluation homomorphism is crucial for the proofs, and we do not know if there is a purely Yangian analog of the Kirillov projector.
	
	Recall that by Kostant's result \cite{Kostant}, the category $\mathrm{Wh}(\g)$ of \emph{Whittaker modules}, i.e., the $\mathrm{U}(\g)$-modules on which $\n_-^{\psi}$ acts locally nilpotently, is equivalent to the category of $\Z_{\g}$-modules. Moreover, the following \emph{translation} property holds: for any $G$-representation $V$, there is a~canonical isomorphism
	\[
	\bigl(\n_-^{\psi} \backslash V \otimes \mathrm{U}(\g)\bigr)^{N_-} \cong V \otimes \bigl(\n_-^{\psi} \backslash \mathrm{U}(\g)\bigr)^{N_-} = V \otimes \Z_{g}.
	\]
	With the help of the Kirillov projector, we show a mirabolic analog of this result for the category~$\mathrm{Wh}(\m_N)$ of (right) Whittaker modules over the mirabolic subalgebra, see Theorem~\ref{thm::mirabolic_skryabin}.
	
	\begin{thmx}
		The functor $(-)^{N_-} \colon \mathrm{Wh}(\m_N) \rightarrow \Vect$ of $\n_-^{\psi}$-invariants is an equivalence of categories such that for any representation $V$ of the mirabolic subgroup $\mathrm{M}_N$, there is a canonical isomorphism
	\smash{$
		\bigl(\n_-^{\psi} \backslash V \otimes \mathrm{U}(\m_N)\bigr)^{N_-} \cong V$}.
	\end{thmx}
	
	In particular, the natural tensor structure on the mirabolic analog $\res_{\m_N}^{\psi}$ of the Kostant--Whittaker reduction discussed above induces a tensor structure on the composition
	\[
	\Rep(\mathrm{M}_N) \xrightarrow{\free} \mathrm{HC}(\mathrm{M}_N) \xrightarrow{\res_{\m_N}^{\psi}} \Vect,
	\]
	which, in particular, produces a solution to the QYBE. We show that, up to a small difference between the mirabolic subgroup in $\GL_N$ and the maximal parabolic subgroup of matrices in $\mathrm{SL}_N$ preserving the linear span of the last basis vector, this construction produces a deformation of the rational analog of the Cremmer--Gervais solution as in \cite{EndelmanHodges}. Extending this construction to~$\GL_N$, the other main result of the paper can be summarized as follows, see Theorem~\ref{thm::kw_tensor_structure} and Section~\ref{subsect:vector_rep_constant}.
	
	\begin{thmx}
		For every vector of parameter $\vec{u}$, upon trivialization induced by the Kirillov projector $P_{\m_N}(\vec{u})$, the monoidal structure on the functor $\Rep(\GL_N) \rightarrow \BiMod{\Z_{\gl_N}}{\Z_{\gl_N}}$ of the Kostant--Whittaker reduction restricted to free Harish-Chandra bimodules produces a constant solution~$R^{\CG}(\vec{u})$ to the QYBE such that for $\vec{u} = (N-1,\dots,1)$, it coincides with the rational analog of the Cremmer--Gervais $R$-matrix.
	\end{thmx}
	
	More precisely, for $\vec{u} = (N-1,\dots,1)$, we show that its semiclassical limit is the inverse of trace pairing with the \emph{principal nilpotent element}
	\[
	\m_N \wedge \m_N \rightarrow \C, \qquad x \wedge y \mapsto \Tr(e\cdot[x,y]), \qquad e := E_{12} + \dots + E_{N-1,N},
	\]
	that coincides (up to the aforementioned difference) with the degeneration of the classical Cremmer--Gervais $r$-matrix by the results of \cite{GerstenhaberGiaquinto}. Moreover, in Section~\ref{subsect:vector_rep_constant} we compute explicitly the quantum $R$-matrix for the dual of the vector representation of $\GL_N$ and show that it also coincides with the rational Cremmer--Gervais solution, now on the quantum level \cite{EndelmanHodges}. Surprisingly, this calculation is essentially based on the action of the \emph{degenerate affine Hecke algebra} on certain Harish-Chandra bimodules, for instance, the resulting $R$-matrix can be simply interpreted as the action of the permutation operator on its polynomial representation.
	
	For an arbitrary vector of parameters, this construction produces a family of constant solutions to QYBE that we call $h$-deformed rational Cremmer--Gervais $R$-matrices, see Definition~\ref{def::cremmer_gervais}. Its semiclassical limit is obtained by replacing $e$ with $e-h$, where $h = u_1 E_{11} + \dots + u_{N-1} E_{N-1,N-1}$, in the definition above.
	
	Finally, as application, we construct a categorical version of the vertex-IRF transformation. It turns out that it can be interpreted via a finite analog of the \emph{Miura transform} \cite{GinzburgKazhdan}: using the Harish-Chandra homomorphism $\Z_{\g} \rightarrow \U(\h)$, one can construct the maps
	\[
	\bigl(\n_-^{\psi} \backslash X\bigr)^{N_-} \otimes_{\Z_{\g}} \mathrm{U}(\h) \rightarrow \n_-^{\psi} \backslash X /\n \leftarrow (X/\n)^{N},
	\]
	that become isomorphisms after extending the scalars to the rational functions on $\h^*$. Moreover, the resulting map $\res^{\psi}_{\g} (X) \rightarrow \res_{\g}(X)$ is compatible with the tensor structure. In particular, for~${G = \GL_N}$, we have the following result, see Theorems~\ref{thm::the_vertex-irf} and~\ref{thm:the_vertex_irf_transformation}.
	
	\begin{thmx}
		The natural transformation between the functor of the Kostant--Whittaker reduction and the parabolic restriction induces a gauge transformation between the standard solution of the QDYBE and the $h$-deformed rational Cremmer--Gervais solution.
	\end{thmx}
	
	We compute this gauge transformation for the special vector of parameters $\vec{u} = (N-1,\dots,1)$ and for the dual of the vector representation. In accordance with the
	earlier constructions~\cite{BalogDabrowskiFeher,BuffenoirRocheTerras,CremmerGervais}, it turns out that, up to a diagonal factor, it is given by the Vandermonde matrix.
	
	We expect that there should be a $q$-deformation of the constructions of the paper that would produce the actual Cremmer--Gervais solution. Indeed, there is a notion of quantum Harish-Chandra bimodules (for instance, see \cite[Section 2.5]{KalmykovSafronov} and references therein) that involves the quantum group $\mathrm{U}_q(\g)$ as well as a generalization of the parabolic restriction. We expect that, similarly to Sevostyanov's construction \cite{SevostyanovWhittakerModules} of the center of $\mathrm{U}_q(\g)$ via a $q$-analog of the Whittaker reduction, there should be a $q$-analog of the Kostant--Whittaker reduction functor as well as a~$q$-analog of the Kirillov projector that uses the evaluation homomorphism from the quantum affine algebra $\mathrm{U}_q\bigl(\hat{\gl}_N\bigr)$. Then the tensor structure on the former under the trivialization of the latter should produce the original Cremmer--Gervais solution.
	
	The properties of the Kirillov projector were predicted using SageMath \cite{SageMath}. The author's code that implements the action of the projector on certain Harish-Chandra bimodules is available here: \url{https://github.com/art-kalm/whittaker}.

{\bf Structure of the paper.} In Section~\ref{sect::dynamical_quantum_groups}, we give a general setup for the paper and review the categorical approach to (standard) dynamical quantum groups via the category of Harish-Chandra bimodules and parabolic restriction functor as in \cite{KalmykovSafronov}. Proofs are mostly omitted unless there is a version that we need later in the paper and it is different from loc.\ cit. In Section~\ref{sect::mirabolic}, we define the mirabolic subgroup and a family of $r$-matrices that deform the degenerate version of the Cremmer--Gervais $r$-matrix. In Section~\ref{sect::basics}, we review Yangians and related constructions such as quantum minors. We prove various technical results that we use later in the paper. In Section~\ref{sect::kirillov_projector}, we introduce the Kirillov projector and prove its defining properties. In Section~\ref{sect::kw_red}, we recall the Kostant--Whittaker reduction functor, introduce its mirabolic version, and show some properties thereof. We also obtain quantization of the family of rational Cremmer--Gervais solutions. In Section~\ref{sect::vertex_irf}, we define a categorical version of the vertex-IRF transformation and show that an equivalence between parabolic restriction and Kostant--Whittaker reduction gives the vertex-IRF transformation in the classical sense.
	
	All the objects and constructions in the paper involve the asymptotic parameter $\hbar$; to simplify the exposition, we have set $\hbar=1$ in the introduction.

	\section{Background: dynamical quantum groups}
	\label{sect::dynamical_quantum_groups}
	
	\renewcommand{\theequation}{\arabic{section}.\arabic{equation}}
	
	In this section, we will give a brief overview of the classical extremal projector and its relation to dynamical quantum groups through the category of Harish-Chandra bimodules.
	
	\subsection{General setup}
	\label{subsect::general_setup}
	For the rest of the paper, we work over an algebraically closed field $\C$ of characteristic zero. All categories and functors we will consider are $\C$-linear; likewise, unless otherwise stated, all tensor products are taken over $\C$. Often, we will use the algebra of polynomial functions $\C[\hbar]$ in variable $\hbar$ and denote by $\Mod{\C[\hbar]}$ the corresponding monoidal category of $\C[\hbar]$-modules with the standard tensor product $\otimes_{\C[\hbar]}$ over $\C[\hbar]$. For any vector space $V$ over $\C$, we denote by~${V[\hbar]:= V \otimes_{\C} \C[\hbar]}$.
	
Throughout this paper, we work with locally presentable categories (we refer to \cite{AdamekRosicky} and \cite[Section 2]{BCJF} for more details). We will also use Deligne's tensor product of locally presentable categories introduced in \cite{Deligne}.

	\subsection{Harish-Chandra bimodules}
	\label{subsect::hc_bimod}
	For details and proofs in a more general setting, we refer the reader to \cite[Section 2]{KalmykovSafronov}.
	
	Let $G$ be an affine algebraic group and $\g$ be its Lie algebra with bracket $[\,,\,]$.
	
	\begin{Definition}
		\label{def::asympt_uea}
		The \emph{asymptotic universal enveloping algebra} $\U(\g)$ of $\g$ is a tensor algebra over~$\C[\hbar]$ generated by the vector space $\g$ with the relations
$
		xy - yx = \hbar [x,y]$, $ x,y\in \g$.
		The \emph{non-asymptotic universal enveloping algebra} $\mathrm{U}(\g)$ is a $\C$-algebra generated by $\g$ with the relations~${xy - yx = [x,y]}$.
	\end{Definition}
	
	\begin{Remark}
		In principle, most of the cited literature deals with the non-asymptotic version~$\mathrm{U}(\g)$. However, unless specified otherwise, all the constructions and proofs can be translated \emph{mutatis mutandis} to the asymptotic setting.
	\end{Remark}
	
	Denote by $\Rep(G)$ the braided monoidal category of $G$-representations over $\C$, where the braiding is given by the permutation map
	\[
	\sigma_{VW} \colon\ V \otimes W \rightarrow W \otimes V, \qquad v\otimes w \mapsto w \otimes v.
	\]
	For any $V \in \Rep(G)$ and $\xi\in \g$, we denote by $\ad_{\xi} \colon V \rightarrow V$ the differential of the $G$-action along $\xi$. The adjoint action of $G$ on $\g$ extends to an action on $\U(\g)$, in particular, the latter is naturally an object of $\Rep(G)$ such that the differential is
	\begin{equation*}
		\ad_{\xi}(x) = \frac{\xi x - x \xi}{\hbar},\qquad \xi\in \g, x\in \U(\g).
	\end{equation*}
	
	\begin{Definition}[{\cite{BezrukavnikovFinkelberg}}]
		An \emph{asymptotic Harish-Chandra bimodule} is a left $\U(\g)$-module $X$ in the category $\Rep(G)$. In other words, it has a structure of a $G$-representation and a left $\U(\g)$-module such that the action morphism
		$
		\U(\g) \otimes X \rightarrow X
		$
		is a homomorphism of $G$-representations. The category of Harish-Chandra bimodules is denoted by $\HC(G)$.
	\end{Definition}
	
	In what follows, we will refer to asymptotic Harish-Chandra bimodules as simply Harish-Chandra bimodules.
	
	There is a natural right $\U(\g)$-module structure on any Harish-Chandra bimodule $X$. Namely, for $\xi \in\g$, define
	\begin{equation}
		\label{eq::left_action}
		x\xi := \xi x - \hbar \ad_{\xi}(x),\qquad x\in X,
	\end{equation}
	and extend it to a right $\U(\g)$-action. Therefore, the category $\HC(G)$ is a subcategory of $\U(\g)$-bimodules, hence is equipped with a tensor structure
$
	X \otimes^{\HC(G)} Y := X \otimes_{\U(\g)} Y$.
	
	\begin{Remark}
		The subcategory of $\HC(G)$ where $\hbar$ acts by identity is identified with the standard category of Harish-Chandra bimodules \cite[Definition 5.2]{BernsteinGelfand}. The one where $\hbar$ acts by zero can be identified with the category of $G$-equivariant quasi-coherent sheaves on $\g^*$.
	\end{Remark}
	
	There is a natural functor of the so-called \emph{free $($right$)$ Harish-Chandra bimodules}
	\[
	\free \colon \Rep(G) \rightarrow \HC(G), \qquad V \mapsto V \otimes \U(\g),
	\]
	such that the \emph{right} $\U(\g)$-action is free and the $G$-action is diagonal. By \cite[equation~(28)]{KalmykovSafronov}, it is isomorphic to the functor of free \emph{left} Harish-Chandra bimodules $V \mapsto \U(\g) \otimes V$; however, in what follows, we will use the right version. One can check that this functor is monoidal. In fact, all Harish-Chandra bimodules can be ``constructed'' from the free ones.
	\begin{Proposition}[{\cite[Proposition 2.7]{KalmykovSafronov}}]
		\label{prop::hc_generators}
		The category $\HC(G)$ is generated by $\free(V)$ for $V \in \Rep(G)$.
	\end{Proposition}
	
	\subsection[Constant and dynamical R-matrices]{Constant and dynamical $\boldsymbol{R}$-matrices}
	\label{subsect::dynamica_r_mat}
	For introduction to the theory of dynamical quantum groups, see \cite{EtingofSchiffmann}; for the interpretation in terms of the category of Harish-Chandra bimodules, see \cite{KalmykovSafronov}.
	
	\begin{Definition}
		\label{def::all_constant}\qquad
		\begin{itemize}\itemsep=0pt
			\item Let $F_{VW} \in \End_{\C[\hbar]}(V[\hbar]\otimes_{\C[\hbar]} W[\hbar]) $ be a collection of invertible $\C[\hbar]$-linear maps natural in $V,W\in \Rep(G)$. The \emph{twist} equation is
			\begin{equation}
				\label{eq::twist}
				F_{U \otimes V,W} \circ (F_{UV} \otimes \id_W) = F_{U,V\otimes W} \circ (\id_U \otimes F_{VW}),
			\end{equation}
			where both sides of the equation lie in $\End_K(U[\hbar] \otimes_{\C[\hbar]} V[\hbar] \otimes_{\C[\hbar]} W[\hbar])$. A solution is called a \emph{$($Drinfeld$)$ twist}.
			\item Let $R_{VW}\in \End_{\C[\hbar]}(V[\hbar]\otimes_{\C[\hbar]} W[\hbar]) $ be a collection of invertible $\C[\hbar]$-linear maps natural in $V,W\in \Rep(G)$. The \emph{quantum Yang--Baxter equation} is
			\begin{equation}
				\label{eq::qybe}
				R_{UV} R_{UW} R_{VW} = R_{VW} R_{UW} R_{UV},
			\end{equation}
			where both sides of the equation lie in $\End_K(U[\hbar] \otimes_{\C[\hbar]} V[\hbar] \otimes_{\C[\hbar]} W[\hbar])$. A solution is called a \emph{quantum $R$-matrix}.
			\item Let $r_{VW}\in \End_{\C}(V\otimes W) $ be a collection of linear maps natural in $V,W\in \Rep(G)$. The \emph{classical Yang--Baxter equation} is
			\begin{equation}
				\label{eq::cybe}
				[r_{UV},r_{UW}] + [r_{UV},r_{VW}] + [r_{UW}, r_{VW}] = 0,
			\end{equation}
			where the left-hand side lies in $\End_{\C}(U \otimes V \otimes W)$. A solution is called a \emph{classical $r$-matrix}.
		\end{itemize}
	\end{Definition}
	
	The following result is standard and follows directly from the corresponding definitions.
	\begin{Proposition}\qquad
		\begin{enumerate}\itemsep=0pt
			\item[$(1)$] A tensor structure on the functor
$
			\Rep(G) \rightarrow \Mod{\C[\hbar]}$, $ V \mapsto V[\hbar]$,
			is equivalent to the data of a Drinfeld twist \eqref{eq::twist} {\rm\cite[\emph{Definition} 2.4.1]{EGNO}}.
			\item[$(2)$] Let $R_{VW} = \sigma_{VW}^{-1} F_{WV}^{-1} \sigma_{VW} F_{VW}$, where $\sigma_{VW} \colon V \otimes W \rightarrow W \otimes V$ is the permutation map. Then $R_{VW}$ is a quantum $R$-matrix \eqref{eq::qybe} {\rm\cite[\emph{Definition} 8.1.1]{EGNO}}.
			\item[$(3)$] Assume that $R_{VW}$ has the form $R_{VW} = \id_{V \otimes W} + \hbar r_{VW} + O\bigl(\hbar^2\bigr)$. Then $r_{VW}$ is a classical $r$-matrix \eqref{eq::cybe} {\rm\cite[\emph{Proposition} 9.4]{EtingofSchiffmannQuantumGroups}}.
		\end{enumerate}
	\end{Proposition}
	
	There is a version of these equations involving dynamical parameter $\lambda \in \h^*$. Let $H \subset G$ be a torus and $\h \subset \g$ be its Lie algebra. Let $V_1,\dots,V_n$ be a collection of $H$-representations, and denote by $\End_{H}(V_1 \otimes \dots \otimes V_n)$ the space of endomorphisms commuting with the induced diagonal $H$-action. For any $\End_{H}(V_1 \otimes \dots \otimes V_n)$-valued rational function $f(\lambda,\hbar)$ on $\h^*\times \mathbb{A}^1$, we will use the following notation:
	\[
	f(\lambda-\hbar h_{V_i},\hbar) (v_1 \otimes \dots \otimes v_n) := f(\lambda - \hbar \mu_i,\hbar) (v_1\otimes \dots \otimes v_n),
	\]
	if vector $v_i$ is of weight $\mu_i \in \h^*$. In what follows, we will not indicate dependence on $\hbar$, for instance, we denote $f(\lambda):=f(\lambda,\hbar)$.
	
	\begin{Definition}\qquad
		\label{def::all_dynamical}
		\begin{itemize}\itemsep=0pt
			\item Let $J_{VW}(\lambda)$ be a collection of invertible $\End_{H}(V \otimes W)$-valued functions on $\h^*\times \mathbb{A}^1$ natural in $V,W\in \Rep(G)$. The \emph{dynamical twist} equation is
			\begin{equation}
				\label{eq::dyn_twist}
				J_{U \otimes V,W}(\lambda) \circ (J_{UV}(\lambda+\hbar h_W) \otimes \id_W) = J_{U,V\otimes W} \circ (\id_U \otimes J_{VW}(\lambda)).
			\end{equation}
			A solution is called a \emph{dynamical $($Drinfeld$)$ twist}.
			\item Let $R_{VW}(\lambda)$ be a collection of invertible $\End_{H}(V \otimes W)$-valued functions on $\h^*\times \mathbb{A}^1$ natural in $V,W\in \Rep(G)$. The \emph{quantum dynamical Yang--Baxter equation} is
			\begin{equation}
				\label{eq::qdybe}
				R_{UV}(\lambda) R_{UW} (\lambda + \hbar h_V) R_{VW}(\lambda) = R_{VW}(\lambda+\hbar h_U) R_{UW}(\lambda) R_{UV}(\lambda + \hbar h_W).
			\end{equation}
			A solution is called a \emph{quantum dynamical $R$-matrix}.
			\item Let $r_{VW}(\lambda)$ be a collection of $\End_{H}(V \otimes W)$-valued functions on $\h^*$ natural in $V,W\in \Rep(G)$. The \emph{classical dynamical Yang--Baxter equation} is
			\begin{gather}
					-\sum\limits_i \left( (x_i)_{U} \frac{\partial r_{VW}}{\partial x^i} - (x_i)_V \frac{\partial r_{UW}}{\partial x^i} + (x_i)_W \frac{\partial r_{UV}}{\partial x^i}\right) +\nonumber \\
					\qquad{}+[r_{UV}(\lambda),r_{UW}(\lambda)] + [r_{UV}(\lambda),r_{VW}(\lambda)] + [r_{UW}(\lambda), r_{VW}(\lambda)] = 0,\label{eq::cdybe}
			\end{gather}
			where $\{x_i\} \subset \h$ is a basis in $\h$, and $\bigl\{x^i\bigr\}$ is the dual basis in $\h^*$. A solution is called a~\emph{classical dynamical $r$-matrix}.
		\end{itemize}
	\end{Definition}
	
	\begin{Remark}
		In particular, a quantum $R$-matrix from Definition~\ref{def::all_constant} (respectively, Drinfeld twist, classical $r$-matrix) that commutes with the $H$-action defines a quantum dynamical $R$-matrix from Definition~\ref{def::all_dynamical} (respectively, dynamical Drinfeld twist, classical dynamical $r$-matrix) which does not depend on the dynamical parameter.
	\end{Remark}
	
	It turns out that there is a Tannakian interpretation of dynamical $R$-matrices in terms of functors to the category of Harish-Chandra bimodules. Namely, consider the category $\HC(H)$ of Harish-Chandra bimodules over $H$. Since $\h$ is commutative, the universal enveloping algebra~$\U(\h)$ is equal to the algebra of polynomial functions $\C\bigl[\h^* \times \mathbb{A}^1\bigr]$ with a trivial $H$-action, where~$\mathbb{A}^1$ corresponds to $\hbar$. Representations of $H$ are completely reducible with irreducibles being one-dimensional parameterized by weights $\Lambda(H)$; in particular, as a plain category,
	\[
	\HC(H) \cong \bigoplus\limits_{\mu \in \Lambda(H)} \RMod_{\C[\h^* \times \mathbb{A}^1]}[\mu],
	\]
	where \smash{$\RMod_{\C[\h \times \mathbb{A}^1]}$} is the category of right $\C\bigl[\h^* \times \mathbb{A}^1\bigr]$-modules and $[\mu]$ simply indicates the corresponding graded component. For $X\in \RMod_{\C[\h^* \times \mathbb{A}^1]}[\mu]$, the left $\C\bigl[\h^* \times \mathbb{A}^1\bigr]$-action is
	\begin{equation}
		\label{eq::hc_left_action}
		f(\lambda) \otimes x \mapsto x f(\lambda+\hbar \mu).
	\end{equation}
	For every graded component, the tensor product is
	\begin{align}
			\otimes^{\HC(H)}\colon\ \RMod_{\C[\h^* \times \mathbb{A}^1]} [\mu] \otimes \RMod_{\C[\h^* \times \mathbb{A}^1]} [\nu] & \rightarrow \RMod_{\C[\h^* \times \mathbb{A}^1]} [\mu+\nu],\nonumber \\
(X[\mu],Y[\nu]) & \mapsto (X \otimes_{\C[\h^* \times \mathbb{A}^1]} Y)[\mu+\nu].\label{eq::hc_h_tensor_product}
	\end{align}
	
	In what follows, we will need a \emph{generic} version of this category. Denote by
	\begin{equation}
		\label{eq::generic_weights}
		\bigl(\h^* \times \mathbb{A}^1\bigr)^{\gen} = \bigl\{(\lambda,\hbar) \in \h^* \times \mathbb{A}^1 \mid \bigl\langle \lambda, \alpha^{\vee} \bigr\rangle \not\in \hbar \mathbb{Z},\, \alpha \in \Delta_+\bigr\}
	\end{equation}
	the set of generic weights, where $\Delta_+$ is the set of positive roots, and let $\U(\h)^{\gen} \subset \Frac\bigl(\cO\bigl(\h^* \times \mathbb{A}^1\bigr)\bigr)$ be the algebra of rational functions regular on $\bigl(\h^* \times \mathbb{A}^1\bigr)^{\gen}$. Define
	\begin{equation}
		\label{eq::generic_hc_definition}
		\HC(H)^{\gen} := \bigoplus\limits_{\mu\in \Lambda(H)} \RMod_{\U(\h)^{\gen}} [\mu]
	\end{equation}
	with the left action and the tensor product given by obvious generalization of formulas \eqref{eq::hc_left_action} and~\eqref{eq::hc_h_tensor_product}, respectively. Similarly to non-generic version, there is a functor of free Harish-Chandra bimodules
$
	\free\colon \Rep(H) \rightarrow \HC(H)^{\gen}$, $ V \mapsto V \otimes \U(\h)^{\gen}$,
	such that if $V = \bigoplus_{\mu\in \Lambda} V[\mu]$ is the weight decomposition, then each summand $V[\mu] \otimes \U(\h)^{\gen}$ is a free right $\U(\h)^{\gen}$-module lying in $\RMod_{\U(\h)^{\gen}}[\mu]$. By direct sum decomposition \eqref{eq::generic_hc_definition}, one can see that for any $V\in \Rep(H)$, endomorphisms $\End_{\HC(H)^{\gen}} (\free(V))$ can be identified with $\End_H(V)$-valued rational functions on $\h^* \times \mathbb{A}^1$ regular on generic weights $\bigl(\h^* \times \mathbb{A}^1\bigr)^{\gen}$. We use this identification for the rest of the paper.
	
	\begin{Theorem}[{\cite{KalmykovSafronov}}]
		\label{thm::dynamical_tannakian}
		Consider the composition
		\[
		\Rep(G) \xrightarrow{\mathrm{forget}} \Rep(H) \xrightarrow{\free} \HC(H)^{\gen},
		\]
		where the first arrow is the forgetful functor.
		\begin{enumerate}\itemsep=0pt
			\item[$(1)$] A monoidal structure on $\Rep(G) \rightarrow \HC(H)^{\gen}$ is equivalent to the data of a dynamical twist \eqref{eq::dyn_twist} regular on the generic weights.
			\item[$(2)$] Define $R_{VW}(\lambda)$ as the composition
			\[
			V\otimes W \xrightarrow{J_{VW}(\lambda)} V \otimes W \xrightarrow{\sigma_{VW}} W \otimes V \xrightarrow{J_{WV}^{-1}(\lambda)} W \otimes V \xrightarrow{\sigma_{VW}^{-1}} V \otimes W,
			\]
			where $\sigma_{VW} \colon V \otimes W \rightarrow W \otimes V$ is the permutation map. Then $R_{VW}(\lambda)$ is a quantum dynamical $R$-matrix \eqref{eq::qdybe}.
			\item[$(3)$] Assume that $R_{VW}(\lambda)$ admits Taylor series around $\hbar=0$ such that
			\[R_{VW} (\lambda)= \id_{V \otimes W} + \hbar r_{VW}(\lambda) + O\bigl(\hbar^2\bigr). \]
			Then $r_{VW}(\lambda)$ is a classical dynamical $r$-matrix \eqref{eq::cdybe}.
		\end{enumerate}
	\end{Theorem}
	
	\begin{Remark}
		Observe that the results in loc.\ cit.\  are formulated and proved for the \emph{left} Harish-Chandra bimodules, in particular, the weight shifts in Definition~\ref{def::all_dynamical} are different. However, they can easily adapted to the setting of \emph{right} modules that we use in the paper.
	\end{Remark}

	From now on, let $G$ be a reductive group. As in the non-dynamical case, there is a standard solution to classical dynamical Yang--Baxter equation. Namely, choose a Borel subgroup $B\subset G$ and a maximal torus $H \subset B$. It induces a triangular decomposition $\g = \n_- \oplus \h \oplus \n$ with respect to the set of positive roots $\Delta_+$. For $\alpha\in \Delta_+$, denote by $\{f_{\alpha}, h_{\alpha}, e_{\alpha}\}$ the corresponding $\sl_2$-triple of the Chevalley base such that $\langle \alpha,h_{\alpha} \rangle =2$.
	
	Define (see \cite{EtingofVarchenkoExchange})
	\begin{equation}
		\label{eq::standard_sol_cdybe}
		r(\lambda) = \sum_{\alpha \in \Delta_+} \frac{f_{\alpha} \wedge e_{\alpha}}{\langle\lambda,h_{\alpha} \rangle }.
	\end{equation}
	One can check that it indeed satisfies \eqref{eq::cdybe}. In what follows, we will recall its quantization in Section~\ref{subsect::parabolic_reduction}.
	
	\subsection{Extremal projector}
	\label{subsect::extremal_proj}
	
	The original construction can be found in \cite{AsherovaSmirnovTolstoy}, see also \cite{ZhelobenkoSalgebras}.
	
	Denote by $\U(\g)^{\gen} = \U (\h)^{\gen} \otimes_{\U(\h)} \U(\g)$ the localized version of the universal enveloping algebra. For $\alpha\in \Delta_+$, define the following element in a completion of $\U(\g)^{\gen}\bigl[\hbar^{-1}\bigr]$:
	\begin{equation}
		\label{eq::extremal_proj_term}
		P_{\alpha}(t) = \sum\limits_{k\geq 0} \frac{(-1)^k \hbar^{-k}}{k!} \frac{1}{\prod_{j=1}^k (h_{\alpha} + \hbar(t+j))} f_{\alpha}^k e_{\alpha}^k,
	\end{equation}
	Choose a normal ordering $<$ on $\Delta_+$ (i.e., such that $\alpha + \beta$ lies between $\alpha$ and $\beta$). Denote by $\rho$ the half-sum of positive roots. Define
	\begin{equation}
		\label{eq::extremal_proj}
		P := \prod_{\alpha \in \Delta_+}^{<} P_{\alpha}(t_{\alpha}),
	\end{equation}
	where the product is taken in the normal ordering, and $t_{\alpha} = h_{\alpha}(\rho)$. By the results of \cite{ZhelobenkoCocycles}, it lies in the space of formal power series of weight zero of the form
	\[
	\sum\limits_{\substack{k = (k_{\alpha} \mid \alpha\in \Delta_+) \\ l = (l_{\alpha} \mid \alpha\in \Delta_+)}} h_{kl} \cdot \prod_{\alpha \in \Delta_+} f_{\alpha}^{k_{\alpha}} \prod_{\alpha \in \Delta_+} e_{\alpha}^{l_{\alpha}},
	\]
	where each product is taken with respect to the normal ordering and $h_{kl} \in \U(\h)^{\gen}\bigl[\hbar^{-1}\bigr]$ are some elements.
	
	\begin{Theorem}[{\cite{AsherovaSmirnovTolstoy}}]
		\label{thm::extremal_proj}
		The element $P$ satisfies
$
		e_{\alpha} P = P f_{\alpha} = 0 $ $\forall \alpha \in \Delta_+$.
		The action of $P$ is well defined on left $\U(\g)^{\gen}$-modules with a locally nilpotent action of $\n$ satisfying $\n X \subset \hbar X$ and whose weight are generic.
	\end{Theorem}
	
	For any left $\U(\g)^{\gen}$-module $X$ such that $\n X \subset \hbar X$, observe that the action of $\hbar^{-1} \n$ is well-defined and closed under commutator, in particular, it defines an action of the Lie algebra $\n$. Assume that integrates to an action of the group $N$. Denote by $\n_- \backslash X$ the quotient by the left action of $\n_-$, and by $X^{N}$ the space of corresponding $N$-invariants. Then the extremal projector defines an inverse of the projection
	\begin{equation}
		\label{eq::extremal_proj_iso_inv_quot}
		X^{N} \rightarrow \n_- \backslash X
	\end{equation}
	for any left module with generic weights in the sense of \eqref{eq::generic_weights} and a locally nilpotent $\n$-action such that $\n X \subset \hbar X$.
	
	\subsection{Parabolic restriction functor}
	\label{subsect::parabolic_reduction}
	
	The main reference for the setup of this paper is \cite{KalmykovSafronov}, but the history of the subject is certainly much richer, for instance, see \cite{EtingofSchiffmann}.
	
	In what follows, the Borel subgroup $B$ acts trivially on $\U(\h)^{\gen}$. We consider any right $\U(\h)$-module as a $\U(\b)$-module via the quotient map $\b \rightarrow \h$ along $\n\subset \b$.
	\begin{Definition}
		A \emph{universal generic category $\cO$} is the category $\cO^{\univ,\gen}$ of $(\U(\g),\U(\h)^{\gen})$-bimodules in the category $\Rep(B)$ such that the differential of the $B$-action $\ad_{\xi}$ on an object $X$ along any $\xi\in \b$ satisfies
		\begin{equation}
			\label{eq::o_univ_diagonal_relation}
			\xi x - x\xi = \hbar \ad_{\xi}(x),\qquad x\in X.
		\end{equation}
		The \emph{universal generic Verma module} is
$
		M^{\univ} := \U(\g)^{\gen} \otimes_{\U(\n)} \C$.
	\end{Definition}
	
	We have two functors
	\begin{equation}
		\label{eq::act_h}
		\act_H \colon \ \HC(H)^{\gen} \rightarrow \cO^{\univ,\gen},\qquad X \mapsto M^{\univ} \otimes_{\U(\h)^{\gen}} X,
	\end{equation}
	and
	\begin{equation}
		\label{eq::act_g}
		\act_G \colon \ \HC(G) \rightarrow \cO^{\univ,\gen}, \qquad X \mapsto X \otimes_{\U(\g)} M^{\univ}.
	\end{equation}
	
	\begin{Proposition}
		For any $X\in \HC(G)$, the tensor product $X \otimes_{\U(\g)} M^{\univ}$ lies in $\cO^{\univ,\gen}$.
		\begin{proof}
			Denote by $X/\n:=X/X\n$ the quotient by the right action of $\n$. Being a quotient of a $B$-representation $X$ by a $B$-stable subset $X\n$, it is a $B$-representation as well, and \eqref{eq::o_univ_diagonal_relation} is satisfied for $X/\n$ due to the definition of the left action \eqref{eq::left_action}. In particular, the tensor product~${X/\n \otimes_{\U(\h)} \U(\h)^{\gen}}$ is an object of $\cO^{\univ,\gen}$. However, $X \otimes_{\U(\g)} M^{\univ}$ is naturally isomorphic to $X/\n \otimes_{\U(\h)} \U(\h)^{\gen} $.
		\end{proof}
	\end{Proposition}
	
	Let $X\in \cO^{\univ,\gen}$. Observe that due to \eqref{eq::o_univ_diagonal_relation}, the left action of $\n$ on $X$ satisfies $\n X \subset \hbar X$. Therefore, the action of the extremal projector \eqref{eq::extremal_proj} is well defined.
	
	\begin{Proposition}
		\label{prop::act_equiv}
		The functor $\act_H$ is an equivalence.
		\begin{proof}
			We repeat an argument from \cite[Theorem 4.17]{KalmykovSafronov}. By an obvious asymptotic analog of~\cite[Proposition 4.6]{KalmykovSafronov}, the functor $(-)^{N}\colon \cO^{\univ,\gen} \rightarrow \HC(H)^{\gen}$ is right adjoint to $\act_H$. The unit of the adjunction is given by
			\smash{$
			X \mapsto \bigl(M^{\univ} \otimes_{\U(\h)^{\gen}} X\bigr)^{N}$}.
			The extremal projector gives an isomorphism with the quotient functor as in \eqref{eq::extremal_proj_iso_inv_quot}
$
			{}_{\n_-}(-)\colon \cO^{\univ,\gen} \rightarrow \HC(H)$, $ M \mapsto \n_- \backslash M$.
			So, the functor $(-)^{N}$ is exact. We have
			\[
			\bigl(M^{\univ} \otimes_{\U(\h)^{\gen}} X\bigr)^{N} \cong \n_- \backslash \bigl(M^{\univ} \otimes_{\U(\h)^{\gen}} X\bigr),
			\]
			and the composition
			$
			X \mapsto \n_- \backslash \bigl(M^{\univ} \otimes_{\U(\h)^{\gen}} X\bigr)
			$
			is an isomorphism by the PBW theorem. Therefore, $\act_H$ is fully faithful.
			
			We have $M^{N} = 0$ if and only if $M=0$. But $(-)^{N}$ is exact, therefore, it is conservative. Since its left adjoint $\act_H$ is fully faithful, it is an equivalence.
		\end{proof}
	\end{Proposition}
	
	\begin{Definition}
		\label{def::parabolic_reduction}
		The \emph{parabolic restriction} functor $\res_{\g}\colon \HC(G) \rightarrow \HC(H)$ is the composition
		\[
		\HC(G) \xrightarrow{\act_G} \cO^{\univ,\gen} \xrightarrow{(-)^{N}} \HC(H)^{\gen}.
		\]
	\end{Definition}
	Explicitly, it is given by quantum Hamiltonian reduction and extension of scalars
	\[
	X \mapsto ((X/ \n) \otimes_{\U(\h)} \U(\h)^{\gen})^{N}.
	\]
	Equivalently, we can extend the scalars on the $\HC(G)$-side first and then take the quantum Hamiltonian reduction
$
	X \mapsto (X^{\gen}/ \n )^{N}$, $ X^{\gen} := \U(\g)^{\gen} \otimes_{\U(\g)} X$.
	There is a natural lax monoidal structure (in the sense of \cite[Section XI.2]{Maclane} where it is called simply ``monoidal structure'') on $\res_{\g}$
	\[
		(X^{\gen}/\n)^{N} \otimes_{\U(\h)^{\gen}} (Y^{\gen}/\n)^{N} \rightarrow ((X \otimes_{\U(\g)} Y)^{\gen} / \n)^{N}, \qquad [x] \otimes [y] \mapsto [x\otimes y].
	\]
	
	\begin{Theorem}[{\cite[Corollary 4.18]{KalmykovSafronov}}]
		\label{thm::res_exact_monoidal}
		The functor $\res_{\g}$ is colimit-preserving and monoidal.
	\end{Theorem}

	Denote by $\forget_G^H\colon \Rep(G) \rightarrow \Rep(H)$ the forgetful functor. We have a diagram
	\begin{equation}
		\label{eq::dyn_comm_diag}
\begin{split}&
		\xymatrix{
			\Rep(G) \ar[d]_{\forget_G^H} \ar[r]^{\free_G} & \HC(G) \ar[d]^{\res_{\g}} \\
			\Rep(H) \ar[r]^{\free_H} & \HC(H)^{\gen}
		}
\end{split}
	\end{equation}
	
	\begin{Theorem}[{\cite{KhoroshkinOgievetsky,ZhelobenkoSalgebras}}]
		\label{thm::dynamical_trivialization}
		The extremal projector defines a natural isomorphism
		\[
		P_V\colon V \otimes \U(\h)^{\gen} \rightarrow \bigl(V \otimes \U(\g)^{\gen} /\n\bigr)^{N}.
		\]
		In other words, the diagram \eqref{eq::dyn_comm_diag} is commutative.
		\begin{proof}
			We will use a slightly different homomorphism as in the citation. By the PBW theorem, we have a vector space isomorphism
			\[
			V \otimes \U(\g)^{\gen} \cong \U(\n_-) \otimes V \otimes U(\h)^{\gen} \otimes \U(\n)
			\]
			for any $G$-representation $V$. In particular,
$
			\n_- \backslash V \otimes \U(\g)^{\gen} / \n \cong V \otimes U(\h)^{\gen}$.
			Therefore, the extremal projector gives a natural isomorphism as in \eqref{eq::extremal_proj_iso_inv_quot}
$
			(V \otimes \U(\g)^{\gen} / \n)^{N} \xrightarrow{\sim} V \otimes U(\h)^{\gen}$.
			Explicitly, the inverse is given by
$
			v \otimes f(\lambda) \mapsto P v \cdot f(\lambda)\in \U(\g)^{\gen} \otimes V / \n
$
			for any $f(\lambda)\in \U(\h)^{\gen}$ and $v\in V$.
		\end{proof}
	\end{Theorem}
	
	In particular, one can translate the natural tensor structure on $\res_{\g}$ to the composition
	\[
	\Rep(G) \xrightarrow{\forget_G^H} \Rep(H) \xrightarrow{\free_H} \HC(H)^{\gen}.
	\]
	
	\begin{Theorem}
		\label{thm::dynamical_quantization}\qquad
		\begin{enumerate}\itemsep=0pt
			\item[$(1)$] There is a collection of maps $J_{VW}(\lambda)$ natural in $V,W\in \Rep(G)$, such that the diagram
			\[
			\xymatrix{
				(V \otimes \U(\h)^{\gen}) \otimes_{\U(\h)^{\gen}} (W \otimes \U(\h)^{\gen}) \ar[r]^-{J_{VW}(\lambda)} \ar[d]^{P_V \otimes P_W} & V \otimes W \otimes \U(\h)^{\gen} \ar[d]^{P_{V\otimes W}} \\
				\res_{\g}(V \otimes \U(\g)) \otimes_{\U(\h)^{\gen}} \res_{\g}(W \otimes \U(\g)) \ar[r] & \res_{\g}(V \otimes W \otimes \U(\g)) \\
			}
			\]
			is commutative, where the lower arrow is the natural tensor structure on $\res_{\g}$. In particular, the collection of $J_{VW}(\lambda)$ satisfies the dynamical twist equation \eqref{eq::dyn_twist}, and the product~${R^{\dyn}_{VW}(\lambda) := J_{WV}(\lambda)^{-1} J_{VW}(\lambda)}$ is a quantum dynamical $R$-matrix from Definition~{\rm\ref{def::all_dynamical}}.
			
			\item[$(2)$] The map $J_{VW}(\lambda)$ has the form
			\[
			J_{VW}(\lambda) \in \id_{V \otimes W} + \hbar \mathrm{U}(\b_-)^{>0} \otimes \mathrm{U}(\b_+)^{>0},
			\]
			where the upper subscript $>0$ means the augmentation ideal.
			
			\item[$(3)$] The coefficient $j_{VW}(\lambda)$ of the first $\hbar$-power of $J_{VW}(\lambda)$ is given by
			\[
			j_{VW}(\lambda) = \frac{f_{\alpha} \otimes e_{\alpha}}{\langle \lambda,h_{\alpha}\rangle }.
			\]
			In particular, the dynamical $R$-matrix $R^{\dyn}_{VW}(\lambda)$ quantizes the standard classical dynamical $r$-matrix \eqref{eq::standard_sol_cdybe}
			\[
			r(\lambda) := j(\lambda) - j^{21}(\lambda) = \sum_{\alpha\in \Delta_+} \frac{f_{\alpha}\wedge e_{\alpha}}{\langle \lambda,h_{\alpha}\rangle }.
			\]
		\end{enumerate}
		\begin{proof}
			We use the argument from \cite{Khoroshkin}. In what follows, we drop the tensor product sign, i.e., for any $x\in V \otimes \U(\g)^{\gen}$, $y\in W \otimes \U(\g)^{\gen}$, we denote
			\[
			xy := [x\otimes_{\U(\g)^{\gen}} y ] \in (V \otimes \U(\g)^{\gen}) \otimes_{\U(\g)^{\gen}} (W \otimes \U(\g)^{\gen}) \cong V \otimes W \otimes \U(\g)^{\gen}.
			\]
			Also, for such $x$ and $y$, we denote by $xP$ (respectively $yP$) their images in the quotient $V \otimes \U(\g)^{\gen} /\n$ (respectively $W \otimes \U(\g)^{\gen} /\n$). This is a well-defined notation thanks to the property~${eP = 0}$ for any $e\in \n$ from Theorem~\ref{thm::extremal_proj}.
			
			Observe that by Theorem~\ref{thm::dynamical_trivialization}, every element in $(V \otimes \U(\g)^{\gen} / \n)^{N}$ can written as $PvP f(\lambda)$, where $f(\lambda) \in \U(\h)^{\gen}$. Then the statement can be presented as follows: for every $v\in V$, $w\in W$, we have
			\begin{equation}
				\label{eq::middle_p_twist_dyn}
				Pv Pw P= PJ_{V,W}(\lambda) (v\otimes w) P.
			\end{equation}
			Since the collection of $J_{VW}(\lambda)$ defines a tensor structure on $\Rep(G) \rightarrow \HC(H)^{\gen}$, it automatically satisfies the dynamical twist equation by Theorem~\ref{thm::dynamical_tannakian}, likewise, the operator
			\[
				R_{VW}(\lambda) := J_{WV}(\lambda)^{-1} J_{VW}(\lambda)
			\]
			automatically satisfies the quantum dynamical Yang--Baxter equation.
			
			Let us write $P$ in the PBW basis
			\[
			P = 1 + \sum_i \hbar^{k_i} f_i e_i a_i(\lambda),
			\]
			where $a_i(\lambda) \in \U(\h)^{\gen}$, $f_i \in \U(\n_-)$, $e_i \in \U(\n)$, and $k_i$ are some negative numbers. It follows from the construction \eqref{eq::extremal_proj_term} that
$
			P\in 1 + \bigl(\n_- \U(\g)^{\gen} \cap \U(\g)^{\gen} \n\bigr)$.
			Therefore, we have $f_i =1$ if and only if $e_i=1$. Then we proceed as follows: consider the middle~$P$ in \eqref{eq::middle_p_twist_dyn}. Using the fact that $P$ commutes with $\U(\h)^{\gen}$, we can push $a_i(\lambda)$ to the right. Using Theorem~\ref{thm::extremal_proj}, we can push $f_i$ to the left until it meets $P$ and becomes zero. Likewise, we can push $e_i$ to the right until it meets $P$ and becomes zero as well.
			
			Let us demonstrate how it works when $f_i$ and $e_i$ are Lie algebra elements. Then we are dealing with the term $\hbar^{-1} f_i e_i a_i(\lambda)$. Let us compute $\hbar^{-1} Pv f_i e_i a_i(\lambda) wP$. First, we push $a_i(\lambda)$ to the right
$
			a_i(\lambda) w P = wPa_i(\lambda + \hbar \wt(w))v$.
			By $Pf_i = 0$, we have
			\[
			Pvf_i = Pf_i v - \hbar P\ad_{f_i}(v) = - \hbar P\ad_{f_i}(v).
			\]
			Likewise, by $e_i P = 0$, we get
			$
			e_i w P= \hbar \ad_{e_i}(w) P$.
			Therefore,
			\[
			\hbar^{-1} Pv f_i e_i a_i(\lambda) wP = -\hbar P \ad_{f_i}(v) \ad_{e_i}(v) P a_i(\lambda + \hbar \wt(w)),
			\]
			which has the necessary form.
			
			In general, we see that each term $\hbar^{-l}f_{\alpha}^l e_{\alpha}^l$ in \eqref{eq::extremal_proj_term} acts with a minimal power of $\hbar^l$, and the second part of the theorem follows. One can easily check that the first power of $\hbar$ is given by the sum of actions of $k=1$ terms of \eqref{eq::extremal_proj_term} in the product \eqref{eq::extremal_proj}, which is
			\[
			j(\lambda) = \sum\limits_{\alpha\in \Delta_+} \frac{f_{\alpha} \otimes e_{\alpha}}{\langle \lambda, h_{\alpha} \rangle },
			\]
			as required.
		\end{proof}
	\end{Theorem}
	The last part motivates the following definition, see \cite{EtingofVarchenkoExchange}.
	\begin{Definition}
		\label{def::standard_dyn_twist}
		The \emph{standard dynamical twist} is the collection $J_{VW}(\lambda)$. The \emph{standard quantum dynamical $R$-matrix} is $R^{\dyn}_{VW}(\lambda)$.
	\end{Definition}
	
	\begin{Example}
		\label{ex:vector_rep_dynamical}
		Let us consider $G= \GL_N$. Denote $\lambda_i = E_{ii} + \hbar (N-i)$. Let $V= \bigl(\C^N\bigr)^*$ the dual of the vector representation
		\[
		V = \mathrm{span}(\phi_1,\dots,\phi_N), \qquad \mathrm{ad}_{E_{ij}}(\phi_k) = -\delta_{ik} \phi_j.
		\]
		Observe that there is a $\gl_N$-equivariant embedding $V \rightarrow \gl_{N+1}$ given by $\phi_i \mapsto E_{N+1,i}$ for $1\leq i \leq N$. Therefore, by formulas \cite[equations~(3.2) and (3.3)]{MolevTransvectorAlgebras}, we have
		\begin{align}
			\label{eq::n_invariant_vectors}
			P \phi_i={}& \sum\limits_{i < i_1 < \dots < i_s \leq N} E_{i_1 i} E_{i_2 i_1} \cdots E_{i_s i_{s-1}} \phi_{i_s} (\lambda_i - \lambda_{i_1})^{-1} \cdots (\lambda_i - \lambda_{i_s})^{-1}\nonumber \\
&{}\in V \otimes \U(\gl_N)^{\gen} /\n.
		\end{align}
		To compute the tensor structure, we need to express
		\[
		P \phi_i P \phi_j = P\phi_i \sum\limits_J E_J \phi_{j_s} \prod_{k\in J} (\lambda_j - \lambda_k)^{-1}
		\]
		in terms of $\{P\phi_a \phi_b \mid a,b=1,\dots, N\}$; here $E_J = E_{j_1 j} \cdots E_{j_s j_{s-1}}$ for $J = \{j_1 <\dots < j_s\}\subset \{j+1,\dots,N\}$. Observe that $E_J \in \U(\n_-)$; since, for any $x\in \n_-$ and $\phi\in V$, we have
		\[
		P\phi x = Px \phi - \hbar P \ad_x(\phi) = -\hbar P\ad_x(\phi)
		\]
		by Theorem~\ref{thm::extremal_proj}, the only nonzero terms are the ones with $J=\varnothing$ and $J=\{i\}$ if $i>j$. Therefore,
		\[
		P\phi_i P\phi_j = P\phi_i \phi_j + \hbar \delta_{i>j} P \phi_j \phi_i \cdot (\lambda_j - \lambda_i)^{-1},
		\]
		and the tensor map $J_{VV}(\lambda)$ is
		\[
			J_{VV}(\lambda) = \mathrm{Id} + \hbar \sum_{i>j} E^V_{ji} \otimes E^V_{ij} \cdot (\lambda_j - \lambda_i)^{-1},
		\]
		where $\bigl\{ E_{ij}^V \bigr\} \subset \End(V)$ are matrix units satisfying $E^V_{ij}\cdot \phi_k = \delta_{jk} \phi_i$. By formula \cite[equation~(36)]{EtingofVarchenkoExchange}, we obtain
		\begin{align*}
			R^{\dyn}_{VV}(\lambda) ={}& \sum\limits_i E^V_{ii} \otimes E^V_{ii} + \sum_{i\neq j} E^V_{ji} \otimes E^V_{ij} \cdot (\lambda_i - \lambda_j)^{-1} + \sum_{i <j} E^V_{ii} \otimes E^V_{jj}\\& - \sum_{i>j} \frac{(\lambda_j - \lambda_i)^2-\hbar^2}{(\lambda_j-\lambda_i)^2} E^V_{jj} \otimes E^V_{ii}.
		\end{align*}
	\end{Example}
	
	\section{Mirabolic subgroup}
	\label{sect::mirabolic}
	
	Denote by $\GL_N$ the group of invertible $N\times N$-matrices and by $\gl_N$ its Lie algebra. We choose a natural basis $\{E_{ij} \mid 1 \leq i,j \leq N\}$ of matrix units in $\gl_N$ with the commutator
	\[
	[E_{ij}, E_{kl}] = \delta_{jk} E_{il} - \delta_{li} E_{kj}.
	\]
	
	For the rest of the paper, we will adopt the following notations. For any $k\leq N$, we consider~${\gl_k \subset \gl_N}$ embedded as the \emph{upper left corner}:
	\[
	\gl_k = \begin{pmatrix}
		* & \dots & * & 0 & \dots & 0 \\
		\vdots & \ddots & \vdots & \vdots & \ddots & \vdots \\
		* & \dots & * & 0 & \dots & 0 \\
		0 & \dots & 0 & 0 & \dots & 0 \\
		\vdots & \ddots & \vdots & \vdots & \ddots & \vdots \\
		0 & \dots & 0 & 0 & \dots & 0
	\end{pmatrix}
	\begin{array}{ c }
		\vphantom{\begin{array}{ c }
				b_N \\ 
				\vdots \\ 
				z_N 
		\end{array}} \\
		\left.\kern-\nulldelimiterspace
		\vphantom{\begin{array}{ c }
				b_N \\ 
				\vdots \\ 
				z_N 
		\end{array}}
		\right\}\text{$N-k.$}
	\end{array}
	\]
	Likewise, for any $k \leq N$, we denote by ${}_k \gl_N$ the \emph{left $k$-th truncation} of $\gl_N$: it is the subalgebra isomorphic to $\gl_{N-k}$, but embedded as the \emph{lower right corner}:
	\[
	{}_k \gl_N = \begin{pmatrix}
		0 & \dots & 0 & 0 & \dots & 0 \\
		\vdots & \ddots & \vdots & \vdots & \ddots & \vdots \\
		0 & \dots & 0 & 0 & \dots & 0 \\
		0 & \dots & 0 & * & \dots & * \\
		\vdots & \ddots & \vdots & \vdots & \ddots & \vdots \\
		0 & \dots & 0 & * & \dots & *
	\end{pmatrix}
	\begin{array}{ c }
		\vphantom{\begin{array}{ c }
				b_N \\ 
				\vdots \\ 
				z_N 
		\end{array}} \\
		\left.\kern-\nulldelimiterspace
		\vphantom{\begin{array}{ c }
				b_N \\ 
				\vdots \\ 
				z_N 
		\end{array}}
		\right\}\text{$N-k$.}
	\end{array}
	\]
	We will also use combinations of these notations. For instance, here is an example of $ {}_1 \gl_4$ inside~$\gl_5$
	\[
	\begin{pmatrix}
		0 & 0 & 0 & 0 & 0 \\
		0 & * & * & * & 0 \\
		0 & * & * & * & 0 \\
		0 & * & * & * & 0 \\
		0& 0 & 0 & 0 & 0
	\end{pmatrix}.
	\]
	
	It turns out that almost all the constructions of the paper involve not the whole group $\GL_N$, but its almost parabolic subgroup
	\begin{Definition}
		The \emph{mirabolic subgroup} $\rM_N$ is the subgroup of $\GL_N$ preserving the last basis vector:
		\[
		\rM_N = \begin{pmatrix}
			* & \dots & * & 0 \\
			\vdots & \ddots & \vdots & \vdots \\
			* & \dots & * & 0 \\
			* & \dots & * & 1
		\end{pmatrix}.
		\]
		The \emph{mirabolic subalgebra} $\m_N$ is the Lie algebra of $\rM_N$ identified with the space of matrices whose last column is zero
		\[
		\m_N = \begin{pmatrix}
			* & \dots & * & 0 \\
			\vdots & \ddots & \vdots & \vdots \\
			* & \dots & * & 0 \\
		\end{pmatrix}.
		\]
	\end{Definition}
	
	We will also adopt notation ${}_i \m_j \subset \m_N$ as in the case of $\gl_N$. Here is an example of ${}_1 \m_4 \subset \m_5$
	\[
	{}_1 \m_4 = \begin{pmatrix}
		0 & 0 & 0 & 0 & 0 \\
		0 & * & * & 0 & 0 \\
		0 & * & * & 0 & 0 \\
		0 & * & * & 0 & 0 \\
		0& 0 & 0 & 0 & 0
	\end{pmatrix}.
	\]
	
	For the rest of the paper, we will use the notations
	\[
			\n_- = \mathrm{span}(E_{ij} \mid N \geq i > j \geq 1), \qquad
			\b = \mathrm{span}(E_{kl} \mid 1 \leq k \leq l \leq N-1)
	\]
	unless specified otherwise. In particular, $\b$ refers to the intersection of $\m_N$ with the Borel subalgebra of $\gl_N$. Likewise, we will denote by ${}_k\n_-$ and ${}_k \b$ the corresponding truncated subalgebras in ${}_k \m_N$.
	
	Let $e = E_{12} + \dots + E_{N-1,N} \in \gl_N$ be the so-called \emph{principal nilpotent element} and $\vec{u}=(u_1,\dots,u_{N-1})\in \C^{N-1}$ be a vector of parameters. Denote by
	\begin{equation}
		\label{eq::principal_nilpotent}
		e_{\vec{u}} := e - u_1 E_{11} - \dots - u_{N-1} E_{N-1,N-1}
	\end{equation}
	While $e_{\vec{u}} \not\in \m_N$, we can still define the following 2-form on $\m_N$
	\[
		\omega(\vec{u})\colon\ \m_N \wedge \m_N \rightarrow k, \qquad x\wedge y \mapsto \Tr(e_{\vec{u}}\cdot [x,y]).
	\]
	
	\begin{Proposition}
		\label{prop::cg_family_equation}
		The form $\omega(\vec{u})$ is non-degenerate, and its inverse $r^{\CG}(\vec{u})$ is uniquely specified by the condition
		\begin{align}
			r^{\CG}(\vec{u})(E_{i+k,i}^*) ={}& -E_{i,i+k-1} - (u_i - u_{i+k-1}) \cdot r^{\CG}(\vec{u})(E_{i+k-1,i}^*)\nonumber\\
&+ \delta_{i>1} r^{\CG}(\vec{u})(E_{i+k-1,i-1}^*).\label{eq::cg_family_equation}
		\end{align}
		for any $1 \leq i < i+k \leq N$, where $E_{ij}^*$ is the dual basis and we consider $r^{\CG}(\vec{u})$ as a map $ \m_N^* \rightarrow \m_N$.
		\begin{proof}
			Note that both $\n_-$ and $\b$ are isotropic subspaces of $\m_N$. Since $\omega(\vec{u})$ is skew-symmetric, it is enough to construct an inverse only of one map, say $\omega(\vec{u})\colon \b \rightarrow \n_-$. Observe that
			\[
			\omega(\vec{u})(E_{i,i+k-1}) = -E_{i+k,i}^* + \delta_{i>1} E_{i+k-1,i-1}^* - (u_i - u_{i+k-1}) \cdot E_{i+k-1,i}^*.
			\]
			Hence, to satisfy $\bigl(r^{\CG}(\vec{u}) \circ \omega(\vec{u}) \bigr) (E_{i,i+k-1}) = E_{i,i+k-1}$, we get the equation \eqref{eq::cg_family_equation}. It can be solved inductively which proves the existence part.
		\end{proof}
	\end{Proposition}
	
	It follows by general arguments (see \cite{EtingofSchiffmannQuantumGroups}) that $r^{\CG}(\vec{u})$ satisfies the classical Yang--Baxter equation \eqref{eq::cybe} for any $\vec{u}$. In particular, for $\vec{u}=(0,\dots,0)$, we get
	\[
	r^{\CG}(0) = -\sum\limits_{N \geq i > j \geq 1} E_{ij} \wedge \sum\limits_{k=1}^j E_{k,k+i-1-j}.
	\]

	This is a certain degeneration of the Cremmer--Gervais $r$-matrix mentioned in the introduction, see \cite{EndelmanHodges}. This motivates the following definition.
	\begin{Definition}
		\label{def::cg_class_r_mat}
		An $h$-\emph{deformed rational Cremmer--Gervais $r$-matrix} is $r^{\CG}(\vec{u})$ for the vector of parameters $\vec{u}\in \C^{N-1}$.
	\end{Definition}
	
	As any classical $r$-matrix, it defines a Poisson--Lie structure on the mirabolic subgroup $\rM_N$ (in this case, even symplectic one); moreover, by considering $r^{\CG}(\vec{u})$ as an element of $\gl_N \wedge \gl_N$, it also defines a Poisson--Lie structure on the whole group $\GL_N$. We will construct its quantization in Section~\ref{sect::kw_red}.

	\section{Yangian}
	\label{sect::basics}
	
	In this section, we recall some facts about the Yangian of the Lie algebra $\gl_N$ and prove some results that we use later in the paper. For details, we refer the reader to \cite{MolevYangians}; note that we use an $\hbar$-version of constructions in question (see Section~\ref{subsect::general_setup}), while in loc.\ cit.\  it is specialized to~${\hbar=1}$. However, all the proofs can be easily generalized to the case of an arbitrary $\hbar$.
	
	\subsection{General definition and properties}

The \emph{$($asymptotic$)$ Yangian} $\rY(\gl_N)$ is an associative algebra generated over $\C[\hbar]$ by the elements~\smash{$t_{ij}^{(r)}$}, where $1\leq i,j \leq N$ and $r\geq 1$. In what follows, we also define \smash{$t_{ij}^{(0)} := \delta_{ij}$}. Introduce the formal Laurent series in $u^{-1}$
	\[
	t_{ij}(u) := \sum\limits_{k \geq 0} t_{ij}^{(r)} u^{-r}.
	\]
	We combine them into a generating matrix
	\[
	T(u) := (t_{ij}(u)) \in \rY(\gl_N) \otimes \End\bigl(\C^N\bigr)\bigl[\bigl[u^{-1}\bigr]\bigr].
	\]
	Let $P\in \End\bigl(\C^N \otimes \C^N\bigr)$ be the permutation matrix. Consider the \emph{Yang $R$-matrix}
	\[
	R(u-v) = \id - \frac{\hbar P}{u-v} \in \End\bigl(\C^N \otimes \C^N\bigr).
	\]
	Then the defining relation of $\rY(\gl_N)$ can be presented in the form
	\[
	R(u-v) T_1(u) T_2(v) = T_2(v) T_1(u) R(u-v),
	\]
	where the equality takes place in $\rY(\gl_N) \otimes \End\bigl(\C^n \otimes \C^N\bigr)$ and $T_1(u) := T(u) \otimes \id$, $T_2(v) = \id \otimes T(v)$. Equivalently, it is given by
	\begin{equation}
		\label{eq::yang_def_rel}
		(u-v) [t_{ij}(u), t_{kl}(v)] = \hbar \frac{t_{kj}(u) t_{il}(v) - t_{kj}(v) t_{il}(u)}{u-v}.
	\end{equation}

	Consider the following $\U(\gl_N)$-valued matrix:
	\begin{equation*}
		E = \begin{pmatrix}
			E_{11} & \dots & E_{1,n-1} & E_{1N} \\
			\vdots & \ddots & \vdots & \vdots \\
			E_{N-1,1} & \dots & E_{N-1,N-1} & E_{N-1,N} \\
			E_{N1} & \dots & E_{N,N-1} & E_{NN}
		\end{pmatrix}.
	\end{equation*}	
	One can define a homomorphism $\rY(\gl_N) \rightarrow \U(\gl_N)$, called the \emph{evaluation homomorphism}, by~${T(u) \mapsto \id + u^{-1} E}$.
	
	Under the evaluation homomorphism, we will be mainly dealing with a slightly different version of the $T$-matrix. Namely, define~${
		L(u) = (L_{ij}(u)) := u\cdot \id + E}$.
	Obviously, we have $L(u) = u\cdot \ev(T(u))$. Almost all the properties of $T(u)$ that we use in the paper are satisfied \emph{mutatis mutandis} by $L(u)$; we will explicitly indicate the cases where it is not true or not trivial.
	
	Let $\{a_1,\dots, a_m\}, \{b_1, \dots, b_m\}$ be two sets of indices. Define a \emph{quantum minor} as the sum over all permutations of $1,2,\dots,m$, see \cite[equations~(1.54) and (1.55)]{MolevYangians}
	\begin{align*}
			t_{b_1\dots b_m}^{a_1\dots a_m}(u) :={} & \sum\limits_{\sigma} \mathrm{sgn}(\sigma) t_{a_{\sigma(1)}b_1}(u) t_{a_{\sigma(2)} b_2}(u-\hbar) \cdots t_{a_{\sigma(m)} b_m} (u-\hbar m + \hbar) \\
			={} & \sum\limits_{\sigma} \mathrm{sgn}(\sigma) t_{a_1 b_{\sigma(1)}} (u-m\hbar+\hbar)\cdots t_{a_{m-1} b_{\sigma(m-1)}} (u-\hbar) t_{a_m b_{\sigma(m)}}(u).
	\end{align*}
	
	Similarly, we can define a quantum minor for $L(u)$
	\begin{align}
			L_{b_1\dots b_m}^{a_1\dots a_m}(u) :={} & \sum\limits_{\sigma} \mathrm{sgn}(\sigma) L_{a_{\sigma(1)}b_1} (u) L_{a_{\sigma(2)} b_2} (u-\hbar) \cdots L_{a_{\sigma(m)} b_m} (u-\hbar m + \hbar)\nonumber \\
			={} & \sum\limits_{\sigma} \mathrm{sgn}(\sigma) L_{a_1 b_{\sigma(1)}}(u-m\hbar+\hbar)\cdots L_{a_{m-1} b_{\sigma(m-1)}} (u-\hbar) L_{a_m b_{\sigma(m)}}(u).		\label{eq::quantum_minor}
	\end{align}
	Observe that
	\begin{equation}
		\label{eq::t_l_minor_relation}
		L_{b_1\dots b_m}^{a_1\dots a_m}(u) = u(u-\hbar)\cdots (u-\hbar m + \hbar)t_{b_1\dots b_m}^{a_1\dots a_m}(u).
	\end{equation}
	
	We will need the following properties of quantum minors.
	\begin{Proposition}\qquad
		\label{prop::minor_basic}
		\begin{enumerate}\itemsep=0pt
			\item[$(1)$]\label{prop::minor_basic_perm} For any permutation $\sigma$, we have
			\[
			L_{b_1\dots b_m}^{a_{\sigma(1)}\dots a_{\sigma(m)}}(u) = \mathrm{sgn}(\sigma) L_{b_1\dots b_m}^{a_1 \dots a_m} = L_{b_{\sigma(1)} \dots b_{\sigma(m)}}^{a_1\dots a_m}(u).
			\]
			In particular, if $a_i = a_j$ or $b_i=b_j$ for some $i\neq j$, then $L_{b_1\dots b_m}^{a_1 \dots a_m} (u) = 0$.
			\item[$(2)$]\label{prop::minor_basic_sum} A quantum minor can be decomposed as {\rm\cite[\emph{Proposition} 1.6.8]{MolevYangians}}
			\begin{align*}
				L_{b_1 \dots b_m}^{a_1 \dots a_m}(u)& =\sum\limits_{l=1}^m (-1)^{m-l} L_{b_1 \dots b_{m-1}}^{a_1 \dots \hat{a}_l \dots a_m}(u) L_{a_l b_m}(u-\hbar m + \hbar) \\
				&= \sum\limits_{l=1}^m (-1)^{m-l} L_{b_1\dots \hat{b}_l \dots b_m}^{a_1\dots a_{m-1}}(u-\hbar) L_{a_m b_l}(u) \\
				&= \sum\limits_{l=1}^m (-1)^{l-1} L_{a_l b_1} (u) L_{b_2\dots b_m}^{a_1 \dots \hat{a}_l \dots a_m} (u-\hbar) \\
				&= \sum\limits_{l=1}^m (-1)^{l-1} L_{a_1 b_l}(u-\hbar m + \hbar) L_{b_1 \dots \hat{b}_l \dots b_m}^{a_2 \dots a_m}(u).
			\end{align*}
		\end{enumerate}
	\end{Proposition}
	
	Quantum minors enjoy the following commutation properties.
	\begin{Proposition}
		\label{prop::minor_comm}
		For $k\neq l$, we have
		\begin{equation*}
			\bigl[E_{kl}, L_{b_1\dots b_m}^{a_1 \dots a_m}(u)\bigr] = \hbar \left(\sum\limits_{i=1}^m \delta_{a_i,l} L_{b_1\dots b_m}^{a_1\dots k \dots a_m}(u) - \sum\limits_{i=1}^m \delta_{k,b_i} L_{b_1\dots l \dots b_m}^{a_1\dots a_m}(u)\right),
		\end{equation*}
		where $k$, resp.\ $l$, are on the $i$-th position.
		\begin{proof}
			It follows from \cite[Proposition 1.7.1]{MolevYangians} that
			\begin{equation*}
				(u-v) \bigl[L_{kl}(u), L_{b_1 \dots b_m}^{a_1 \dots a_m}(v)\bigr] =\hbar \left( \sum\limits_{i=1}^m L_{a_i l}(u) L_{b_1 \dots b_m}^{a_1 \dots k \dots a_m}(v) - \sum\limits_{i=1}^{m} L_{b_1 \dots l \dots b_m}^{a_1 \dots a_m}(v) L_{k b_i}(u) \right),
			\end{equation*}
			where the indices $k$ and $l$ replace the corresponding indices $a_i$ and $b_i$ respectively. Since $L_{ij}(u) = \delta_{ij} \cdot u + E_{ij}$, the statement follows by comparing the coefficients of $u$ on both sides.
		\end{proof}
	\end{Proposition}
	
	A particular case of a quantum minor will be important.
	
	\begin{Definition}
		\label{def::qdet}
		The \emph{quantum determinant} of $L(u)$ (respectively $T(u)$) is
$
		\qdet(L(u)) = L_{1\dots N}^{1\dots N}(u)
$
		\big(respectively $\qdet(T(u)) = T_{1\dots N}^{1\dots N}(u)$\big). In what follows, we will also call $\qdet(L(u))$ the \emph{quantum characteristic polynomial}.
	\end{Definition}
	
	It follows from Proposition~\ref{prop::minor_comm} that $\qdet(L(u))$ is central
$
	[L_{ij}(u), \qdet(L(v))] = 0$  $\forall i,j$.
	
	Observe the matrix $T(u)$ has the form
	\[
	T(u) \in \id + u^{-1} \rY(\gl_N) \otimes \End\bigl(\C^N\bigr)\bigl[\bigl[u^{-1}\bigr]\bigr],
	\]
	hence is invertible. We can explicitly construct its inverse as follows.
	\begin{Definition}[{\cite[Definition 1.9.1] {MolevYangians}}]
		Denote by \smash{$\hat{t}_{ij}(u) := (-1)^{i+j} t_{1\dots \hat{i} \dots N}^{1\dots \hat{j} \dots N}(u)$} the complementary quantum minor. The \emph{quantum comatrix} is $
		\hat{T}(u) := \bigl(\hat{t}_{ij}(u)\bigr)$.
	\end{Definition}
	
	Likewise, one can define $\hat{L}_{ij}(u) := L_{1\dots \hat{i} \dots N}^{1\dots \hat{j} \dots N}(u)$.
	
	\begin{Proposition}[{\cite[Proposition 1.9.2]{MolevYangians}}]
		\label{prop::inverse_comatrix}
		The quantum comatrix satisfies
		\[
		\hat{T}(u) T(u - \hbar N+ \hbar) = \qdet(T(u))\cdot \id.
		\]
		In particular, the inverse of $T(u)$ is given by
		\[
		T(u)^{-1} = \qdet(T(u+\hbar N - \hbar)^{-1}\cdot \hat{T}(u+\hbar N - \hbar).
		\]
	\end{Proposition}
	
	Under the evaluation homomorphism, the inverse of $\ev(T(u))$ can be explicitly presented in terms of $E$ as
	\begin{equation}
		\label{eq::inverse_series}
		\ev\bigl(T(u)^{-1}\bigr) = \bigl(\id + u^{-1} E\bigr)^{-1} = \sum\limits_{k\geq 0} (-1)^k u^{-k} E^k.
	\end{equation}
	
	\subsection{Technical lemmas} In what follows, we will prove some technical statements that we will use later in the paper.
	
	\begin{Proposition}
		\label{prop::minor_ac_comm}
		For any $k\leq l$, we have
		\begin{equation*}
			L_{1\dots k}^{1\dots k}(u) L_{1,\dots,k-1,k}^{1,\dots,k-1,l}(u-\hbar) = L_{1,\dots,k-1,k}^{1,\dots,k-1,l}(u) L_{1\dots k}^{1\dots k}(u-\hbar).
		\end{equation*}
		\begin{proof}
			By \cite[Theorem 1.12.1]{MolevYangians}, the map
			\begin{equation*}
				\rY(\gl_{l-k+1}) \rightarrow \rY(\gl_{l}), \qquad t_{ij}(u) \mapsto t_{1\dots k-1,k-1+j}^{1\dots k-1,k-1+i}(u)
			\end{equation*}
			is an algebra homomorphism. The defining commutation relation \eqref{eq::yang_def_rel} gives
			\[
			[t_{11}(u),t_{l-k+1,1}(v)] = \hbar \frac{t_{l-k+1,1}(u) t_{11}(v) - t_{l-k+1,1}(v) t_{11}(u)}{u-v}.
			\]
			Substituting $u-v=\hbar$, we conclude.
		\end{proof}
	\end{Proposition}
	
	\begin{Proposition}
		\label{prop::minor_ac2_comm}
		For any $a>d$, we have
		\[
		L_{1\dots d-1}^{1\dots d-1}(u) L_{1\dots d-1}^{a,2,\dots,d-1}(u-\hbar) = L_{1 \dots d-1}^{a,2,\dots d-1}(u) L_{1\dots d-1}^{1\dots d-1}(u-\hbar).
		\]
		\begin{proof}
			We can suppose that the statement takes place in $\gl_a$. Consider the automorphism of~$\gl_a$ cyclically permuting the indices
			\[
			\phi\colon\ \{1,2,\dots,a \} \mapsto \{a,1,\dots,a-1\}, \qquad E_{ij} \mapsto E_{\phi(i),\phi(j)}.
			\]
			After some permutation of indices in the minors, the statement of the lemma becomes
			\[
			L_{1,\dots,d-2,a}^{1,\dots,d-2,a}(u) L_{1,\dots,d-2,a}^{1,\dots,d-2,a-1}(u-\hbar) = L_{1,\dots,d-2,a}^{1,\dots,d-2,a-1}(u) L_{1,\dots,d-2,a}^{1,\dots,d-2,a}(u-\hbar).
			\]
			Then we can proceed as in the proof of Proposition~\ref{prop::minor_ac_comm} reducing to
			\begin{gather*}
				[t_{a-d+2,a-d+2}(u),t_{a-d+1,a-d+2}(v)] \\
				\qquad=\hbar \frac{t_{a-d+2,a-d+1}(u) t_{a-d+2,a-d+2}(v) - t_{a-d+2,a-d+1}(v) t_{a-d+2,a-d+2}(u)}{u-v},
			\end{gather*}
			the defining commutation relation.
		\end{proof}
	\end{Proposition}
	
	Let $\psi\colon \n_- \rightarrow \C$ be the character such that $\psi(E_{i+1,i})=1$ for all $1\leq i \leq N-1$, and
	\[
	\n_-^{\psi} = \mathrm{span}(x - \psi(x)\mid x\in \n_-) \subset \U(\gl_N)
	\]
	be the shift of $\n_-$ as in Section~\ref{sect::kirillov_projector}.
	
	\begin{Proposition}
		\label{prop::minor_quotient}
		For all $c\leq N$,
		\[
		L_{1\dots c-1}^{1\dots \hat{l} \dots c}(u) \equiv L_{1\dots l-1}^{1\dots l-1}(u)\ \mathrm{mod}\ \n_-^{\psi} \U(\gl_N).
		\]
		\begin{proof}
			For $l=c$ the lemma is obvious. By Proposition~\ref{prop::minor_basic}, we can permute the upper indices
			\[
			L_{1\dots c-1}^{1\dots \hat{l} \dots c}(u) = (-1)^{c} L_{1\dots c-1}^{c,1,\dots,\hat{l},\dots, c-1}(u).
			\]
			By part (\ref{prop::minor_basic_sum}), the right-hand side is equal to
			\[
			\sum\limits_{k=1}^{c-1} (-1)^{c+k-1} E_{ck} L_{1\dots \hat{k} \dots c-1}^{1\dots \hat{l} \dots c-1}(u).
			\]
			But $E_{lk} \equiv 0$ for $k\neq l-1$ and $E_{l,l-1} \equiv 1$. Hence, the right-hand side is equivalent to $L_{1\dots c-2}^{1\dots \hat{l} \dots c-1}(u)$. The statement follows from the obvious induction.
		\end{proof}
	\end{Proposition}
	
	\begin{Proposition}
		\label{prop::inv_qdet_sum}
		For any $1\leq i,j \leq c$, we have
		\[
		\frac{L_{1\dots \hat{i} \dots c}^{1\dots \hat{j} \dots c}(v) L_{1\dots c}^{1\dots c}(u) - L_{1\dots \hat{i} \dots c}^{1\dots \hat{j} \dots c}(v) L_{1\dots c}^{1\dots c}(u)}{u-v} = \sum\limits_{l=1}^c L_{1\dots \hat {i} \dots c}^{1\dots \hat{l} \dots c}(v) L_{1\dots \hat{l} \dots c}^{1\dots \hat{j} \dots c}(u).
		\]
		
		\begin{proof}
			We can assume that the statement takes place in $\gl_c$. Then, for any $k$, $l$
			\[
				L_{1\dots \hat{k} \dots c}^{1\dots \hat{l} \dots c}(u) = \hat{L}_{kl}(u), \qquad
				L_{1\dots c}^{1\dots c}(u) = \qdet(L(u)).
			\]
			By \eqref{eq::t_l_minor_relation}, we have
			\begin{gather*}
				\hat{L}_{ab}(u) = u(u-\hbar)\cdots (u-\hbar N + 2\hbar) \hat{t}_{ab}, \\
				\qdet(L(u)) = u(u-\hbar) \cdots (u-\hbar N+\hbar) \qdet(T(u))
			\end{gather*}
			for any $a$, $b$. Therefore, if we divide both sides by
			\[
			\qdet(T(u)) \qdet(T(v)) \prod\limits_{i=0}^{N-2} (u - \hbar i) \prod\limits_{i=0}^{N-2} (v - \hbar i),
			\]
			(recall that $\qdet(T(u))$ is central), we get
			\begin{gather*}
				\frac{(u-\hbar N+\hbar) \qdet(T(v))^{-1} \hat{t}_{ij}(v) - (v-\hbar N+ \hbar) \qdet(T(u))^{-1}\hat{t}_{ij}(u)}{u-v} \\
				\qquad= \sum\limits_{l=1}^N \bigl(\qdet(T(v))^{-1} \hat{t}_{il}(v)\bigr) \bigl(\qdet(T(u))^{-1} \hat{t}_{lj}(u)\bigr).
			\end{gather*}
			Let us shift the variables $u \mapsto u + \hbar N-\hbar$, $v \mapsto v + \hbar N-\hbar$
			\begin{gather*}
				\frac{u\cdot \qdet(T(v + \hbar N - \hbar))^{-1} \hat{t}_{ij}(v + \hbar N-\hbar) - v \cdot \qdet(T(u + \hbar N-\hbar))^{-1}\hat{t}_{ij}(u + \hbar N-\hbar)}{u-v} \\
				\qquad= \sum\limits_{l=1}^N \qdet(T(v + \hbar N-\hbar))^{-1} \hat{t}_{il}(v + \hbar N-\hbar)\\
\phantom{\qquad=}{}\cdot \qdet(T(u + \hbar N-\hbar))^{-1} \hat{t}_{lj}(u + \hbar N-\hbar).
			\end{gather*}
			Denote by $\tilde{t}_{ij}(u)$ the matrix entry of $T(u)^{-1}$. Then by Proposition~\ref{prop::inverse_comatrix}, we have
			\[
			\tilde{t}_{ij}(u) = \qdet(T(u + \hbar N-\hbar))^{-1} \hat{t}_{ij} (u + \hbar N-\hbar ).
			\]
			Therefore, the equality above is equivalent to
			\[
			\frac{u\cdot \tilde{t}_{ij}(v) - v \cdot \tilde{t}_{ij}(u)}{u-v} = \sum\limits_{l=1}^N \tilde{t}_{il}(v) \tilde{t}_{lj}(u).
			\]
			In terms of series coefficients
			\[
			\tilde{t}_{ij}(u) = \sum\limits_{r\geq 0} \tilde{t}_{ij}^{(r)} u^{-r},
			\]
			this equality is equivalent to
			\begin{equation}
				\label{eq::inverse_matrix_product}
				\tilde{t}_{ij}^{(r)} = \sum\limits_{l=1}^{N} \tilde{t}_{il}^{(r')} \tilde{t}_{lj}^{(s')}
			\end{equation}
			for every $r' + s' = r$. Now let us recall that the inverse is given by the series \eqref{eq::inverse_series} under the evaluation homomorphism.\footnote{I would like to thank Vasily Krylov for suggesting the argument!} In particular, the element \smash{$\tilde{t}_{ij}^{(r)}$} is the $ij$-entry of the matrix power~$(-E)^r$. Therefore, equation \eqref{eq::inverse_matrix_product} is just the tautological formula \smash{$(-E)^{r'+s'} = (-E)^{r'} \cdot (-E)^{s'}$} in terms of the matrix entries.
		\end{proof}
	\end{Proposition}
	
	\section{Kirillov projector}
	\label{sect::kirillov_projector}
	Recall the notion of the mirabolic subalgebra from Section~\ref{sect::mirabolic}. In this section, we define the Kirillov projector and show its key properties.
	
	Denote by $\n_- = \mathrm{span}(E_{ji}\mid j>i) $ the negative nilpotent subalgebra of $\m_N$. Let $\psi\colon \n_- \rightarrow \C$ be a non-degenerate character of $\n_-$; we will always assume that $\psi(E_{i+1,i}) = 1$ for any $i\leq N-1$. For any $x\in \n_-$, denote by \smash{$x^{\psi} := x - \psi(x)$}. Define the shift \smash{$\n_-^{\psi} = \mathrm{span}\bigl(x^{\psi} \mid x\in \n_-\bigr) \subset \U(\gl_N)$}. In this section, we will consider only those right $\m_N$-modules $X$ that satisfy \smash{$X \n_-^{\psi} \subset \hbar X$}. In particular, the action of \smash{$\hbar^{-1} \n_-^{\psi}$} is closed under the commutator and defines an action of the Lie algebra $\n_-$.
	
	\begin{Definition}
		\label{def::whittaker_module}
		A right $\U(\m_N)$-module $X$ such that $X\n_-^{\psi} \subset \hbar X$ is \emph{Whittaker} if the action of~\smash{$\hbar^{-1} \n_-^{\psi}$} integrates to an action of the negative unipotent subgroup $N_-$. A \emph{Whittaker vector} in a~Whittaker module is an element invariant under $N_-$. The space of Whittaker vectors is denoted by $X^{N_-}$.
	\end{Definition}

	Recall the notion of a quantum minor \eqref{eq::quantum_minor}. For $i>j$ and $u\in \C[\hbar]$, consider the following element in a certain completion of $\U(\m_N)\bigl[\hbar^{-1}\bigr]$
	\begin{equation}
		\label{eq::kirillov_proj_term}
		P^{\psi}_{ij}(u) := \sum\limits_{k\geq 0} (-1)^{(i+j)k} \hbar^{-k} \bigl(-E_{ij}^{\psi}\bigr)^k \frac{L_{j\dots i-1}^{j \dots i-1}(u) \cdots L_{j\dots i-1}^{j \dots i-1}(u - \hbar k+\hbar) }{k!}.
	\end{equation}
	For instance, its action is well defined on any \emph{right} $\U(\m_N)$-module $X$ where $\n_-^{\psi}$ acts locally nilpotently and satisfies \smash{$X\n_-^{\psi} \subset \hbar X$}.
	
	In what follows, we denote by $\vec{u}=(u_1,\dots,u_{N-1}) \in \C[\hbar]^{{N-1}}$ a vector of variables. Let us introduce the main object of the paper.
	\begin{Definition}
		\label{def::kirillov_projector}
		The \emph{Kirillov projector} $P_{\m_N}^{\psi}(\vec{u})$ is the element
		\begin{align}
				P_{\m_N}^{\psi}(\vec{u})={}&\bigl(P^{\psi}_{N,N-1}(u_{N-1}) \bigr) \bigl( P_{N-1,N-2}^{\psi}(u_{N-2}) P_{N,N-2}^{\psi}(u_{N-2}) \bigr)\nonumber\\& \cdots \bigl( P_{21}^{\psi}(u_1) \cdots P_{N1}^{\psi}(u_1)\bigr).\label{eq::kirillov_proj}
		\end{align}
	\end{Definition}

	Denote by
	\begin{equation}
		\label{eq::borel_u}
		\b^{\vec{u}} = \mathrm{span}(E_{ij} + \delta_{ij} \cdot u_i \mid 1\leq i \leq j \leq N-1).
	\end{equation}
	
	This is the main theorem of the paper.
	\begin{Theorem}
		\label{thm::main_thm}
		For any right Whittaker module $W$ over $\U(\m_N)$ satisfying $W\n_-^{\psi} \subset \hbar W$, the Kirillov projector defines a unique linear operator $P(\vec{u})\colon W \rightarrow W$ $($acting on the right$)$ such that
		\begin{gather*}
			P(\vec{u}) (x-\psi(x)) = 0,\qquad x\in \n_-, \\
			(E_{ij} + \delta_{ij}\cdot u_i) P(\vec{u}) = 0,\qquad 1\leq i \leq j \leq N-1,
		\end{gather*}
		and the normalization condition
		$
		w P(\vec{u}) = w
		$
		for any Whittaker vector $w\in W^{N_-}$. In particular, it defines a canonical isomorphism
		\begin{equation}
			\label{eq::kirillov_proj_iso_inv_quot}
			W^{N_-} \rightarrow W/ W \b^{\vec{u}},
		\end{equation}
		with the quotient $W/\b^{\vec{u}} := W/ W\b^{\vec{u}}$.
		\begin{proof}
			Uniqueness is clear: assume we have another operator $P_2(\vec{u})$ satisfying these assumptions. By Proposition~\ref{prop::minor_basic},
			\[
			L_{j\dots i-1}^{j \dots i-1}(u_{j}) = (-1)^{i-j} L_{j+1,\dots,i-1,j}^{j,\dots,i-1}(u_j) = \sum\limits_{l=j}^{i-1} (-1)^{l} L_{j\dots \hat{l} \dots i-1}^{j+1 \dots i-1}(u_j-\hbar) L_{jl}(u_j).
			\]
			Observe that $L_{jl}(u_j) = E_{jl} + \delta_{jl}\cdot u_j$, hence $L_{j\dots i-1}^{j\dots i-1}(u_j) P_2 (\vec{u}) = 0$ and
$
			P_{ij}^{\psi}(u_j) P_2(\vec{u}) = P_2(\vec{u})
$
			for every $i > j$. Therefore,
$
			P_{\m_N}(\vec{u}) P_2(\vec{u}) = P_2(\vec{u})$.
			Likewise, for any vector~$v\in W$, the vector $v P_{\m_N}(\vec{u})$ is Whittaker. Therefore, by the normalization condition on $P_2(\vec{u})$, we have
$
			vP_{\m_N}(\vec{u}) P_2(\vec{u}) = vP_{\m_N}(\vec{u})$.
			In particular, as operators on $W$, we get
$
			P_{\m_N}(\vec{u}) P_2(\vec{u}) = P_{\m_N}(\vec{u})$,
			and so
			$
			P_2(\vec{u}) = P_{\m_N}(\vec{u})$.
			The properties for $P_{\m_N}(\vec{u})$ are proven in Theorems~\ref{thm::proj_n_inv} and~\ref{thm::proj_left_inv}.
		\end{proof}
	\end{Theorem}
	
	\begin{Remark}
		\label{rmk::kirillov_proj_classical}
		Classically, these properties can be interpreted as follows. The element $P_{\m_N}(\vec{u})$ defines a projection
\smash{$
		\n_-^{\psi} \backslash \U(\m_N) \rightarrow \bigl(\n_-^{\psi} \backslash \U(\m_N)\bigr)^{N_-}$}.
		By the second property, it induces an isomorphism
\smash{$
		\n_-^{\psi} \backslash \U(\m_N) / \b^{\vec{u}} \xrightarrow{\sim} \bigl(\n_-^{\psi} \backslash \U(\m_N)\bigr)^{N_-}$}.
		Composing with a natural projection $\U(\m_N) \rightarrow \U(\m_N)/\b^{\vec{u}}$, we get a map
\smash{$
		\U(\m_N) \rightarrow \n_-^{\psi} \backslash \U(\m_N) /\b^{\vec{u}}$}.
		We can take the limit $\hbar \rightarrow 0$ which gives~\smash{$
		\Sym(\m_N) \cong \cO(\m_N^*) \rightarrow \n_-^{\psi} \backslash\cO(\m_N^*)/ \b^{\vec{u}}$}.
		Observe that the target is the space of functions on the closed subvariety of $\m_N^*$ specified by the equations $E_{ij} = \psi(E_{ij})$ for $E_{ij}\in \n_-$ and ${E_{kl} = -\delta_{kl}\cdot u_k}$ for $1 \leq k \leq l$, which is just the point~${e_{\vec{u}}\in \m_N^*}$ of \eqref{eq::principal_nilpotent}. Moreover, by the tensor property from Theorem~\ref{thm::mirabolic_kw_trivialization}, this is an algebra map when $\hbar=0$. Therefore, the classical limit of the Kirillov projector gives a map of varieties~${
		\pt \rightarrow \m_N^*}$, $ \pt \mapsto e_{\vec{u}} \in \m_N^*$.
	\end{Remark}
	
	\subsection[Kirillov projector: right n-invariance]{Kirillov projector: right $\boldsymbol{\n_-^{\psi}}$-invariance}
	\label{subsect::proj_right_inv}
	
	In this subsection, we will prove the first property of Theorem~\ref{thm::main_thm}.
	
	\begin{Theorem}
		\label{thm::proj_n_inv}
		The Kirillov projector satisfies
		\[
		P_{\m_N}^{\psi}(\vec{u}) (x - \psi(x)) = 0\qquad \forall x\in \n_-.
		\]
		In particular, for any Whittaker module $X$, it defines a canonical projection
$
		P_{\m_N}^{\psi}(\vec{u}) \colon X \mapsto X^{N_-}
$
		to the space of Whittaker vectors.
	\end{Theorem}
	
	We will prove it by two-fold induction. In what follows, we will always identify $u = u_1$.
	
	\subsubsection{Global induction}
	
	The first one is the induction on the dimension $N$ of $\m_N$. The base case is $N=2$.
	
	\begin{Proposition}
		\label{prop::gl2_proj_inv}
		We have $P_{21}^{\psi}(u) (E_{21} - 1) = 0$.
		\begin{proof}
			Recall that $L_1^1(u) = E_{11} + u$. Then
$
			L_1^1(u) E_{21} = E_{21} L_1^1(u-\hbar)$.
			Therefore,
			\begin{align}
					P_{21}^{\psi}(u) E_{21}={}& \sum\limits_{k\geq 0} \hbar^{-k} (E_{21} -1)^k \frac{L_1^1(u) \cdots L_1^1(u - \hbar k+\hbar)}{k!} E_{21} \nonumber \\
					={}& \sum\limits_{k\geq 0} \hbar^{-k} (E_{21} -1)^k E_{21} \frac{L_1^1(u-\hbar) \cdots L_1^1(u - \hbar k)}{k!} \nonumber \\
					={}& \sum\limits_{k\geq 0} \hbar^{-k} (E_{21}-1)^{k+1} \frac{L_1^1(u-\hbar) \cdots L_1^1(u - \hbar k)}{k!} \nonumber \\
					&+\sum\limits_{k\geq 0} \hbar^{-k} (E_{21}-1)^{k} \frac{L_1^1(u-\hbar) \cdots L_1^1(u - \hbar k)}{k!}.	\label{eq::proj_gl2_e21}
			\end{align}
			The constant term (i.e., with no powers of $E_{21}-1$) is 1. Combining the coefficients of $(E_{21}-1)^{k+1}$, we get
			\begin{gather*}
				\hbar^{-k} \frac{L_1^1(u-\hbar) \cdots L_1^1(u - \hbar k)}{k!} \cdot \left( 1 + \frac{L_1^1(u-\hbar k-\hbar))}{\hbar(k+1)} \right) \\
				\qquad= \hbar^{-k} \frac{L_1^1(u-\hbar) \cdots L_1^1(u - \hbar k)}{k!} \cdot \frac{E_{11} + u}{\hbar(k+1)} \\
				\qquad= \hbar^{-k-1} \frac{L_1^1(u)\cdots L_1^1(u-\hbar k)}{(k+1)!}.
			\end{gather*}
			Therefore, equation \eqref{eq::proj_gl2_e21} becomes
			\begin{align*}
				P_{21}^{\psi}(u) E_{21} = 1 + \sum\limits_{k \geq 0} \hbar^{-k-1} (E_{21}-1)^{k+1} \frac{L_1^1(u) \cdots L_1^1(u-\hbar k)}{(k+1)!} = P_{21}^{\psi}(u).\tag*{\qed}
			\end{align*}\renewcommand{\qed}{}
		\end{proof}
	\end{Proposition}
	
	Then we proceed to the induction step. Recall the notation ${}_1\m_N$ from Section~\ref{sect::mirabolic}
	\begin{equation}
		\label{eq::mirabolic_lower_corner}
		{}_1\m_{N} =
		\begin{pmatrix}
			0 & 0 &\dots & 0 & 0 \\
			0 & * & \dots & * & 0 \\
			\vdots & \ddots & \vdots & \vdots \\
			0 & * & \dots & * & 0
		\end{pmatrix}
	\end{equation}
	and ${}_1\n_-^{\psi} = \mathrm{span}(E_{ij} \mid i > j \geq 2)$. Denote by ${}_1\vec{u} = (u_2,\dots,u_{N-1})$ the truncated vector of parameters. Assume we proved that the element
	\[
	P_{{}_1\m_N}({}_1\vec{u}) = \bigl(P_{N,N-1}^{\psi}(u_{N-1}) \bigr) \cdots \bigl(P_{32}^{\psi}(u_2) \cdots P_{N2}^{\psi}(u_2) \bigr)
	\]
	is invariant under ${}_1 \n_-^{\psi}$. Denote by
$
	\calP_b(\vec{u}) = P_{{}_1\m_N}({}_1\vec{u})P_{21}^{\psi} \cdots P_{b1}^{\psi}(u)
$
	the $b$-truncation of the Kirillov projector for $\m_N$. Then it would be enough to prove that
\smash{$
	\calP_N (\vec{u}) E_{ij}^{\psi} = 0$}
	for all $i>j$.
	
	\subsubsection{Local induction}
	The second induction is on the truncation $b$.
	\begin{Proposition}
		\label{prop::local_ind}
		For any $b \geq i > j$, we have
		$
		\calP_b(\vec{u}) E^{\psi}_{ij} = 0$.
	\end{Proposition}
	Again, the base case $b=2$ is Proposition~\ref{prop::gl2_proj_inv}. For the rest of the subsection, we assume that we proved for it for $b$. Obviously, $[E_{b+1,1},E_{ij}]=0$ for any $b\geq i > j$. By centrality of $L_{1\dots b}^{1\dots b}(u)$, we also have $\bigl[L_{1\dots b}^{1\dots b}(u),E_{ij}\bigr]=0$ for the same $i$, $j$. Therefore,
	\[
	\calP_{b+1}(\vec{u}) E^{\psi}_{ij} = \calP_{b}(\vec{u}) P_{b+1,1}^{\psi}(u) E^{\psi}_{ij} = \calP_{b}(\vec{u}) E_{ij}^{\psi} P_{b+1,1}^{\psi}(u) = 0
	\]
	by the local induction assumption. Hence, it is enough to prove for $i=b+1$ and $j=b$, since other elements are generated by the corresponding commutators.
	
	For the proof, we need the following lemma.
	\begin{Lemma}
		\label{lm::nilp_to_minor}
		For any $c\leq d \leq b < a$, we have
		\[
		\calP_d(\vec{u}) E_{ac}^{\psi} = (-1)^{c+1} \calP_d(\vec{u}) E_{a1} L_{1\dots c-1}^{1\dots c-1}(u).
		\]
	\end{Lemma}
	
	We will prove it by induction on $d$. Base is $d=2$: the case $c=1$ is trivial, and for $c=2$, we use
	\[
	\bigl[\bigl(E_{21}^{\psi}\bigr)^k,E_{a2}^{\psi}\bigr] = -k\hbar E_{a1} \bigl(E_{21}^{\psi}\bigr)^{k-1},
	\]
	so that
	\begin{align*}
		P_{{}_1\m_N}({}_1\vec{u}) P_{21}^{\psi}(u_1) E_{a2}
		={}& P_{{}_1\m_N}({}_1\vec{u})\sum\limits_{k\geq 0} \bigl(E_{21}^{\psi}\bigr)^k \hbar^{-k} \frac{\prod_{i=0}^{k-1} L_1^1(u_1 - i\hbar)}{k!} E_{a2} \\
		={}&P_{{}_1\m_N}({}_1\vec{u}) E_{a2}^{\psi} P_{21}^{\psi}(u_1) \\
		&- P_{{}_1\m_N}({}_1\vec{u}) \sum\limits_{k\geq 0} \bigl(E_{21}^{\psi}\bigr)^k E_{a1} \hbar^{-k} \frac{\prod_{i=0}^{k}L_1^1(u_1 - i\hbar + \hbar)}{k!} \\
		={}& -\calP_2(\vec{u}) E_{a1} L_1^1(u),
	\end{align*}
	where we used $P_{{}_1\m_N}({}_1\vec{u}) E_{a2}^{\psi} = 0$ by global induction assumption.
	
	\begin{Lemma}
		\label{lm::minor_der}
		For $c \leq d \leq b < a$ and any $k\geq 0$, we have
		\[
		\calP_d(\vec{u}) \bigl(L_{1\dots c}^{a,2,\dots,c}(v) + (-1)^c \psi(E_{ac})\bigr) = \calP_d(\vec{u}) E_{a1} \frac{L_{1\dots c}^{1\dots c}(v) - L_{1\dots c}^{1\dots c}(u)}{v-u}.
		\]
		In particular,
		\[
		\calP_d(\vec{u}) \partial_u^k L_{1\dots c}^{a,2,\dots, c}(u) + (-1)^c \delta_{k,0} \psi(E_{ac}) \calP_d(\vec{u})= \frac{1}{k+1} \calP_d(\vec{u}) E_{a1} \bigl(\partial_u^{k+1} L_{1\dots c}^{1\dots c}(u)\bigr),
		\]
		where by $\partial_u$ we mean the derivative with respect to $u$.
		\begin{proof}
			By Proposition~\ref{prop::minor_basic},
			\[
			L_{1\dots c}^{a,2,\dots,c}(v) = \sum\limits_{l=1}^c (-1)^{l-1} E_{al} L_{1\dots \hat{l} \dots c}^{2 \dots c}(v).
			\]
			Since $E_{al} = E_{al}^{\psi} + \psi(E_{al})$ and $\psi(E_{al})=0$ except for, possibly, $l=c=a-1$, we can rewrite the sum as
			\[
			\sum\limits_{l=1}^{c-1} (-1)^{l-1} E^{\psi}_{al} L_{1\dots \hat{l} \dots c}^{2 \dots c}(v) + (-1)^{c-1} \psi(E_{ac}) L_{1\dots c-1}^{2\dots c}(v).
			\]
			
			Therefore, by Lemma~\ref{lm::nilp_to_minor},
			\[
			\calP_d(\vec{u}) L_{1\dots c}^{a,2,\dots, c}(v) = \calP_d(\vec{u}) E_{a1} \sum\limits_{l=1}^c L_{1\dots l-1}^{1\dots l-1}(u) L_{1\dots \hat{l} \dots c}^{2\dots c}(v) + (-1)^{c-1} \calP_d(\vec{u}) \psi(E_{ac}) L_{1\dots c-1}^{2\dots c}(v).
			\]
			Observe that $[E_{a1},E_{ij}]=0$ for any $c\geq i > j$. By Proposition~\ref{prop::minor_quotient}, the classes of $L_{1\dots l-1}^{1\dots l-1}(u)$ and~\smash{$L_{1\dots c-1}^{1\dots \hat{l} \dots c}(u)$} in the left quotient by the shift $\mathrm{span}(E_{ij} - \psi(E_{ij})\mid c\geq i > j \geq 1)$ are equal for every $l$; therefore, by Proposition~\ref{prop::local_ind},
			\[
			\calP_d(\vec{u}) E_{a1} L_{1\dots l-1}^{1\dots l-1}(u) = \calP_d(\vec{u}) E_{a1} L_{1\dots c-1}^{1\dots \hat{l} \dots c}(u)
			\]
			and the right-hand side is equal to
			\[
			\calP_d(\vec{u}) E_{a1} \sum\limits_{l=1}^{c} L_{1\dots c-1}^{1\dots \hat{l} \dots c}(u) L_{1\dots \hat{l} \dots c}^{2\dots c}(v) + (-1)^{c-1} \psi(E_{ac}) \calP_d(\vec{u}).
			\]
			Recall Proposition~\ref{prop::inv_qdet_sum}
			\[
			\sum\limits_{l=1}^{c} L_{1\dots c-1}^{1\dots \hat{l} \dots c}(u) L_{1\dots \hat{l} \dots c}^{2\dots c}(v) = \frac{L_{1\dots c-1}^{2\dots c}(v) L_{1\dots c}^{1\dots c}(u) - L_{1\dots c-1}^{2\dots c}(u) L_{1\dots c}^{1\dots c}(v)}{u-v}.
			\]
			Again, by Proposition~\ref{prop::minor_quotient}, the class $L_{1\dots c-1}^{2\dots c}(u)$ in the quotient is just 1, so,
			\[
			\calP_d(\vec{u}) \bigl(L_{1\dots c}^{a,2,\dots,c}(v) + (-1)^c \psi(E_{ac})\bigr) = \calP_d(\vec{u}) E_{a1} \frac{L_{1\dots c}^{1\dots c}(v) - L_{1\dots c}^{1\dots c}(u)}{v-u}.
			\]
			The second part of the lemma follows by setting $v=u+t$ and comparing the Taylor series in $t$ on both sides.
		\end{proof}
	\end{Lemma}
	
	\begin{proof}[Proof of Lemma~\ref{lm::nilp_to_minor}]
		First, let us show that it will be enough to prove it for $c=d$. Indeed, let $c<d$, then
		\begin{align*}
		\calP_d(\vec{u}) E_{ac}^{\psi} &= \hbar^{-1}\bigl(\calP_d(\vec{u}) E_{ad}^{\psi} E_{dc}^{\psi} - \calP_d(\vec{u}) 	E_{dc}^{\psi} E_{ad}^{\psi}\bigr) \\
		&	= \hbar^{-1} \calP_d(\vec{u}) E_{ad}^{\psi} E_{dc}^{\psi} = (-1)^{d+1} \hbar^{-1}\calP_d(\vec{u}) E_{a1} L_{1\dots d-1}^{1\dots d-1}(u) E_{dc}^{\psi},
		\end{align*}
		where we used Proposition~\ref{prop::local_ind} in the second equality. By Proposition~\ref{prop::minor_comm},
		\[
\bigl[L_{1\dots d-1}^{1\dots d-1}(u),E_{dc}^{\psi}\bigr] = - \hbar L_{1\dots d-1}^{1\dots d \dots d-1}(u) = (-1)^{d-c} \hbar L_{1\dots d-1}^{1\dots \hat{c} \dots d}(u),
\]
		so
		\[
		(-1)^{d+1} \hbar^{-1} \calP_d(\vec{u}) E_{a1} L_{1\dots d-1}^{1\dots d-1}(u) E_{dc}^{\psi} = (-1)^{c+1} \calP_d(\vec{u}) E_{a1} L_{1\dots d-1}^{1\dots \hat{c} \dots d}(u).
		\]
		The statement then follows from Proposition~\ref{prop::minor_quotient} as in the proof of Lemma~\ref{lm::minor_der}.
		
		So, let us assume that $c=d$. Then
		\begin{gather*}
			\calP_{d-1}(\vec{u}) P^{\psi}_{d1}(u) E_{ad}^{\psi} \\
		\qquad	= \calP_{d-1}(\vec{u})\sum\limits_{k\geq 0} \hbar^{-k} (-1)^{(d+1)k} \bigl(-E_{d1}^{\psi}\bigr)^k \frac{\prod\limits_{i=0}^{k-1} L_{1\dots d-1}^{1\dots d-1}(u-i\hbar )}{k!} E_{ad}^{\psi} \\
			\qquad= \calP_{d-1}(\vec{u}) \Biggl( E_{ad}^{\psi} P^{\psi}_{d1}(u) + E_{a1} \sum\limits_{k\geq 1} \hbar^{-k+1} (-1)^{(d+1)k} \bigl(-E_{d1}^{\psi}\bigr)^{k-1} \frac{\prod\limits_{i=0}^{k-1} L_{1\dots d-1}^{1\dots d-1}(u-i\hbar)}{(k-1)!} \Biggr) \\
		\qquad	=  \calP_{d-1}(\vec{u}) \bigl( E_{ad}^{\psi} P^{\psi}_{d1}(u) + (-1)^{d+1} E_{a1} P^{\psi}_{d1}(u-\hbar) L_{1\dots d-1}^{1\dots d-1}(u) \bigr),
		\end{gather*}
		where we used
		\[
		\bigl[\bigl(-E_{d1}^{\psi}\bigr)^k, E_{ad}^{\psi}\bigr] = k \hbar \bigl(-E_{d1}^{\psi}\bigr)^{k-1} E_{a1}.
		\]
		It follows from the definition that $\bigl[\calP_{d-1}(\vec{u}),E_{ad}^{\psi}\bigr]=0$, therefore, the first summand is zero. Hence,
		\[
		\calP_{d}(\vec{u}) E_{ad}^{\psi} = (-1)^{d+1} \calP_{d-1}(\vec{u})E_{a1} P^{\psi}_{d1}(u-\hbar) L_{1\dots d-1}^{1\dots d-1}(u).
		\]
		Therefore, the statement of the lemma would follow from
		\begin{equation}
			\label{eq::lm_nilp_to_minor_shift}
			\calP_{d-1}(\vec{u}) P_{d1}^{\psi}(u) E_{a1} = \calP_{d-1}(\vec{u}) E_{a1} P_{d1}^{\psi}(u-\hbar).
		\end{equation}
		Let us commute $E_{a1}$ past $P_{d1}^{\psi}(u)$. We need three formulas:
		\begin{itemize}\itemsep=0pt
			\item It follows from Proposition~\ref{prop::minor_comm} that
			\begin{equation*}
				\bigl[L_{1\dots d-1}^{1\dots d-1}(u), E_{a1}\bigr] = - \hbar L_{1 \dots d-1}^{a,2,\dots,d-1}(u).
			\end{equation*}
			
			\item From Proposition~\ref{prop::minor_ac2_comm}, we get
			\[
			L_{1\dots d-1}^{1\dots d-1}(v) L_{1\dots d-1}^{a,2,\dots,d-1}(v-\hbar) = L_{1 \dots d-1}^{a,2,\dots d-1}(v) L_{1\dots d-1}^{1\dots d-1}(v-\hbar).
			\]
			
			\item By Proposition~\ref{prop::minor_comm},
$\bigl[L^{a,2,\dots d-1}_{1 \dots d-1}(u),E_{a1}\bigr] = 0$.
		\end{itemize}
		Using them, one can show that
		\[
		\left[\prod\limits_{i=0}^{k-1} L_{1\dots d-1}^{1\dots d-1}(u - \hbar i),E_{a1}\right] = -k\hbar L_{1\dots d-1}^{a,2,\dots,d-1}(u) \prod\limits_{i=1}^{k-1} L_{1\dots d-1}^{1\dots d-1}(u-\hbar i)
		\]
		Then one uses the induction on $d$. By Lemma~\ref{lm::minor_der}, we have
		\begin{gather*}
			 \calP_{d-1}(\vec{u})P_{d1}^{\psi}(u) E_{a1} \\
			\qquad= \calP_{d-1}(\vec{u})\sum\limits_{k\geq 0} \hbar^{-k} (-1)^{(d+1)k} \bigl(-E_{d1}^{\psi}\bigr)^k \frac{\prod L_{1\dots d-1}^{1\dots d-1}(u-i\hbar)}{k!} E_{a1} \\
			\qquad= \calP_{d-1}(\vec{u})E_{a1} \sum\limits_{k\geq 0} (-1)^{(d+1)k} \hbar^{-k} \bigl(-E_{d1}^{\psi}\bigr)^k \frac{\prod_{i=0}^{k-1} L_{1\dots d-1}^{1\dots d-1}(u-i\hbar)}{k!} \\
			\phantom{\qquad=}- \calP_{d-1}(\vec{u}) L_{1\dots d-1}^{a,2,\dots,d-1}(u) \sum\limits_{k\geq 1} \hbar^{-k+1} (-1)^{(d+1)k} \bigl(-E_{d1}^{\psi}\bigr)^k \frac{\prod_{i=1}^{k-1} L_{1\dots d-1}^{1\dots d-1}(u-i\hbar)}{(k-1)!} \\
			\qquad= \calP_{d-1}(\vec{u})E_{a1} \sum\limits_{k\geq 0} (-1)^{(d+1)k} \hbar^{-k} \bigl(-E_{d1}^{\psi}\bigr)^k \frac{\prod_{i=0}^{k-1} L_{1\dots d-1}^{1\dots d-1}(u-i\hbar)}{k!} \\
		\phantom{\qquad=}	- \calP_{d-1}(\vec{u}) E_{a1} \partial_u L_{1\dots d-1}^{1\dots d-1}(u) \sum\limits_{k\geq 1} \hbar^{-k+1} (-1)^{(d+1)k} \bigl(-E_{d1}^{\psi}\bigr)^k \frac{\prod_{i=1}^{k-1} L_{1\dots d-1}^{1\dots d-1}(u-i\hbar)}{(k-1)!}.
		\end{gather*}
		Recall \eqref{eq::lm_nilp_to_minor_shift}. By term-by-term comparison, the lemma would follow from
		\begin{gather}
				\calP_{d-1}(\vec{u})\left[ E_{a1} \bigl(-E_{d1}^{\psi}\bigr)^k \frac{L_{1\dots d-1}^{1\dots d-1}(u)}{k} - \hbar E_{a1} \partial_u L_{1\dots d-1}^{1\dots d-1}(u) \bigl(-E_{d1}^{\psi}\bigr)^k \right] \nonumber \\
				\qquad=\calP_{d-1}(\vec{u}) E_{a1} \bigl(-E_{d1}^{\psi}\bigr)^k \frac{L_{1\dots d-1}^{1\dots d-1}(u-k\hbar)}{k}.		\label{eq::taylor_form}
		\end{gather}
		For $d=2$, it is obviously satisfied. Assume $d>2$. Then we can identify $E_{d1}^{\psi} = E_{d1}$. Let us study the term \smash{$\hbar E_{a1} \partial_u L_{1\dots d-1}^{1\dots d-1}(u) \bigl(-E_{d1}^{\psi}\bigr)^k$}. Taking derivative of Proposition~\ref{prop::minor_comm}, we obtain
		\[
		\bigl[\partial_u L_{1\dots d-1}^{1\dots d-1}(u),\bigl(-E_{d1}^{\psi}\bigr)^k\bigr] = \hbar k \partial_u L_{1,2,\dots,d-1}^{d,2,\dots,d-1}(u) \bigl(-E_{d1}^{\psi}\bigr)^{k-1}.
		\]
		Hence, by Lemma~\ref{lm::minor_der}, we get
		\begin{gather}
\hbar \calP_{d-1}(\vec{u})E_{a1} \partial_u L_{1\dots d-1}^{1\dots d-1}(u) \bigl(-E_{d1}^{\psi}\bigr)^k\nonumber \\
			\qquad	=\hbar \calP_{d-1}(\vec{u}) \bigl( E_{a1} \bigl(-E_{d1}^{\psi}\bigr)^k \partial_u L_{1\dots d-1}^{1\dots d-1}(u) + k \hbar^2 E_{a1} \partial_u L_{1,2,\dots,d-1}^{d,2,\dots,d-1}(u) \bigl(-E_{d1}^{\psi}\bigr)^{k-1} \bigr)\nonumber \\
				\qquad= \hbar \calP_{d-1}(\vec{u})\bigl(E_{a1} \bigl(-E_{d1}^{\psi}\bigr)^k \partial_u L_{1\dots d-1}^{1\dots d-1}(u) + k \hbar^2 \partial_u L_{1,2,\dots,d-1}^{d,2,\dots,d-1}(u) E_{a1} \bigl(-E_{d1}^{\psi}\bigr)^{k-1} \bigr)\nonumber\\
			\qquad	= \hbar \calP_{d-1}(\vec{u})\biggl( E_{a1} \bigl(-E_{d1}^{\psi}\bigr)^k \partial_u L_{1\dots d-1}^{1\dots d-1}(u) \nonumber \\
\phantom{\qquad	=}{} - \frac{k\hbar^2}{2} \bigl(-E_{d1}^{\psi}\bigr) \partial^2_u L_{1\dots d-1}^{1\dots d-1}(u) E_{a1} \bigl(-E_{d1}^{\psi}\bigr)^{k-1}\biggr).\label{eq::derivative_past_k_power}
		\end{gather}
		
		We will need the following technical statement.
		\begin{Lemma}
			For any $l$, we have
			\begin{gather*}
				\calP_{d-1}(\vec{u})\bigl(-E_{d1}^{\psi}\bigr) \partial^l_u L_{1\dots d-1}^{1\dots d-1}(u) (-E_{a1}) \bigl(-E_{d1}^{\psi}\bigr)^{k-1} \\
				\qquad=\calP_{d-1}(\vec{u}) \biggl((-E_{a1}) \bigl(-E_{d1}^{\psi}\bigr)^k \partial_u^l L_{1\dots d-1}^{1\dots d-1}(u) \\
\phantom{\qquad=}{}  - k\hbar \bigl(-E_{d1}^{\psi}\bigr) \frac{\partial_u^{l+1} L_{1\dots d-1}^{1\dots d-1}(u)}{l+1} (-E_{a1}) \bigl(-E_{d1}^{\psi}\bigr)^{k-1} \biggr).
			\end{gather*}
			\begin{proof}
				Here is a series of equalities where we repeatedly use Lemma~\ref{lm::minor_der}
				\begin{gather*}
					\calP_{d-1}(\vec{u})\bigl(-E_{d1}^{\psi}\bigr)\partial^l_u L_{1\dots d-1}^{1\dots d-1}(u) (-E_{a1}) \bigl(-E_{d1}^{\psi}\bigr)^{k-1} \\
					\qquad= \calP_{d-1}(\vec{u})\bigl(-E_{d1}^{\psi}\bigr)(-E_{a1}) \partial^l_u L_{1\dots d-1}^{1\dots d-1}(u) \bigl(-E_{d1}^{\psi}\bigr)^{k-1} + \hbar \partial_u^l L_{1\dots d-1}^{a,2,\dots,d-1}(u) \bigl(-E_{d1}^{\psi}\bigr)^{k} \\
				\qquad	= \calP_{d-1}(\vec{u})(-E_{a1}) \bigl(-E_{d1}^{\psi}\bigr)^{k} \partial^l_u L_{1\dots d-1}^{1\dots d-1}(u) + (k-1)\hbar \partial_u^l L_{1,2,\dots,d-1}^{d,2,\dots,d-1}(u)(-E_{a1}) \bigl(-E_{d1}^{\psi}\bigr)^{k-1} \\
				\phantom{\qquad=}{}	-\calP_{d-1}(\vec{u}) \hbar (-E_{a1}) \frac{\partial_u^{l+1} L_{1\dots d-1}^{1\dots d-1}(u)}{l+1} \bigl(-E_{d1}^{\psi}\bigr)^k \\
					\qquad=\calP_{d-1}(\vec{u})(-E_{a1}) \bigl(-E_{d1}^{\psi}\bigr)^{k} \partial^l_u L_{1\dots d-1}^{1\dots d-1}(u) \\
				\phantom{\qquad=}{}- (k-1)\hbar \bigl(-E_{d1}^{\psi}\bigr) \frac{\partial_u^{l+1} L_{1\dots d-1}^{1\dots d-1}(u)}{l+1} (-E_{a1}) \bigl(-E_{d1}^{\psi}\bigr)^{k-1} \\
				\phantom{\qquad=}{}	-\calP_{d-1}(\vec{u}) \hbar (-E_{a1}) \frac{\partial_u^{l+1} L_{1\dots d-1}^{1\dots d-1}(u)}{l+1} \bigl(-E_{d1}^{\psi}\bigr)^k.
				\end{gather*}
				Therefore, the claim would follow from
				\[
				\calP_{d-1}(\vec{u})(-E_{a1}) \partial_u^{l+1} L_{1\dots d-1}^{1\dots d-1}(u) \bigl(-E_{d1}^{\psi}\bigr) = \calP_{d-1}(\vec{u})\bigl(-E_{d1}^{\psi}\bigr) \partial_u^{l+1} L_{1\dots d-1}^{1\dots d-1}(u) (-E_{a1}).
				\]
				It can be proven be descending induction on $l$. Indeed, for $l$ large, the corresponding derivative is just zero and the equation is trivial. Assume it holds for the $k$-th derivative. Then, on the one hand,
				\begin{gather*}
					\calP_{d-1}(\vec{u})(-E_{a1}) \partial_u^{k-1} L_{1\dots d-1}^{1\dots d-1}(u) \bigl(-E_{d1}^{\psi}\bigr) \\
				\qquad	= \calP_{d-1}(\vec{u})\bigl( (-E_{a1}) \bigl(-E_{d1}^{\psi}\bigr) \partial_u^{k-1} L_{1\dots d-1}^{1\dots d-1}(u) + \hbar \partial_u^{k-1} L_{1,2,\dots, d-1}^{d,2,\dots,d-1}(u) (-E_{a1}) \bigr) \\
					\qquad=\calP_{d-1}(\vec{u}) \left((-E_{a1}) \bigl(-E_{d1}^{\psi}\bigr) \partial_u^{k-1} L_{1\dots d-1}^{1\dots d-1}(u) + \hbar \bigl(-E_{d1}^{\psi}\bigr) \frac{\partial_u^k L_{1\dots d-1}^{1\dots d-1}(u)}{k} (-E_{a1}) \right),
				\end{gather*}
				on the other hand,
				\begin{gather*}
					\calP_{d-1}(\vec{u})\bigl(-E_{d1}^{\psi}\bigr) \partial_u^{k-1} L_{1\dots d-1}^{1\dots d-1}(u) (-E_{a1}) \\
					\qquad= \calP_{d-1}(\vec{u}) \bigl( \bigl(-E_{d1}^{\psi}\bigr) (-E_{a1}) \partial_u^{k-1} L_{1\dots d-1}^{1\dots d-1}(u) + \hbar \partial_u^{k-1} L_{1,2,\dots d-1}^{a,2,\dots,d-1}(u) \bigl(-E_{d1}^{\psi}\bigr) \bigr) \\
					\qquad=\calP_{d-1}(\vec{u}) \left((-E_{a1}) \bigl(-E_{d1}^{\psi}\bigr) \partial_u^{k-1} L_{1\dots d-1}^{1\dots d-1}(u) + \hbar (-E_{a1}) \frac{\partial_u^k L_{1\dots d-1}^{1\dots d-1}(u)}{k} \bigl(-E_{d1}^{\psi}\bigr) \right),
				\end{gather*}
				and we can apply induction.
			\end{proof}
		\end{Lemma}
		By repeated application of the lemma, one can see that \eqref{eq::derivative_past_k_power} is equal to
		\[
		-\frac{1}{k} \calP_{d-1}(\vec{u}) E_{a1} \bigl(-E_{d1}^{\psi}\bigr)^k \sum\limits_{l=1}^{+\infty} \frac{(-\hbar k)^l}{l!} \partial_u^l L_{1\dots k-1}^{1\dots k-1}(u).
		\]
		However, this is nothing but
		\[
		\frac{1}{k} \calP_{d-1}(\vec{u})E_{a1} \bigl(-E_{d1}^{\psi}\bigr)^k \bigl(L_{1\dots d-1}^{1\dots d-1}(u) - L_{1\dots d-1}^{1\dots d-1}(u-k\hbar)\bigr)
		\]
		Substituting it back into \eqref{eq::taylor_form}, we conclude.
	\end{proof}
	
	Finally, we are ready to prove Proposition~\ref{prop::local_ind}.
	
	\begin{proof}[Proof of Proposition~\ref{prop::local_ind}]
		As it was mentioned in the discussion after the proposition statement, it is enough to show
		\smash{$
		\calP_{b+1}(\vec{u}) E_{b+1,b}^{\psi} = 0$}.
		As usual, let us commute \smash{$E_{b+1,b}^{\psi}$} past \smash{$P_{b+1,1}^{\psi}(u)$}. We need three formulas:
		\begin{itemize}\itemsep=0pt
			\item It follows from Proposition~\ref{prop::minor_comm} that
			\begin{equation*}
				\bigl[L_{1\dots b}^{1\dots b}(u), E_{b+1,b}^{\psi}\bigr] = - \hbar L_{1,\dots,b-1,b}^{1,\dots,b-1,b+1}(u).
			\end{equation*}
			
			\item From Proposition~\ref{prop::minor_ac_comm}, we get
			\[
			L_{1\dots b}^{1\dots b}(v) L_{1,\dots,b-1,b}^{1,\dots,b-1,b+1}(v-\hbar) = L_{1,\dots,b-1,b}^{1,\dots,b-1,b+1}(v) L_{1\dots b}^{1\dots b}(v-\hbar).
			\]
			
			\item By Proposition~\ref{prop::minor_comm},
$
			\bigl[L^{1,\dots,b-1,b+1}_{1,\dots,b-1,b}(u),E_{b+1,b}^{\psi}\bigr] = 0$.
		\end{itemize}
		It implies that
		\[
		\left[\prod_{i=0}^{k-1} L_{1\dots b}^{1\dots b}(u-i\hbar), E_{b+1,b}^{\psi} \right] = -k\hbar L_{1,\dots,b-1,b}^{1,\dots,b-1,b+1}(u) \prod_{i=1}^{k-1} L_{1\dots b}^{1\dots b}(u-i\hbar).
		\]
		Hence,
		\begin{gather*}
			\calP_{b}(\vec{u})P_{b+1,1}^{\psi}(u) E_{b+1,b}^{\psi} \\
			\qquad= \calP_{b}(\vec{u})\sum\limits_{k\geq 0} \bigl(-E_{b+1,1}^{\psi}\bigr)^k (-1)^{(b+2)k} \hbar^{-k} \frac{\prod_{i=0}^{k-1} L_{1\dots b}^{1\dots b}(u-i\hbar)}{k!} E_{b+1,b}^{\psi} \\
		\qquad	=\calP_{b}(\vec{u}) E_{b+1,b}^{\psi} P_{b+1,1}^{\psi}(u) \\
		\phantom{\qquad	=}{}	- \calP_{b}(\vec{u})\sum\limits_{k\geq 1} \bigl(-E_{b+1,1}^{\psi}\bigr)^k (-1)^{(b+2)k} \hbar^{-k+1} L_{1,\dots,b-1,b}^{1,\dots,b-1,b+1}(u) \frac{\prod_{i=1}^{k-1} L_{1\dots b}^{1\dots b}(u-i\hbar)}{(k-1)!}.
		\end{gather*}
		Denote
		\begin{gather}
				A := \calP_{b}(\vec{u}) E_{b+1,b}^{\psi} P_{b+1,1}^{\psi}(u),\nonumber \\
				B := -\calP_{b}(\vec{u})\sum\limits_{k\geq 1} \bigl(-E_{b+1,1}^{\psi}\bigr)^k (-1)^{(b+2)k} \hbar^{-k+1} L_{1,\dots,b-1,b}^{1,\dots,b-1,b+1}(u) \frac{\prod_{i=1}^{k-1} L_{1\dots b}^{1\dots b}(u-i\hbar)}{(k-1)!}	.		\label{eq::ab_terms}
		\end{gather}
		Then, naturally, we need to prove that
		$
		A + B = 0$.
		Let us study the $A$-term first. By Lemma~\ref{lm::nilp_to_minor},
		\[
		\calP_{b}(\vec{u}) E_{b+1,b}^{\psi} = (-1)^{b+1} \calP_{b}(\vec{u}) E_{b+1,1} L_{1\dots b-1}^{1\dots b-1}(u).
		\]
		Since we assume $b>2$, we can identify $E_{b+1,1}^{\psi} = E_{b+1,1}$. Therefore,
		\begin{gather*}
\begin{split}
& \calP_{b}(\vec{u}) E_{b+1,b}^{\psi} P_{b+1,1}^{\psi}(u)
=\calP_{b}(\vec{u}) E_{b+1,1} L_{1\dots b-1}^{1\dots b-1}(u) \\
&	\qquad= \calP_{b}(\vec{u}) (-1)^{b+1} E_{b+1,1} \sum\limits_{k\geq 0} L_{1\dots b-1}^{1\dots b-1}(u) (-E_{b+1,1})^k (-1)^{(b+2)k} \hbar^{-k} \frac{\prod_{i=0}^{k-1} L_{1\dots b}^{1\dots b}(u-i\hbar)}{k!}.
\end{split}
		\end{gather*}
		Let us look at the term
		\[
		\calP_{b}(\vec{u})(-E_{b+1,1}) L_{1\dots b-1}^{1\dots b-1}(u) (-E_{b+1,1})^k.
		\]
		Similarly to the proof of Lemma~\ref{lm::nilp_to_minor}, we need the following technical statement.
		\begin{Lemma}
			\label{lm::gt_push_right_shift}
			For any $l\leq b-1$, we have
			\[
			\calP_{b}(\vec{u})(-E_{b+1,1}) L_{1\dots l}^{1\dots l}(u) (-E_{b+1,1})^k = \calP_{b}(\vec{u}) (-E_{b+1,1})^{k+1} L_{1\dots l}^{1\dots l}(u-k\hbar).
			\]
			\begin{proof}
				Indeed, by Lemma~\ref{lm::nilp_to_minor},
				\begin{gather*}
					\calP_{b}(\vec{u})(-E_{b+1,1}) L_{1\dots l}^{1\dots l}(u) (-E_{b+1,1})^k  \\
					\qquad= \calP_{b}(\vec{u})(-E_{b+1,1})^{k+1} L_{1\dots l}^{1\dots l}(u) + k\hbar L_{1,2,\dots,l}^{b+1,2,\dots,l}(u) (-E_{b+1,1})^{k} \\
					\qquad=\calP_{b}(\vec{u})(-E_{b+1,1})^{k+1} L_{1\dots l}^{1\dots l}(u) - k\hbar (-E_{b+1,1}) \partial_u L_{1\dots l}^{1\dots l}(u) (-E_{b+1,1})^k.
				\end{gather*}
				Then we continue the process by pushing the derivative to the right. For instance,
				\begin{gather*}
					\calP_{b}(\vec{u})(-E_{b+1,1})^{k+1} L_{1\dots l}^{1\dots l}(u) - k\hbar (-E_{b+1,1}) \partial_u L_{1\dots l}^{1\dots l}(u) (-E_{b+1,1})^k \\
					\qquad= \calP_{b}(\vec{u})(-E_{b+1,1})^{k+1} L_{1\dots l}^{1\dots l}(u) - k\hbar (-E_{b+1,1})^{k+1} \partial_u L_{1\dots l}^{1\dots l}(u) \\
					\phantom{\qquad=}{}- \calP_{b}(\vec{u})(k\hbar)^2 \partial_u L_{1,2,\dots,l}^{b+1,2,\dots,l}(u) (-E_{b+1,1})^k \\
					\qquad=\calP_{b}(\vec{u})(-E_{b+1,1})^{k+1} L_{1\dots l}^{1\dots l}(u) - (-E_{b+1,1})^{k+1} k\hbar \partial_u L_{1\dots l}^{1\dots l}(u) \\
					\phantom{\qquad=}{}+ \calP_{b}(\vec{u})(k\hbar)^2 (-E_{b+1,1}) \frac{\partial_u^2 L_{1\dots l}^{1\dots l}(u)}{2} (-E_{b+1,1})^k.
				\end{gather*}
				In the end, we arrive at
				\begin{align*}
					\calP_{b}(\vec{u})(-E_{b+1,1}) L_{1\dots l}^{1\dots l}(u) (-E_{b+1,1})^k
					={}&\calP_{b}(\vec{u}) (-E_{b+1,1})^{k+1} \sum\limits_{j=0}^{+\infty} \frac{(-k\hbar)^j}{j!} \partial_u^j L_{1\dots l}^{1\dots l}(u) \\
					={}&\calP_{b}(\vec{u}) (-E_{b+1,1})^{k+1} L_{1\dots l}^{1\dots l}(u-k\hbar),
				\end{align*}
				where the sum is finite.
			\end{proof}
		\end{Lemma}
		
		Therefore,
		\begin{align}
				A = {}& \calP_{b}(\vec{u}) E_{b+1,b}^{\psi} P_{b+1,1}^{\psi}(u) \nonumber\\
				={}& \calP_{b}(\vec{u}) \sum\limits_{k\geq 0} (-E_{b+1,1})^{k+1} (-1)^{(b+2)(k+1)} \hbar^{-k} L_{1\dots b-1}^{1\dots b-1}(u-k\hbar) \frac{\prod_{i=0}^{k-1} L_{1\dots b}^{1\dots b}(u-i\hbar)}{k!}.			\label{eq::a_term}
		\end{align}
		
		Now let us study $B$ of \eqref{eq::ab_terms}. Consider the term
		\[
		\calP_{b}(\vec{u})(-E_{b+1,1})^k L_{1,\dots,b-1,b}^{1,\dots,b-1,b+1}(u).
		\]
		By Proposition~\ref{prop::minor_basic}, we have
		\[
		L_{1,\dots,b-1,b}^{1,\dots,b-1,b+1}(u) = (-1)^{b+1} \sum\limits_{l=1}^{b} (-1)^{l-1} E_{b+1,l} L_{1\dots \hat{l} \dots b}^{1\dots b-1}(u).
		\]
		Recall that $E_{b+1,,b} = E_{b+1,b}^{\psi} + 1$. By Lemmas~\ref{lm::nilp_to_minor} and~\ref{lm::gt_push_right_shift}, we get
		\begin{gather*}
			\calP_{b}(\vec{u})(-E_{b+1,1})^k L_{1,\dots,b-1,b}^{1,\dots,b-1,b+1}(u) \\
			\qquad= \calP_{b}(\vec{u}) \left( \sum\limits_{l=1}^{b} E_{b+1,l}^{\psi} (-1)^{b+l} (-E_{b+1,1})^k L_{1\dots \hat{l} \dots b}^{1\dots b-1}(u) \right)+ (-E_{b+1,1})^k L_{1\dots b-1}^{1\dots b-1}(u) \\
			\qquad= \calP_{b}(\vec{u})(-1)^{b+2} (-E_{b+1,1}) \left( \sum\limits_{l=1}^{b} L_{1\dots l-1}^{1\dots l-1}(u) (-E_{b+1,1})^k L_{1\dots \hat{l} \dots b}^{1\dots b-1}(u) \right)\\
\phantom{\qquad= }{} + (-E_{b+1,1})^k L_{1\dots b-1}^{1\dots b-1}(u) \\
		\qquad	= \calP_{b}(\vec{u})(-1)^{b+2} (-E_{b+1,1})^{k+1}\left( \sum\limits_{l=1}^{b} L_{1\dots l-1}^{1\dots l-1}(u-k\hbar) L_{1\dots \hat{l} \dots b}^{1\dots b-1}(u) \right)\\
\phantom{\qquad= }{} + (-E_{b+1,1})^k L_{1\dots b-1}^{1\dots b-1}(u).
		\end{gather*}
		By Propositions~\ref{prop::minor_quotient} and~\ref{prop::inv_qdet_sum}, we have
		\begin{gather*}
			\calP_{b}(\vec{u}) (-E_{b+1,1})^{k+1} \sum\limits_{l=1}^{b} L_{1\dots l-1}^{1\dots l-1}(u-k\hbar) L_{1\dots \hat{l} \dots b}^{1\dots b-1}(u) \\
			\qquad=\calP_{b}(\vec{u}) (-E_{b+1,1})^{k+1} \frac{L_{1\dots b-1}^{1\dots b-1}(u-k\hbar) L_{1\dots b}^{1\dots b}(u) - L_{1\dots b}^{1\dots b}(u-k\hbar) L_{1\dots b-1}^{1\dots b-1}(u)}{k\hbar}.
		\end{gather*}
		Therefore,
		\begin{align*}
			B ={}&-\calP_{b}(\vec{u})\sum\limits_{k\geq 1} (-E_{b+1,1})^k (-1)^{(b+2)k} \hbar^{-k+1} L_{1,\dots,b-1,b}^{1,\dots,b-1,b+1}(u) \frac{\prod_{i=1}^{k-1} L_{1\dots b}^{1\dots b}(u-i\hbar)}{(k-1)!} \\
			={}&-\calP_{b}(\vec{u})\sum\limits_{k\geq 1} (-E_{b+1,1})^{k+1} (-1)^{(b+2)(k+1)} \frac{\hbar^{-k}}{k!} L_{1\dots b-1}^{1\dots b-1}(u-k\hbar) \prod_{i=0}^{k-1} L_{1\dots b}^{1\dots b}(u-i\hbar) \\
			&+ \calP_{b}(\vec{u})\sum\limits_{k\geq 1} (-E_{b+1,1})^{k+1} (-1)^{(b+2)(k+1)} \frac{\hbar^{-k}}{k!} L_{1\dots b-1}^{1\dots b-1}(u) \prod_{i=1}^{k} L_{1\dots b}^{1\dots b}(u-i\hbar) \\
			&- \calP_{b}(\vec{u})\sum\limits_{k\geq 1}(-E_{b+1,1})^k (-1)^{(b+2)k} \frac{\hbar^{-(k-1)}}{(k-1)!} L_{1\dots b-1}^{1\dots b-1}(u) \prod_{i=1}^{k-1} L_{1\dots b}^{1\dots b}(u-i\hbar).
		\end{align*}
		By comparing the second and the third rows, we see that their sum is $\calP_{b}(\vec{u})(-E_{b+1,1})$. Likewise, recall \eqref{eq::a_term} that
		\[
		A = -\calP_{b}(\vec{u}) \sum\limits_{k\geq 0} (-E_{b+1,1})^{k+1} (-1)^{(b+2)(k+1)} \hbar^{-k} L_{1\dots b-1}^{1\dots b-1}(u-k\hbar) \frac{\prod_{i=0}^{k-1} L_{1\dots b}^{1\dots b}(u-i\hbar)}{k!}.
		\]
		In particular, the sum of $A$ and the first row of $B$ is equal to $-\calP_{b}(\vec{u})(-E_{b+1,1})$. The proposition follows.
	\end{proof}
	
	\begin{Remark}
		\label{rmk::kirillov_proj_image_quot}
		It follows from the proof that instead of $P_{ij}^{\psi}(u_j)$, one can consider its image under the quotient map. Namely, denote by
		\[
		\bigl[L_{j \dots i-1}^{j \dots i-1}(u_j)\bigr] \in \n_-^{\psi} \backslash \U(\m_j) \cong \U(\b_j),
		\]
		where $\b_j \subset \m_j$ is the positive Borel subalgebra. The image of $L_{j \dots i-1}^{j \dots i-1}(u_j)$ in $\U(\b_j)$. Consider an analog of \eqref{eq::kirillov_proj_term}:
		\begin{equation*}
			\bigl[P^{\psi}_{ij}(u)\bigr] := \sum\limits_{k\geq 0} (-1)^{(i+j)k} \hbar^{-k} \bigl(-E_{ij}^{\psi}\bigr)^k \frac{\bigl[L_{j\dots i-1}^{j \dots i-1}(u)\bigr] \cdots \bigl[L_{j\dots i-1}^{j \dots i-1}(u - \hbar k+\hbar)\bigr]}{k!}
		\end{equation*}
		and $[P_{\m_N}(\vec{u})]$ given by a similar formula as \eqref{eq::kirillov_proj}. It follows that it defines the same element as~$P_{\m_N}(\vec{u})$. Indeed, as usual, one can show it by two-fold induction. Consider the element~${
		P_{{}_1\m_N}({}_1\vec{u}) P_{21}^{\psi}(u_1) \dots P_{b1}^{\psi}(u_1)}$.
		We know that it is invariant under the right action of $E_{cd}^{\psi}$ for $1\leq d < c \leq b$. Assume we proved it for $b$. To show the statement for $b+1$, observe that
		\[
		P^{\psi}_{b+1,1}(u) := \sum\limits_{k\geq 0} (-1)^{(b+2)k} \hbar^{-k} \bigl(-E_{b+1,1}^{\psi}\bigr)^k \frac{L_{1\dots b}^{1 \dots b}(u_1) \cdots L_{1\dots b}^{1 \dots b}(u_1 - \hbar(k-1))}{k!}.
		\]
		Since $L_{1\dots b}^{1\dots b}(v)$ is central in $\gl_b$ and $E_{b+1,1}^{\psi}$ commutes with all $E_{cd}^{\psi}$ as above, the statement follows.
	\end{Remark}
	
	\subsection[Kirillov projector: left b-invariance]{Kirillov projector: left $\boldsymbol{\b^{\vec{u}}}$-invariance}
	
	In this subsection, we will show the second property from Theorem~\ref{thm::main_thm}.
	\begin{Theorem}
		\label{thm::proj_left_inv}
		The Kirillov projector satisfies
		\begin{gather*}
				(E_{ii} + u_i) P_{\m_N}^{\psi}(\vec{u}) = 0,\qquad 1\leq i \leq N-1, \qquad
				E_{ij} P_{\m_N}^{\psi}(\vec{u}) =0,\qquad 1 \leq i < j \leq N-1.
		\end{gather*}
	\end{Theorem}
	
	We will prove this statement by induction on dimension $N$. The base case is $N=2$.
	
	\begin{Proposition}
		\label{prop::proj_left_inv_gl2}
		The Kirillov projector $P_{\m_2}^{\psi}(u_1)$ satisfies
	$
			(E_{11} + u) P_{21}^{\psi}(u) = 0$.
		\begin{proof}
			First, observe that
			\begin{equation}
				\label{eq::left_e21_transform}
				E_{21} P^{\psi}_{21}(u-\hbar) = P_{21}(u).
			\end{equation}
			Indeed, we have
			\begin{align*}
				P_{21}^{\psi}(u) E_{21} ={}& \sum\limits_{k\geq 0} \bigl(E_{21}^{\psi}\bigr)^k \hbar^{-k} \frac{L_1^1(u) \cdots L_1^1(u-\hbar k+\hbar)}{k!} E_{21} \\
				= {}& E_{21} \sum\limits_{k\geq 0} \bigl(E_{21}^{\psi}\bigr)^k \hbar^{-k} \frac{L_1^1(u-\hbar) \cdots L_1^1(u-\hbar k)}{k!} \\
				= {}& E_{21} P_{21}^{\psi}(u-\hbar).
			\end{align*}
			But by Theorem~\ref{thm::proj_left_inv}, we have $P_{21}^{\psi}(u)E_{21} = P_{21}^{\psi}(u)$.
			
			Then, we will use the following commutation relation:
			\[
			\bigl[E_{11},\bigl(E_{21}^{\psi}\bigr)^k\bigr] = -k\hbar E_{21} \bigl(E_{21}^{\psi}\bigr)^{k-1}.
			\]
			Therefore,
			\begin{align*}
				(E_{11}+u) P_{21}^{\psi} = {}& E_{11} \sum \bigl(E_{21}^{\psi}\bigr)^k \hbar^{-k} \frac{L_1^1(u) \cdots L_1^1(u-\hbar k+\hbar)}{k!} \\
				={}& \sum \bigl(E_{21}^{\psi}\bigr)^k \hbar^{-k} (E_{11}+u) \frac{L_1^1(u) \cdots L_1^1(u-\hbar k+\hbar)}{k!} \\
				&- E_{21} \sum \bigl(E_{21}^{\psi}\bigr)^{k-1} \hbar^{1-k} \frac{L_1^1(u) \cdots L_1^1(u-\hbar k+\hbar)}{(k-1)!} \\
				={}& P_{21}^{\psi}(u) L_1^1(u) - E_{21} P_{21}^{\psi}(u-\hbar) L_1^1(u) = 0
			\end{align*}
			by \eqref{eq::left_e21_transform}.
		\end{proof}
	\end{Proposition}

	Now we make an induction step. Recall the notation ${}_1\m_{N}$ from \eqref{eq::mirabolic_lower_corner}. Assume we proved Theorem~\ref{thm::proj_left_inv} for \smash{$P_{{}_1\m_{N}}^{\psi}({}_1\vec{u})$}. Observe that
	\[
	P_{\m_N}^{\psi}(\vec{u}) = P_{{}_1\m_{N}}^{\psi}({}_1\vec{u}) P_{21}^{\psi}(u_1) \cdots P_{N1}^{\psi}(u_1).
	\]
	
	Then the statement is reduced to
	\begin{align*}
			(E_{11} + u_1) P_{\m_N}^{\psi}(\vec{u}) = 0, \qquad
			E_{1b} P_{\m_N}^{\psi}(\vec{u}) = 0,\qquad 1< k \leq N-1.
	\end{align*}
	The first property follows from Proposition~\ref{prop::proj_left_inv_gl2}, as $E_{11}$ commutes with $P_{{}_1\m_{N}}^{\psi}({}_1\vec{u})$. For the second property, let us show first an a~priori weaker statement.
	\begin{Proposition}
		\label{prop::left_inv_proj_quotient}
		We have
$
		P_{{}_1\m_{N}}^{\psi}({}_1\vec{u}) E_{1b} P_{ \m_{N}}^{\psi}(\vec{u}) = 0$.
		In other words, the equality $E_{1b} P_{ \m_{N}}^{\psi}(\vec{u})\allowbreak=0$ holds in the left quotient by ${}_1\n_-^{\psi}$.
		\begin{proof}
			By induction on $N$ and by Proposition~\ref{prop::kirillov_proj_e1n}, we have
			\[
			E_{1b}P_{\m_b}(u_1,u_2,\dots,u_{b-1}) = P_{\m_b}(u_1-\hbar,u_2,\dots,u_{b-1}) L_{1\dots b}^{1\dots b} (u_1).
			\]
			We claim that nothing changes in the quotient when we pass to a larger algebra, namely,
			\begin{gather*}
				P_{{}_1\m_N}^{\psi}({}_1\vec{u}) E_{1b} P_{{}_1\m_N}^{\psi}({}_1\vec{u}) P_{21}^{\psi}(u_1) \cdots P_{b1}^{\psi}(u_1) \\
				\qquad=P_{{}_1\m_N}^{\psi}({}_1\vec{u}) P_{21}^{\psi}(u_1-\hbar) \cdots P_{b1}^{\psi}(u_1-\hbar) L_{1\dots b}^{1\dots b}(u_1).
			\end{gather*}
			It follows from the following technical, but in fact simple lemma.
			\begin{Lemma}
				Assume that $x = \sum_{i=0}^{l-1} E_{1,b-i} \alpha_i$ with $\alpha_i \in \U({}_{b-1-i} \gl_{b-1})$. Then
				\[
				P_{{}_1\m_{N}}^{\psi}({}_1\vec{u}) x P_{b-l+1,b-l}^{\psi} (u_{b-l}) \cdots P_{b,b-l}(u_{b-l}) \in \sum_{i=0}^{l} E_{1,b-i} \cdot ({}_{b-1-i} \gl_{b-1})
				\]
				and
				\[
				P_{{}_1\m_{N}}^{\psi}({}_1\vec{u}) x P_{b+1,b-l+1}^{\psi}(u_{b-l+1}) \cdots P_{N,b-l+1}^{\psi}(u_{b-l+1}) = P_{{}_1\m_{N}}^{\psi}({}_1\vec{u}) x.
				\]
				\begin{proof}
					Let us only show the second property. The condition $\alpha_i \in ({}_{b-1-i} \gl_{b-1})$ means that $\alpha_i$ is the sum of products of elements of the form $E_{\alpha\beta}$ with $b-l+1 \leq \alpha,\beta \leq b-1$. It is clear that for any such $\alpha,\beta$, the commutator $[E_{j,b-l+1},E_{\alpha \beta}]$ for $j>b$ will be a nilpotent element of the form $E_{jc}$ with $b-l+1 \leq c \leq b-1$, in particular, the value of $\psi$ on it is zero. Therefore, we have
\smash{$
					P_{{}_1\m_{N}}^{\psi}({}_1\vec{u}) \mathrm{ad}_{E_{j,b-l+1}}(x) = 0$}.
					Since
$
					P_{{}_1\m_{N}}^{\psi}({}_1\vec{u}) x \xi^{\psi} = P_{{}_1\m_{N}}^{\psi}({}_1\vec{u}) \mathrm{ad}_{\xi}(x)
$
					for any $\xi$ in the negative nilpotent subalgebra of ${}_1 \m_N$, we conclude that every \smash{$P_{j,b-l+1}^{\psi}(u_{b-l+1})$} for $j>b$ acts on \smash{$P_{{}_1\m_{N}}^{\psi}({}_1\vec{u}) x$} by identity. Hence,
					\[
					P_{{}_1\m_{N}}^{\psi}({}_1\vec{u}) x P_{b+1,b-l+1}^{\psi}(u_{b-l+1}) \cdots P_{N,b-l+1}^{\psi}(u_{b-l+1})(u_{b-l+1}) = P_{{}_1\m_{N}}^{\psi}({}_1\vec{u}) x.
					\]
					The first property follows by similar analysis of commutators and is left to the reader.
				\end{proof}
			\end{Lemma}
			Therefore, the claim of the proposition would follow from
			\[
			P_{{}_1\m_N}^{\psi}({}_1\vec{u}) P_{21}^{\psi}(u_1-\hbar) \dots P_{b1}^{\psi}(u_1-\hbar) L_{1\dots b}^{1\dots b}(u_1) P_{b+1,1}^{\psi}(u_1) = 0.
			\]
			For brevity, denote by $\vec{u} - \hbar := (u_1-\hbar,u_2,\dots,u_{N-1})$ and $u=u_1$. Adopting notations of Section~\ref{subsect::proj_right_inv}, it can be reformulated as
			\[
			\calP_{b}(\vec{u}-\hbar) L_{1\dots b}^{1\dots b}(u) P_{b+1,1}^{\psi} (u) = 0.
			\]
			Observe that $b \geq 2$, in particular, $E_{b+1,1} = E_{b+1,1}^{\psi}$. We need two facts. The first one: by Lemma~\ref{lm::minor_der}, we have
			\begin{gather*}
				\calP_{b}(\vec{u}-\hbar) L_{1\dots b}^{1\dots b}(u_1) (-E_{b+1,1}) \\
				\qquad=\calP_{b}(\vec{u}-\hbar) (-E_{b+1,1}) L_{1\dots b}^{1\dots b}(u) + \hbar \calP_{b}(\vec{u}-\hbar) L_{1,2,\dots,b}^{b+1,2,\dots,b}(u) \\
				\qquad=\calP_{b}(\vec{u}-\hbar) (-E_{b+1,1}) L_{1\dots b}^{1\dots b}(u) + \calP_b(\vec{u}-\hbar)\hbar E_{b+1,1} \frac{L_{1\dots b}^{1\dots b}(u) - L_{1\dots b}^{1\dots b}(u-\hbar)}{\hbar} \\
				\phantom{\qquad=}{}- (-1)^{b+2} \hbar \calP_{b}(\vec{u}-\hbar) \\
				\qquad= \calP_{b}(\vec{u}-\hbar) (-E_{b+1,1}) L_{1\dots b}^{1\dots b}(u-\hbar) - (-1)^{b+2} \hbar \calP_{b}(\vec{u}-\hbar).
			\end{gather*}
			The second one: by the same lemma,
			\begin{gather*}
				\calP_{b}(\vec{u}-\hbar) (-E_{b+1,1}) L_{1\dots b}^{1\dots b}(u-\hbar) (-E_{b+1,1})^{k} \\
				\qquad= \calP_{b}(\vec{u}-\hbar) (-E_{b+1,1})^{k+1} L_{1\dots b}^{1\dots b}(u-\hbar) + k\hbar \calP_{b}(\vec{u}-\hbar) L_{1,2,\dots,b}^{b+1,2,\dots,b}(u-\hbar) (-E_{b1})^k
 \\
				\qquad= \calP_{b}(\vec{u}-\hbar) (-E_{b+1,1})^{k+1} L_{1\dots b}^{1\dots b}(u-\hbar) - k\hbar (-E_{b+1,1}) \partial_u L_{1\dots b}^{1\dots b}(u-\hbar) (-E_{b+1,1})^k \\
				\phantom{\qquad=}{}- (-1)^{b+2} k\hbar \calP_{b}(\vec{u}-\hbar) (-E_{b+1,1})^k.
			\end{gather*}
			Proceeding as in the proof of Lemma~\ref{lm::gt_push_right_shift}, we conclude that
			\begin{gather*}
				\calP_{b}(\vec{u}-\hbar) (-E_{b+1,1}) L_{1\dots b}^{1\dots b}(u-\hbar) (-E_{b+1,1})^{k} \\
				\qquad=\calP_{b}(\vec{u}-\hbar) (-E_{b+1,1})^{k+1} L_{1\dots b}^{1\dots b}(u-k\hbar-\hbar) - (-1)^{b+2} k\hbar \calP_{b}(\vec{u}-\hbar) (-E_{b+1,1})^k
			\end{gather*}
			Combining the two facts, we conclude that
			\begin{gather*}
				\calP_{b}(\vec{u}-\hbar) L_{1\dots b}^{1\dots b}(u) (-E_{b+1,1})^{k+1} \\
			\qquad	= \calP_{b}(\vec{u}-\hbar) (-E_{b1})^{k+1} L_{1\dots b}^{1\dots b} (u-k\hbar-\hbar) - (-1)^{b+2} (k+1)\hbar \calP_{b}(\vec{u}-\hbar) (-E_{b+1,1})^k.
			\end{gather*}
			Therefore,
			\begin{gather*}
				\calP_{b}(\vec{u}-\hbar) L_{1\dots b}^{1\dots b}(u) P_{b+1,1}^{\psi}(u) \\
				\qquad= \calP_{b}(\vec{u}-\hbar) \left( 1 + \sum\limits_{k\geq 0} L_{1\dots b}^{1\dots b}(u) (-E_{b+1,1})^{k+1} (-1)^{(b+2)(k+1)} \hbar^{-k-1} \frac{\prod_{i=0}^{k} L_{1\dots b}^{1\dots b}(u-i\hbar)}{k!} \right)\\
				\qquad= \calP_{b}(\vec{u}-\hbar) L_{1\dots b}^{1\dots b}(u) \\
			\phantom{\qquad=}{}	+\sum\limits_{k\geq 0} (-E_{b+1,1})^{k+1} L_{1\dots b}^{1\dots b}(u-k\hbar-\hbar) (-1)^{(b+2)(k+1)} \hbar^{-k-1} \frac{\prod_{i=0}^{k} L_{1\dots b}^{1\dots b}(u-i\hbar)}{(k+1)!} \\
			\phantom{\qquad=}{}	- \sum\limits_{k\geq 0} (-E_{b+1,1})^k (-1)^{(b+2)k} \hbar^{-k} \frac{\prod_{i=0}^{k} L_{1\dots b}^{1\dots b}(u-i\hbar)}{k!}.
			\end{gather*}
			Term-by-term comparison gives that this sum is indeed zero, i.e.,
			\[
			\calP_{b}(\vec{u}-\hbar) L_{1\dots b}^{1\dots b}(u) P_{b+1,1}^{\psi}(u) = 0,
			\]
			and the proposition follows.
		\end{proof}
	\end{Proposition}
	
	It turns out that the weaker version implies the stronger one.
	\begin{Proposition}
		For $2 \leq b \leq N-1$, we have
		$
		E_{1b} P_{\m_N}(\vec{u}) = 0$,
		where the equality is understood as of linear operators on any Whittaker module.
		\begin{proof}
			Let $W$ be a Whittaker module and $w\in W$. Observe that $\mathrm{span}(E_{1i}\mid 2 \leq i \leq N-1)$ can be naturally identified with the dual of the vector representation of ${}_1 \m_N$. Therefore, we can apply Remark~\ref{rmk::mirabolic_tensor_structure} and conclude that
			\begin{align*}
				w E_{1b}P_{{}_1\m_N}({}_1\vec{u}) = \sum\limits_{i=2}^{N-1} w_i P_{{}_1\m_N}({}_1\vec{u}) E_{1i} P_{{}_1\m_N}({}_1\vec{u})
			\end{align*}
			for some $w_i \in W$. Hence,
			\begin{align*}
				w E_{1b} P_{\m_N}(\vec{u})
				={}&w E_{1b}P_{{}_1\m_N}({}_1\vec{u}) P_{21}^{\psi}(u_1) \cdots P_{N1}^{\psi}(u_1) \\
				={}& \sum\limits_{i=2}^{N-1} w_i P_{{}_1\m_N}({}_1\vec{u}) E_{1i} P_{{}_1\m_N}({}_1\vec{u}) P_{21}^{\psi}(u_1) \cdots P_{N1}^{\psi}(u_1) \\
				={}&\sum\limits_{i=2}^{N-1} w_i P_{{}_1\m_N}({}_1\vec{u}) E_{1i} P_{\m_N}(\vec{u})
			\end{align*}
			for some $\{w_i\}\subset W$. Therefore, to prove that it is zero, it is enough to show that
			\[
			P_{{}_1\m_N}({}_1\vec{u}) E_{1i} P_{\m_N}(\vec{u}) = 0
			\]
			for every $2\leq i \leq N-1$. But this is exactly Proposition~\ref{prop::left_inv_proj_quotient}.
		\end{proof}
	\end{Proposition}
	
	\subsection{Kirillov projector: other properties}
	\label{subsect::kirillov_proj_last_col}
	
	In what follows, we denote
$
	\vec{u} - \hbar e_i := (u_1,\dots, u_{i-1},u_i - \hbar, u_{i+1},\dots,u_{N-1})$.
	
	\begin{Proposition}
		\label{prop:proj_var_shift}
		The Kirillov projector satisfies
		\[
		L_{i\dots N-1}^{i+1 \dots N}(u_i - \hbar) P_{\m_N}(\vec{u} - \hbar e_i) = P_{\m_N}(\vec{u}).
		\]
		for every $1 \leq i \leq N-1$.
		\begin{proof}
			By inductive construction of the Kirillov projector, it is enough to prove the statement for $i=1$. Recall Theorem~\ref{thm::main_thm}: the operator $P_{\m_N}(\vec{u})$, satisfying
			\begin{gather*}
				P_{\m_N} (\vec{u})(x-\psi(x)) = 0,\qquad x\in \n_-, \\
				(E_{ij} + \delta_{ij}\cdot u_j) P_{\m_N} (\vec{u}) =0,\qquad 1\leq i \leq j \leq N-1
			\end{gather*}
			together with the normalization condition $P_{\m_N} (\vec{u}) \big|_{W^{N_-}} = \id$, is unique. Therefore, it would be enough to prove that the left-hand side satisfies these properties. The first property is clear, and the normalization condition follows from Proposition~\ref{prop::minor_quotient}. For $2 \leq i \leq j \leq N-1$, the second property follows from
$
			\bigl[E_{ij}, L_{1\dots N-1}^{2\dots N} (u_1-\hbar)\bigr] = 0$.
			Therefore, it is enough to consider the case $i=1$. Let $j=1$. Then
			\[
			\bigl[E_{11}, L_{1\dots N-1}^{2\dots N} (u_1-\hbar)\bigr] = -\hbar L_{1\dots N-1}^{2\dots N} (u_1-\hbar),
			\]
			hence
			\begin{align*}
				(E_{11} + u_1) L_{1\dots N-1}^{2\dots N} (u_1-\hbar) P_{\m_N}(\vec{u}-\hbar e_1)
				&=L_{1\dots N-1}^{2\dots N} (u_1-\hbar) (E_{11}+u_1-\hbar) P_{\m_N}(\vec{u}-\hbar e_1)\\& = 0.
			\end{align*}
			Let $j>1$, then
			\[
			\bigl[E_{1j}, L_{1\dots N-1}^{2\dots N}(u_1-\hbar)\bigr] = (-1)^{j} \hbar L_{1\dots N-1}^{1 \dots \hat{j} \dots N}(u_1-\hbar).
			\]
			Since
			\[
			L_{1\dots N-1}^{1 \dots \hat{j} \dots N}(u_1-\hbar) = \sum\limits_{l=1}^{N-1} (-1)^{l-1} L_{1\dots \hat{l} \dots N-1}^{2\dots \hat{j} \dots N}(u_1-2\hbar) L_{1l}(u_1-\hbar)
			\]
			and each term $L_{1l}(u_1-\hbar)$ acts by zero on $P_{\m_N}(\vec{u}-\hbar)$, we conclude that
			\[
			E_{ij}L_{1\dots N-1}^{2\dots N}(u_1-\hbar) P_{\m_N}(\vec{u}-\hbar) = 0.
			\]
			Therefore, the element
$
			L_{1\dots N-1}^{2\dots N}(u_1-\hbar) P_{\m_N}(\vec{u}-\hbar)
$
			satisfies the conditions of Theorem~\ref{thm::main_thm}, thus coincides with $P_{\m_N} (\vec{u})$.
		\end{proof}
	\end{Proposition}
	
	\begin{Corollary}
		\label{prop::kirillov_proj_e1n}
		We have
		\[
		E_{1N} P_{\m_N}(\vec{u})= P_{\m_N} (\vec{u} - \hbar e_1) L_{1\dots N}^{1\dots N}(u_1).
		\]
		\begin{proof}
			Since $L_{1\dots N}^{1\dots N}(u_1)$ is central, we have
			\[
			P_{\m_N}(\vec{u}-\hbar e_1) L_{1\dots N}^{1\dots N}(u_1) = L_{1\dots N}^{1\dots N}(u_1) P_{\m_N}(\vec{u}-\hbar e_1).
			\]
			By Proposition~\ref{prop::minor_basic}, we have
			\[
			L_{1\dots N}^{1\dots N}(u_1) = \sum\limits_{l=1}^N (-1)^{l-1} L_{1\dots \hat{l}\dots N}^{2\dots N} (u_1-\hbar) L_{1l}(u_1).
			\]
			By Theorem~\ref{thm::proj_left_inv}, the only terms that act non-trivially are for $i=1$ and $i=N$. The former gives
			\begin{align*}
				L_{2\dots N}^{2\dots N}(u_1-\hbar) L_{11}(u_1) P_{\m_N}(\vec{u}-\hbar e_1)
				={}&L_{2\dots N}^{2\dots N}(u_1-\hbar) (L_{11}(u_1-\hbar) + \hbar) P_{\m_N}(\vec{u}-\hbar e_1) \\
				={}&\hbar L_{2\dots N}^{2\dots N}(u_1-\hbar)P_{\m_N}(\vec{u}-\hbar e_1),
			\end{align*}
			while the latter
			\[
			(-1)^{N-1} L_{1\dots N-1}^{2\dots N}(u_1-\hbar) E_{1N} P_{\m_N}(\vec{u}-\hbar e_1).
			\]
			We also have
			\begin{align*}
			L_{1\dots N-1}^{2\dots N}(u_1-\hbar) E_{1N} ={}& E_{1N} L_{1\dots N-1}^{2\dots N}(u_1-\hbar) + \hbar (-1)^{N-1} L_{1\dots N-1}^{1\dots N-1}(u_1-\hbar)\\& - \hbar (-1)^{N-1} L_{2\dots N}^{2\dots N}(u_1-\hbar).
			\end{align*}
			Since
			\[
			L_{1\dots N-1}^{1\dots N-1}(u_1-\hbar) P_{\m_N}(\vec{u}-\hbar e_1) = 0,
			\]
			we have
			\begin{align*}
				L_{1\dots N}^{1\dots N}(u_1) P_{\m_N}(\vec{u}-\hbar e_1)
				={}& \bigl[\hbar L_{2\dots N}^{2\dots N}(u_1-\hbar)+ (-1)^{N-1} E_{1N} L_{1\dots N-1}^{2\dots N}(u_1-\hbar) \\
				&+ \hbar L_{1\dots N-1}^{1\dots N-1}(u_1-\hbar) - \hbar L_{2\dots N}^{2\dots N}(u_1-\hbar)\bigr] P_{\m_N}(\vec{u}-\hbar e_1) \\
				= {}&(-1)^{N-1} E_{1N} L_{1\dots N-1}^{2\dots N}(u_1-\hbar) P_{\m_N}(\vec{u}-\hbar e_1).
			\end{align*}
			Therefore, the claim of the proposition follows from
			\[
			L_{1\dots N-1}^{2\dots N}(u_1-\hbar) P_{\m_N}(\vec{u}-\hbar) = P_{\m_N}(\vec{u}),
			\]
			which is Proposition~\ref{prop:proj_var_shift}.
		\end{proof}
	\end{Corollary}

	It seems that there is no closed expression for the action of $P_{\m_N}(\vec{u})$ on the last column $\u_N$ of~$\gl_N$
	\[
	\u_N = \begin{pmatrix}
		0 & \dots & 0 & * \\
		\vdots & \ddots & \vdots & \vdots \\
		0 & \dots & 0 & *
	\end{pmatrix}
	\]
	except for $E_{1N}$. However, it turns out that in the quotient by the ideal $\n_-^{\psi} \U(\gl_N)$, its image is essentially given by the coefficients of the quantum characteristic polynomial that we normalize~as
	\[
		L_{1\dots N}^{1\dots N}(v+N\hbar - \hbar) = v^N + \sum_{i=0}^{N-1} (-1)^{N-i} A_i v^i.
	\]
	\begin{Proposition}\quad
		\label{prop::last_col_center_iso}
		\begin{enumerate}\itemsep=0pt
			\item[$(1)$] For a special choice of parameters $u_i = (N-i)\hbar$, we have
			\[
			E_{kN} P_{\m_N}(N\hbar-\hbar,\dots,\hbar) \equiv A_{k-1} \qquad \bigl(\mathrm{mod}\ \n_-^{\psi} \U(\gl_N)\bigr).
			\]
			
			\item[$(2)$] For arbitrary parameters $\vec{u}$, there are functions $\alpha_{kj} = \alpha_{kj}(\vec{u})$ and $\beta_{kj} = \beta_{kj}(\vec{u})$ such that
			\[
			E_{kN} P_{\m_N} (\vec{u}) \equiv (A_{k-1} - \beta_{kk}) + \sum_{j < k} \alpha_{kj} (A_{j-1}-\beta_{kj}) \qquad \bigl(\mathrm{mod}\ \n_-^{\psi} \U(\gl_N)\bigr).
			\]
			
			\item[$(3)$] The map
			$
				\n_-^{\psi} \backslash \U(\gl_N) / \b^{\vec{u}} \rightarrow \Z_{\gl_N}$
			induced by the Kirillov projector, is an isomorphism of vector spaces. Here $\Z_{\gl_N}$ is the center of the universal enveloping algebra $\U(\gl_N)$.
		\end{enumerate}
		
		\begin{proof}
			In what follows, we use ``$\equiv$'' to denote equivalence modulo the ideal $\n_-^{\psi} \U(\gl_N)$.
			
			Since $L_{1\dots N}^{1\dots N}(v)$ is central,
			\[
			L_{1\dots N}^{1\dots N}(v) \equiv P_{\m_N} (\vec{u}) L_{1\dots N}^{1\dots N}(v) = L_{1\dots N}^{1\dots N}(v) P_{\m_N} (\vec{u}).
			\]
			By Proposition~\ref{prop::minor_basic}, we have
			\[
			L_{1\dots N}^{1\dots N}(v) = \sum\limits_{l=1}^N (-1)^{N-l} L_{1\dots N-1}^{1\dots \hat{l} \dots N}(v) L_{lN}(v-\hbar N+\hbar).
			\]
			By Proposition~\ref{prop::minor_quotient}, we get
			\begin{align*}
				L_{1\dots N}^{1\dots N}(v) \equiv{} &\sum\limits_{l=1}^N (-1)^{N-l} L_{1\dots l-1}^{1\dots l-1}(v) L_{lN}(v-\hbar N+\hbar) \\
				={}&\sum\limits_{l=1}^N (-1)^{N-l} L_{lN}(v-\hbar N+\hbar) L_{1\dots l-1}^{1\dots l-1}(v).
			\end{align*}
			Using Proposition~\ref{prop::minor_basic}, one can show that
			\begin{align*}
				L_{1\dots l-1}^{1\dots l-1}(v) P_{\m_N} (\vec{u})=(v-u_1) (v-u_2-\hbar) \cdots (v-u_{l-1} - \hbar l+2\hbar) P_{\m_N} (\vec{u}),
			\end{align*}
			therefore,
			\begin{align*}
					L_{1\dots N}^{1\dots N}(v) P_{\m_N} (\vec{u}) \equiv{}& \sum_{l=1}^N \left((-1)^{N-l} \prod_{i=1}^{l-1}(v-u_i - \hbar i+\hbar) \right) E_{lN} P_{\m_N} (\vec{u}) \\
					&+(v-N\hbar + \hbar)\prod_{i=1}^{N-1}(v - u_i -\hbar i + \hbar) P_{\m_N} (\vec{u}).
			\end{align*}
			The first and the second part of the proposition then follow by comparing the coefficients of $v$ on both sides. As for the third part, it follows from Remark~\ref{rmk::mirabolic_tensor_structure} that for any $x,y\in \U(\u_N)$,
			\begin{align*}
				xy P_{\m_N}(\vec{u}) = x P_{\m_N}(\vec{u}) y P_{\m_N}(\vec{u}) + (K_{(1)} \cdot x) P_{\m_N}(\vec{u}) (K_{(2)} \cdot y) P_{\m_N}(\vec{u}),
			\end{align*}
			where $K = K_{(1)} \otimes K_{(2)}$ is some nilpotent matrix acting on $\U(\u_N) \otimes \U(\u_N)$, where we identify~$\U(\u_N)$ with the double quotient \smash{$\n_-^{\psi} \backslash \U(\gl_N) / \b^{\vec{u}}$}. In particular, we can invert $1 + K$. It follows from Theorem~\ref{thm::center} that $\Z_{\gl_N}$ is a free commutative algebra on generators $\{A_i\}$. Since we know that the transformation induced by the Kirillov projector is invertible on generators, the third part follows by induction on the degree of generators.
		\end{proof}
	\end{Proposition}
	
	\begin{Remark}
		\label{rmk::special_param}
		Classically, the first part of Proposition~\ref{prop::last_col_center_iso} can be interpreted as follows. Denote by
		$
		P_{\hbar} := P_{\m_N}(N\hbar-\hbar,\dots,\hbar)
		$
		the Kirillov projector for this special choice of parameters. It gives a linear map
		\[
		P_{\hbar} \colon\ \n_-^{\psi} \backslash \U(\gl_N) \rightarrow \bigl(\n_-^{\psi} \backslash \U(\gl_N)\bigr)^{N_-}.
		\]
		Recall that by Kostant's theorem \cite{Kostant}, the target is isomorphic to the center $\Z_{\gl_N}$. Composing with the quotient map \smash{$\U(\gl_N) \rightarrow \n_-^{\psi} \backslash \U(\gl_N)$}, we get a homomorphism
		$
		\U(\gl_N) \rightarrow \Z_{\gl_N}$.
		Taking limit $\hbar \rightarrow 0$ and identifying $\g \cong \g^*$, we obtain a linear map
$
		P_0\colon \cO(\g) \rightarrow \cO\bigl(\mathbb{A}^N\bigr)$
		where we identified $\cO\bigl(\mathbb{A}^N\bigr)$ with the classical limit of $\Z_{\gl_N}$ using the coefficients $A_k$ as the generators. Observe that they become the coefficients of the usual characteristic polynomial. Moreover, it follows from Remark~\ref{rmk::mirabolic_tensor_structure} that this is an algebra map, since
$
		xy P_{\hbar}= x P_{\hbar} y P_{\hbar} + O(\hbar)$.
		In particular, it corresponds to a map between varieties $\mathbb{A}^N \rightarrow \g$. It is not hard to determine which one: since~$P_0$ acts as
		\begin{gather*}
			(E_{ij} - \psi(E_{ij})) P_{0} =0,\qquad i>j, \qquad
			E_{kl} P_0 = 0,\qquad k \leq l, \\
			E_{kN} P_0 = A_{k-1},\qquad 1 \leq k \leq N,
		\end{gather*}
		this map is given by
		\[
		(A_1,\dots,A_N) \mapsto
		\begin{pmatrix}
			0 & \dots & 0 & A_0 \\
			1 & \dots & 0 & A_{N-1} \\
			\vdots & \ddots & \vdots & \vdots \\
			0 & \dots & 1 & A_{N-1}
		\end{pmatrix}
		\]
		(up to transposition induced by the isomorphism $\gl_N \cong \gl_N^*$). Therefore, the Kirillov projector quantizes the companion matrix which is a form of a Kostant slice.
	\end{Remark}
	
	\section{Kostant--Whittaker reduction and quantization}
	\label{sect::kw_red}
	
	In this section, we consider the Whittaker analog of the parabolic reduction functor from Section~\ref{subsect::parabolic_reduction} that we call (following \cite{BezrukavnikovFinkelberg}) the \emph{Kostant--Whittaker reduction}. This is a finite-dimensional analog of the Drinfeld--Sokolov reduction, see \cite{Arakawa}.
	
	\subsection{Mirabolic setting} Consider the category $\HC(\rM_N)$ of Harish-Chandra bimodules over the mirabolic subalgebra~$\m_N$ as in Section~\ref{subsect::hc_bimod}. Denote by $\Wh(\m_N)$ the category of right Whittaker modules over $\m_N$ from Definition~\ref{def::whittaker_module}. It has a distinguished object \smash{$Q_{\m_N} = \n_-^{\psi} \backslash \U(\m_N)$}. Similarly to \eqref{eq::act_h}, we can define a functor
	\begin{equation*}
		\act_{\C[\hbar]}^{\psi} \colon\ \Mod{\C[\hbar]} \rightarrow \Wh(\m_N), \qquad A \mapsto A \otimes_{\C[\hbar]} Q_{\m_N}.
	\end{equation*}
	Then we have an analog of Kostant's equivalence \cite{Kostant}:
	\begin{Theorem}[mirabolic Kostant's equivalence]
		\label{thm::mirabolic_skryabin}
		The functor $\act_{\C[\hbar]}^{\psi}$ is an equivalence.
		\begin{proof}
			We proceed as in Proposition~\ref{prop::act_equiv}. The functor $(-)^{N_-} \colon \Wh(\m_N) \rightarrow \Mod{\C[\hbar]}$ is right adjoint to \smash{$\act_{\C[\hbar]}^{\psi}$}
			\begin{align*}
				\Hom_{\Wh(\m_N)} \bigl(\act_{\C[\hbar]}^{\psi}(A), M\bigr) &= \Hom_{\Wh(\m_N)}(A \otimes_{\C[\hbar]} Q_{\m_N},M) \\
				&\cong \Hom_{\U(\m_N)} (A \otimes_{\C[\hbar]} Q_{\m_N},M) \\
				&\cong \Hom_{\C[\hbar]}\bigl(A,M^{N_-}\bigr).
			\end{align*}
			The unit of the adjunction is given by
$
			A \mapsto (A \otimes_{\C[\hbar]} Q_{\m_N})^{N_-}$.
			The Kirillov projector gives an isomorphism with the quotient functor as in \eqref{eq::kirillov_proj_iso_inv_quot}
			\[
			(-)_{\b^{\vec{u}}}\colon\ \Wh(\m_N) \rightarrow \Mod{\C[\hbar]}, \qquad M \mapsto M/\b^{\vec{u}}.
			\]
			So, the functor $(-)^{N_-}$ is exact. We have
			\[
			(A \otimes_{\C[\hbar]} Q_{\m_N})^{N_-} \cong (A \otimes_{\C[\hbar]} Q_{\m_N})/ \b^{\vec{u}}
			\]
			and the composition
			\[
			A \mapsto (A \otimes_{\C[\hbar]} Q_{\m_N})^{N_-} \cong (A \otimes_{\C[\hbar]} Q_{\m_N})/ \b^{\vec{u}}
			\]
			is an isomorphism by the PBW theorem. The rest follows by the same arguments as in Proposition~\ref{prop::act_equiv}.
		\end{proof}
	\end{Theorem}
	
	Similarly to \eqref{eq::act_g}, one can define the action functor
	\[
	\act_{\m_N}^{\psi}\colon\ \HC(\rM_N) \rightarrow \Wh(\m_N), \qquad X \mapsto Q \otimes_{\U(\m_N)} X.
	\]
	The following definition is an analog of Definition~\ref{def::parabolic_reduction}:
	\begin{Definition}
		The \emph{mirabolic Kostant--Whittaker reduction} functor $\res_{\m_N}^{\psi}  \colon  \HC(\rM_N)  \rightarrow \ \allowbreak\Mod{\C[\hbar]}$ is the composition
		\[
		\HC(\rM_N) \xrightarrow{\act_{\m_N}} \Wh(\m_N) \xrightarrow{(-)^{N_-}} \Mod{\C[\hbar]}.
		\]
	\end{Definition}
	
	Explicitly, it is given by a quantum Hamiltonian reduction \smash{$X \mapsto \bigl(\n_-^{\psi} \backslash X\bigr)^{N_-}$}. There is a~natural lax monoidal structure on \smash{$\res_{\m_N}^{\psi}$}
	\[
	\bigl(\n_-^{\psi} \backslash X\bigr)^{N_-} \otimes_{\C[\hbar]} \bigl(\n_-^{\psi} \backslash Y\bigr)^{N_-} \rightarrow \bigl(\n_-^{\psi} \backslash X\otimes_{\U(\m_N)} Y\bigr)^{N_-}.
	\]
	As in Theorem~\ref{thm::res_exact_monoidal}, we have the following.
	\begin{Theorem}
		The functor $\res_{\m_N}^{\psi}$ is colimit-preserving and monoidal.
		\begin{proof}
			Follows from exactly the same arguments as \cite[Corollary 4.18]{KalmykovSafronov}.
		\end{proof}
	\end{Theorem}
	
	Recall the monoidal functor $\free \colon \Rep(\rM_N) \rightarrow \HC(\rM_N)$ from Section~\ref{subsect::hc_bimod}. On the one hand, we get a monoidal functor $\Rep(\rM_N) \rightarrow \Mod{\C[\hbar]}$ by composing with \smash{$\res_{\m_N}^{\psi}$}. On the other hand, there is the composition of the forgetful functor with the extensions of scalars from $\C$ to $\C[\hbar]$ that we denote by $\forget\colon \Rep(\rM_N) \rightarrow \Mod{\C[\hbar]}$. So, we have a diagram
	\begin{equation}
		\label{eq::mirabolic_res_diag}
\begin{split}
& \xymatrix{
			\Rep(\rM_N) \ar[r]^{\free} \ar[dr]_{\forget} & \HC(\rM_N) \ar[d]^{\res_{\m_N}^{\psi}} \\
			& \Mod{\C[\hbar]}.
		}
\end{split}
	\end{equation}
	The following result is an analog of Theorem~\ref{thm::dynamical_trivialization}.
	\begin{Theorem}
		\label{thm::mirabolic_kw_trivialization}
		The Kirillov projector defines a natural isomorphism
		\[
		P_{\m_N}(\vec{u})_V \colon\ V \rightarrow \bigl(\n_-^{\psi} \backslash V \otimes \U(\m_N)\bigr)^{N_-}.
		\]
		In other words, the diagram \eqref{eq::mirabolic_res_diag} is commutative.
		\begin{proof}
			By \eqref{eq::kirillov_proj_iso_inv_quot}, we have
			\[
			\bigl(\n_-^{\psi} \backslash V \otimes \U(\m_N)\bigr)^{N_-} \cong \n_-^{\psi} \backslash V \otimes \U(\m_N) / \b^{\vec{u}}.
			\]
			By the PBW theorem, the composition
			\[
			V \rightarrow \bigl(\n_-^{\psi} \backslash \U(\m_N) \otimes V\bigr)^{N_-} \cong \n_-^{\psi} \backslash V \otimes \U(\m_N) / \b^{\vec{u}}
			\]
			is an isomorphism.
		\end{proof}
	\end{Theorem}

	In particular, there is a monoidal structure on $\forget$ induced from $\res_{\m_N}^{\psi}$.
	\begin{Theorem}\quad
		\label{thm::mirabolic_tensor_structure}
		\begin{enumerate}\itemsep=0pt
			\item[$(1)$] There is a collection of maps $F_{VW}(\vec{u})\in \End_{\C[\hbar]}(V[\hbar] \otimes_{\C[\hbar]} W[\hbar])$ natural in $V,W\in \Rep(\rM_N)$, such that the diagram
			\[
			\xymatrixcolsep{5pc}\xymatrix{
				V[\hbar] \otimes_{\C[\hbar]} W[\hbar] \ar[r]^{F_{VW}(\vec{u})} \ar[d]_{(P_{\m_N}(\vec{u}))_V \otimes (P_{\m_N}(\vec{u}))_W} & V[\hbar] \otimes_{\C[\hbar]} W[\hbar] \ar[d]^{(P_{\m_N}(\vec{u}))_{V\otimes W}} \\
				V[\hbar] \otimes_{\C[\hbar]} W[\hbar] \ar[r] & V[\hbar] \otimes_{\C[\hbar]} W[\hbar]
			}
			\]
			is commutative, where the lower right arrow is the natural tensor structure on $\res_{\m_N}^{\psi}$. In particular, the collection of $F_{VW}(\vec{u})$ satisfies the twist equation \eqref{eq::twist}, and $R^{\CG}_{VW} (\vec{u}):= F_{WV}(\vec{u})^{-1} F_{VW}(\vec{u})$ is a solution to the quantum Yang--Baxter equation.
			\item[$(2)$] The map $F_{VW}(\vec{u})$ has the form
			\[
			F_{VW} (\vec{u})\in \id_{V\otimes W} + \hbar \mathrm{U}(\n_-)^{>0} \otimes \mathrm{U}(\b)^{>0},
			\]
			where the upper subscript $>0$ means the augmentation ideal.
			\item[$(3)$] The inverse $F_{VW}(\vec{u})^{-1}$ also satisfies
			\[
			F_{VW}(\vec{u})^{-1} \in \id_{V\otimes W} + \hbar \mathrm{U}(\n_-)^{>0} \otimes \mathrm{U}(\b)^{>0}.
			\]
		\end{enumerate}
		\begin{proof}
			We proceed exactly as in Theorem~\ref{thm::dynamical_quantization}. For brevity, denote $P = P_{\m_N}(\vec{u})$. Every element in \smash{$\bigl(\n_-^{\psi} \backslash V\otimes \U(\m_N) \bigr)^{N_-}$} can be presented as $PvP$ for $v\in V$ by Theorem~\ref{thm::mirabolic_kw_trivialization}. Therefore, we need to show that
			\begin{equation}
				\label{eq::middle_p_twist}
				PvPwP = PF(v\otimes w) P.
			\end{equation}
			
			First, let us choose the following PBW basis: $\bigl\{E_{ij}^{\psi}\bigr\}$ for $\U(\n_-)$ and
			\[
			(E_{kl}, E_{kk} + u_k \mid 1 \leq k < l \leq N-1)
			\]
			for $\U(\b)$ (exactly in this order for the latter). It is clear from the construction \eqref{eq::kirillov_proj} and Remark~\ref{rmk::kirillov_proj_image_quot} that
			\begin{equation}
				\label{eq::kirillov_proj_normalization}
				P \in 1 + \n_-^{\psi} \U(\m_N) \cap \U(\m_N) \b^{\vec{u}}
			\end{equation}
			(recall the notation $\b^{\vec{u}}$ from \eqref{eq::borel_u}). Let us write it in this PBW basis
			\begin{equation*}
				P = 1 + \hbar^{k_i}\sum\limits_i f^{\psi}_i e_i^u,
			\end{equation*}
			where $f^{\psi}_i \in \U(\n_-)$, $e_i^u \in \U(\b)$, and $k_i$ is some negative number. This is well defined, since a~part of $P$ that acts non-trivially is actually finite. It is clear from \eqref{eq::kirillov_proj_normalization} that $e^u_i=1$ if and only if $f^{\psi}_i=1$.
			
			Now consider the middle $P$ in \eqref{eq::middle_p_twist}. As in Theorem~\ref{thm::dynamical_quantization}, we push the $f_i^{\psi}$-term to the left until it meets $P$ and becomes zero. Likewise, we push the $e_i^u$-term to right until it meets $P$ and becomes zero as well. Since every term
			\[
			\hbar^{-k} \bigl(-E_{ij}^{\psi}\bigr) (-1)^{(i+j)k} \frac{\prod_{i=0}^{k-1} L_{j \dots i-1}^{j \dots i-1}(u-\hbar i + \hbar)}{k!}
			\]
			in \eqref{eq::kirillov_proj_term} acts with the power at least $\hbar^k$, the second part of the theorem follows. The third part follows from the fact that $F_{VW}(\vec{u})$ is upper-triangular in the natural filtration on $V\otimes W$ induced from a $\n_-$-filtration on $V$ and $\b$-filtration on $W$.
		\end{proof}
	\end{Theorem}
	
	\begin{Remark}
		\label{rmk::mirabolic_tensor_structure}
		In general, let $X$ be any Whittaker module and $V$ any representation of $\rM_N$. Then for any $x\in X$ and $v\in V$,
		\[
		PxvP = \sum\limits_i Px_i P v_i P
		\]
		for some $x_i\in X$, $v_i \in V$. Indeed, recall that by Definition~\ref{def::whittaker_module}, the action of $\hbar^{-1} \n_-^{\psi}$ integrates to that of $N_-$. Therefore, we can proceed exactly as in the proof of Theorem~\ref{thm::mirabolic_tensor_structure} and use the fact that the analog of $F_{VW}(\vec{u})$ in this case is invertible as well.
	\end{Remark}
	
	\begin{Example}
		Let $N=2$. The Kirillov projector is given by
		\[
		\sum\limits_{k\geq 0} \bigl(E_{21}^{\psi}\bigr)^k \hbar^{-k} \frac{L_1^1(u) \cdots L_1^1(u-\hbar k + \hbar)}{k!}.
		\]
		Then the procedure described in the proof of Theorem~\ref{thm::mirabolic_tensor_structure} gives the following formula for the twist:
		\begin{equation}
			\label{eq::cg_twist_gl2}
			F(u)= \sum\limits_{k\geq 0} \frac{(-\hbar)^k}{k!} E_{21}^k \otimes \prod_{i=0}^{k-1} (E_{11} - i).
		\end{equation}
		Observe that it does not depend on $u$. A version of it appeared in \cite{GerstenhaberGiaquintoSchack}, see also \cite{KhoroshkinStolinTolstoy}. An additional parameter in the formulas from loc.\ cit.\  ($t$ in the former and $\lambda$ in the latter) corresponds to a~choice of a non-degenerate character $\psi\colon \n_- \rightarrow \C$ that we assumed to be $\psi(E_{21})=1$ for simplicity.
	\end{Example}
	
	\begin{Definition}
		\label{def::cremmer_gervais}
		Let $h = \sum_{i=1}^{N-1} u_i E_{ii} \in \m_N$ be a diagonal matrix corresponding to a vector of parameters $\vec{u} \in \C^{N-1}$. The $h$\emph{-deformed rational Cremmer--Gervais twist} $F^{\CG}(\vec{u})$ is the collection of isomorphisms $F_{VW}(\vec{u})$. The $h$\emph{-deformed rational Cremmer--Gervais $R$-matrix} is~$R^{\CG}(\vec{u}):= \bigl(F_{WV}^{\CG}(\vec{u})\bigr)^{-1} F_{VW}^{\CG}(\vec{u})$.
	\end{Definition}
	
	To justify the name, let us compute the semiclassical limit of $R_{VW}^{\CG}(\vec{u})$:
	\begin{Theorem}
		\label{eq::cg_quantization}
		The tensor structure of Theorem~{\rm\ref{thm::mirabolic_tensor_structure}} quantizes a family of rational Cremmer--Gervais $r$-matrices $r^{\CG}(\vec{u})$ from Definition~{\rm\ref{def::cg_class_r_mat}}.
		\begin{proof}
			We use notations and conventions from the proof of Theorem~\ref{thm::mirabolic_tensor_structure}. Denote by
			\[
			f_{VW} = \sum\limits_{N \geq i > j \geq 1} E_{ij} \otimes \alpha_{ij},
			\]
			where $\alpha_{ij} \in \b \subset \m_N$. It will follow from the proof that the semiclassical limit depends only on the zeroth $\hbar$-power of $\vec{u}$, therefore, we can assume without loss of generality that $\vec{u}$ does not depend on $\hbar$.
			
			Recall that by Proposition~\ref{prop::cg_family_equation}, it is enough to prove that these elements satisfy the relation
			\[
			\alpha_{i+k,i} = -E_{i,i+k-1} - (u_i-u_{i+k-1}) \alpha_{i+k-1,i} + \delta_{i > 1} \alpha_{i+k-1,i-1}.
			\]
			It will follow from the following lemma.
			\begin{Lemma}
				Let
				\[
				PvPwP = Pv + P\biggl( \sum\limits_{N \geq i > j \geq 1} \ad_{E_{ij}}(v) \beta_{ij} + O\bigl(\hbar^2\bigr)\biggr)wP.
				\]
				Then, up to higher-order terms,
				\begin{enumerate}\itemsep=0pt
					\item[$(1)$] $\beta_{ij}$ are sums of diagonal minor $L_{a\dots b}^{a\dots b}(u_a)$ with some coefficients;
					\item[$(2)$] $\beta_{ij}$ satisfy the relation
					\begin{align}
							\beta_{i+k,i} ={}& (-1)^{k+1}\bigl(L_{i\dots i+k-1}^{i\dots i+k-1}(u_i) - (u_i - u_{i+k-1}) L_{i \dots i+k-2}^{i \dots i+k-2}(u_i)\bigr)\nonumber \\
							&-(u_i - u_{i+k-1}) \beta_{i+k-1,i} + \beta_{i+k-1,i-1}.\label{eq::beta_relation}
					\end{align}
				\end{enumerate}
				\begin{proof}
					Before proving the statement, let us show that
					\[
					Pv L_{k\dots l}^{k \dots l}(u_k) x P \in \hbar \sum\limits_{i=k}^l (-1)^{i-k} \prod_{m=i+1}^l (u_k - u_m) \cdot \ad_{E_{ki}}(x) P + O\bigl(\hbar^2\bigr)
					\]
					for any $x\in \U(\m_N) \otimes W$ (in particular, nonzero $\hbar$-degrees in $\vec{u}$ do not contribute to the semiclassical limit). Indeed, by Proposition~\ref{prop::minor_basic}, we have
					\[
					L_{k \dots l}^{k \dots l}(u_k) = \sum\limits_{i=k}^l (-1)^{i-k} L_{k \dots \hat{i} \dots l}^{k+1 \dots l}(u_k-\hbar) L_{ki}(u_k).
					\]
					Observe that for every $i$,
					\[
					L_{ki}(u_k) wP = \hbar \ad_{E_{ki}}(w) P + wL_{ki}(u_k)P = \hbar \ad_{E_{ki}}(w)
					\]
					by the left invariance of the Kirillov projector from Theorem~\ref{thm::proj_left_inv}. Also observe that for every~${\xi \in \n_-}$,
					$
					Pv \xi^{\psi} = - \hbar P \ad_{\xi}(v)
					$
					by the right invariance Theorem~\ref{thm::proj_n_inv}. Therefore, it is enough to consider the terms $L_{k \dots \hat{i} \dots l}^{k+1 \dots l}(u_k-\hbar)$ in the quotient \smash{$\n_-^{\psi} \backslash \U(\m_N)$}. Again, by Proposition~\ref{prop::minor_basic},
					\[
					L_{k \dots \hat{i} \dots l}^{k+1 \dots l}(u_k-\hbar) = \sum\limits_{j=k+1}^l (-1)^{j-k-1} L_{jk}(u-\hbar l + k \hbar - \hbar) L_{k+1\dots \hat{i} \dots l}^{k+1\dots \hat{j} \dots l}(u_k - \hbar).
					\]
					In particular, all terms except for $j=k+1$ will be zero in the quotient. Therefore, the class of $L_{k \dots \hat{i} \dots l}^{k+1 \dots l}(u_k)$ is equal to that of $L_{k+1\dots \hat{i} \dots l}^{k+2 \dots l}(u_k-\hbar)$. By obvious induction, we can conclude that its class is equal to
					\begin{equation}
						\label{eq::class_quot}
						L_{k \dots \hat{i} \dots l}^{k+1 \dots l}(u_k-\hbar) \equiv L_{i+1\dots l}^{i+1\dots l}(u_k-\hbar) \mod\bigl(\n_-^{\psi} \U(\m_N)\bigr).
					\end{equation}
					Moreover,
					\[
					L_{i+1\dots l}^{i+1\dots l}(u_k-\hbar) = \prod\limits_{m=i+1}^l (u_k - u_m) + O(\hbar).
					\]
					Therefore,
					\begin{equation}
						\label{eq::action_minor}
						Pv L_{k\dots l}^{k \dots l}(u_k) x P \in \hbar \sum\limits_{i=k}^l (-1)^{i-k} \prod_{m=i+1}^l (u_k - u_m) \cdot \ad_{E_{ki}}(x) P + O\bigl(\hbar^2\bigr).
					\end{equation}
					
					We will prove the statement of the lemma by induction on $N$ in $\m_N$. The base $N=2$ follows from \eqref{eq::cg_twist_gl2}. Assume we proved it for $N-1$, in particular, we can apply it to ${}_1\m_N$. Let
					\[
					P' := P_{{}_1\m_N} (u_2,\dots, u_{N-1}).
					\]
					Then we consider
					\[
					PvPwP = PvwP + \sum\limits_{N\geq i > j \geq 2} \ad_{E_{ij}}(v) ({}_1\beta_{ij}) P_{21}^{\psi}(u_1) \cdots P_{N1}^{\psi}(u_1) wP,
					\]
					where $({}_1\beta_{ij})$ are the coefficients of the lemma for $P'$. Let us study the action of $P_{k1}^{\psi}(u_1)$. It is clear that it is enough to consider only degree one truncations
					\[
					\tilde{P}_{k1}^{\psi}(u_1) := 1 + (-1)^{k+1} \bigl(-E_{k1}^{\psi}\bigr) L_{1\dots k-1}^{1\dots k-1}(u_1)
					\]
					of $P_{k1}^{\psi}(u_1)$, since
					\[
					\prod_{i=0}^{k-1} L_{1\dots k-1}^{1\dots k-1}(u_1 - i\hbar)wP \in O\bigl(\hbar^k\bigr).
					\]
					By our assumption, each ${}_1\beta_{ij}$ is a sum of diagonal minors of the form $L_{a\dots b}^{a\dots b}(u_a)$. So, let $2 \leq a < b$ be arbitrary. We have
					\begin{gather*}
						P\ad_{E_{ij}}(v) L_{a \dots b}^{a\dots b}(u_a) \cdot (-1)^{k+1} \bigl(-E_{k1}^{\psi}\bigr) L_{1\dots k-1}^{1\dots k-1}(u_1) wP \\
						\qquad= (-1)^{k+1} P \ad_{E_{k1}} \ad_{E_{ij}}(v) L_{a\dots b}^{a\dots b}(u_a) L_{1\dots k-1}^{1\dots k-1}(u_1) wP \\
						\phantom{\qquad=}{}+ (-1)^{k+1} P \ad_{E_{ij}}(v) \ad_{E_{k1}} \bigl(L_{a\dots b}^{a \dots b}(u_a)\bigr) L_{1\dots k-1}^{1\dots k-1}(u_1)wP.
					\end{gather*}
					Observe that
					\[
					L_{a\dots b}^{a\dots b}(u_a) L_{1\dots k-1}^{1\dots k-1}(u_1) wP \in O\bigl(\hbar^2\bigr).
					\]
					Since $L_{1\dots k-1}^{1\dots k-1}(u_1)w \in O(\hbar)$, it is enough to consider the image of $\ad_{E_{k1}}\bigl(L_{a\dots b}^{a\dots b}(u_a)\bigr) $ in the quotient~${\n_- \backslash \U(\m_N) \cong \U(\b)}$.
					We have
					\[
					\ad_{E_{k1}}\bigl(L_{a\dots b}^{a\dots b}(u_a)\bigr) = (-1)^{k-a-1}\delta_{a \leq k \leq b} L_{1,a,\dots,\hat{k},\dots,b}^{a,\dots,b}(u_a).
					\]
					By Proposition~\ref{prop::minor_basic},
					\[
					L_{1,a,\dots,\hat{k},\dots,b}^{a,\dots,b}(u_a) = \sum\limits_{i=a}^b (-1)^{i-a} E_{i1} L_{a\dots \hat{k} \dots b}^{a\dots \hat{i} \dots b}(u_a-\hbar).
					\]
					Therefore, its image in the quotient is zero unless $i=a=2$, hence for $a\not= 2$,
					\begin{equation}
						\label{eq::not_two_no_contribution}
						P\ad_{E_{ij}}(v) L_{a\dots b}^{a\dots b}(u_a) \cdot (-1)^{k+1} \bigl(-E_{k1}^{\psi}\bigr) L_{1\dots k-1}^{1\dots k-1}(u_1) wP \in O\bigl(\hbar^2\bigr)
					\end{equation}
					So, we need to determine the image of \smash{$L_{2\dots \hat{k} \dots b}^{3 \dots b}(u_2-\hbar)$} in the quotient, which is $L_{k+1\dots b}^{k+1\dots b}(u_2 - (k-1)\hbar)$ by \eqref{eq::class_quot}. Therefore, up to higher-order terms, we have for $a=2$
					\begin{gather*}
						(-1)^{k+1}P\ad_{E_{ij}}(v)\ad_{E_{k1}}\bigl(L_{2\dots b}^{2\dots b}(u_2)\bigr) L_{1\dots k-1}^{1\dots k-1}(u_1)w P \\
						\qquad=P\ad_{E_{ij}}(v) L_{k+1\dots b}^{k+1\dots b}(u_2-(k-1)\hbar) L_{1\dots k-1}^{1\dots k-1}(u_1) wP.
					\end{gather*}
					
					Next, let $b\geq 1$ be arbitrary and $k>b+1$. Let us study the term
					\[
					P\ad_{E_{ij}}(v) L_{1 \dots b}^{1\dots b}(u_1) \cdot (-1)^{k+1} \bigl(-E_{k1}^{\psi}\bigr) L_{1\dots k-1}^{1\dots k-1}(u_1) wP.
					\]
					By exactly the same argument, it is enough to look at the image of
					\[
					\ad_{E_{k1}}\bigl(L_{1\dots b}^{1\dots b}(u_1)\bigr) = L_{1,2,\dots,b}^{k,2,\dots,b}(u_1) = \sum\limits_{i=1}^b (-1)^{i-1} E_{ki} L_{1\dots \hat{i} \dots b}^{2\dots b}(u_1)
					\]
					in the quotient (again, we use Proposition~\ref{prop::minor_basic}). However, by our assumption on $k$ and $b$, they are just zero. Therefore, the terms of this kind do not contribute to degree one coefficient, and the first part of the lemma follows.
					
					As for the second part, let assume first that $i>2$. It follows from \eqref{eq::not_two_no_contribution} that the coefficient~$\beta_{i+k,i}$ does not depend on $N$. Therefore, it is enough to prove the statement for $\beta_{Ni}$.
					
					By induction assumption, we have
					\begin{align*}
							{}_1\beta_{Ni}
							={}&(-1)^{N-i+1}\bigl(L_{i\dots N-1}^{i\dots N-1}(u_i) - (u_i - u_{N-1}) L_{i \dots N-1}^{i \dots N-1}(u_i)\bigr)\\
&- (u_i - u_{N-1}) \cdot {}_1\beta_{N-1,i} + {}_1\beta_{N-1,i-1}.
					\end{align*}
					It follows from arguments before that
					\[
					\beta_{Ni} = {}_1 \beta_{Ni} + \sum\limits_{k=1}^{N}	(-1)^{k+1} \ad_{E_{k1}} ({}_1 \beta_{Ni}) L_{1\dots k-1}^{1\dots k-1}(u_1).
					\]
					But by \eqref{eq::not_two_no_contribution}, $\ad_{E_{k1}}\bigl(L_{i\dots N-1}^{i \dots N-1}(u_i)\bigr)$ does not contribute to the first power. In particular,
					\[
					\ad_{E_{k1}}({}_1 \beta_{Ni}) =- (u_i - u_{N-1})\cdot \ad_{E_{k1}}({}_1\beta_{N-1,i}) + \ad_{E_{k1}}({}_1\beta_{N-1,i-1}) + O\bigl(\hbar^2\bigr),
					\]
					and the result follows.
					
					For $i=2$, observe that
					\begin{align*}
						{}_1\beta_{N2} = (-1)^{N} L_{2\dots N-1}^{2\dots N-1}(u_2), \qquad
						{}_1 \beta_{N-1,2} = (-1)^{N-1} L_{2\dots N-2}^{2\dots N-2}(u_2).
					\end{align*}
					Therefore,
					\begin{gather*}
						\beta_{N2} = (-1)^{N} L_{2\dots N-1}^{2\dots N-1}(u_2) + (-1)^{N} \sum\limits_{k=1}^{N-1} \prod_{l=k+1}^{N-1} (u_2 - u_{l}) L_{1\dots k-1}^{1\dots k-1}(u_1), \\
						\beta_{N-1,2} = (-1)^{N-1} L_{2\dots N-2}^{2\dots N-2}(u_2) + (-1)^{N-1} \sum\limits_{k=1}^{N-2} \prod_{l=k+1}^{N-2} (u_2 - u_{l}) L_{1\dots k-1}^{1\dots k-1}(u_1),\\
						\beta_{N-1,1} = (-1)^{N} L_{1\dots N-2}^{1\dots N-2}(u_1).
					\end{gather*}
					The statement follows by direct calculation of \eqref{eq::beta_relation}.
					
					Finally, for $i=1$, the statement follows from $
					\beta_{k1} = (-1)^{k+1} L_{1\dots k-1}^{1\dots k-1}(u_1)$.
				\end{proof}
			\end{Lemma}
			
			Now observe that due to \eqref{eq::action_minor}, the semiclassical limit of \eqref{eq::beta_relation} is
			\[
			\beta_{i+k,i} = -E_{i,i+k-1} - (u_i - u_{i+k-1}) \beta_{i+k-1,i} + \beta_{i+k-1,i-1}.
			\]
			But this is exactly \eqref{eq::cg_family_equation}. The theorem follows.
		\end{proof}
	\end{Theorem}

	\subsection{General setting} In this subsection, we formulate and show some properties of the Kostant--Whittaker reduction in the case of the full algebra $\gl_N$. The proof of almost all the statements can be directly translated from the mirabolic setting using the Kirillov projector and will be mostly omitted.
	
	Denote by $\Z_{\gl_N}:= \Z(\U(\gl_N))$ the center of the universal enveloping algebra of $\gl_N$. In this case, it has a very explicit presentation (for instance, see \cite{MolevNazarovOlshanskii}):
	\begin{Theorem}
		\label{thm::center}
		Let
		\[
		L_{1\dots N}^{1\dots N}(v) = v^N + \sum\limits_{i=0}^{N-1} (-1)^{N-i} A_i v^i.
		\]
		Then $\Z_{\gl_N}$ is a commutative polynomial algebra over $K$ freely generated by $\{A_i\}$.
	\end{Theorem}
	
	Denote by $\Wh(\gl_N)$ the category of $(\Z_{\gl_N},\U(\gl_N))$-bimodules which are Whittaker with respect to the right $\U(\gl_N)$-action (the definition is the same as in the mirabolic setting of Definition~\ref{def::whittaker_module}). It has a distinguished object \smash{$Q_{\gl_N} = \n_-^{\psi} \backslash \U(\gl_N)$}. In particular, we have a~functor
	\[
	\act_{\Z_{\gl_N}}^{\psi} \colon\ \BiMod{\Z_{\gl_N}}{\Z_{\gl_N}}\rightarrow \Wh(\gl_N), \qquad X \mapsto X \otimes_{\Z_{\gl_N}} Q_{\gl_N},
	\]
	where $\BiMod{\Z_{\gl_N}}{\Z_{\gl_N}}$ is the category of $\Z_{\gl_N}$-bimodules. The functor of Whittaker invariants is right adjoint
	\[
	(-)^{N_-} \colon\ \Wh(\gl_N) \rightarrow \BiMod{\Z_{\gl_N}}{\Z_{\gl_N}}.
	\]
	
	\begin{Theorem}[Kostant's equivalence]
		\label{thm::gl_skryabin}
		The functor $(-)^{N_-}$ is an equivalence.
		\begin{proof}
			The main difference with the mirabolic setting is the following fact: we need to show that for any $X \in \BiMod{\Z_{\gl_N}}{\Z_{\gl_N}}$, the unit of adjunction
$
			X \mapsto (X \otimes_{\Z_{\gl_N}} Q_{\gl_N})^{N_-}
$
			is an isomorphism. As before, the Kirillov projector gives an isomorphism
			\[
			(X \otimes_{\Z_{\gl_N}} Q_{\gl_N})^{N_-} \cong (X \otimes_{\Z_{\gl_N}} Q_{\gl_N})_{\b^{\vec{u}}}.
			\]
			By the PBW theorem, we have
			\[
			(X \otimes_{\Z_{\gl_N}} Q_{\gl_N})/\b^{\vec{u}} \cong X \otimes_{\Z_{\gl_N}} \U(\u_N)
			\]
			with a suitable $\Z_{\gl_N}$-action on $\U(\u_N)$, see Section~\ref{subsect::kirillov_proj_last_col}. But it follows from Proposition~\ref{prop::last_col_center_iso} that the map $\U(\u_N) \rightarrow \Z_{\gl_N}$ induced by the Kirillov projector is an isomorphism. Therefore, the unit of adjunction
			\[
			X \mapsto (X \otimes_{\Z_{\gl_N}} Q_{\gl_N})^{N_-} \cong (X \otimes_{\Z_{\gl_N}} Q_{\gl_N})/\b^{\vec{u}}
			\]
			is an isomorphism as well. The rest follows by the same arguments as Theorem~\ref{thm::mirabolic_skryabin}.
		\end{proof}
	\end{Theorem}
	
	Recall the category $\HC(\GL_N)$ of Harish-Chandra bimodules over $\gl_N$. We have a functor
	\[
	\act_{\gl_N}^{\psi}\colon\ \HC(G) \rightarrow \Wh(\gl_N), \qquad X \mapsto Q_{\gl_N} \otimes_{\U(\gl_N)} X.
	\]
	
	The following functor was studied in \cite{BezrukavnikovFinkelberg}.
	\begin{Definition}
		The \emph{Kostant--Whittaker reduction} functor $\res_{\gl_N}^{\psi} \colon \HC(G) \rightarrow \BiMod{\Z_{\gl_N}}{\Z_{\gl_N}}$ is the composition
		\[
		\HC(G) \xrightarrow{\act_{\gl_N}} \Wh(\gl_N) \xrightarrow{(-)^{N_-}} \BiMod{\Z_{\gl_N}}{\Z_{\gl_N}}.
		\]
	\end{Definition}
	
	Explicitly, it is given by the quantum Hamiltonian reduction $X \mapsto \bigl(\n_-^{\psi} \backslash X\bigr)^{N_-}$. There is a~natural lax monoidal structure on \smash{$\res_{\gl_N}^{\psi}$}
	\[
	\bigl(\n_-^{\psi} \backslash X\bigr)^{N_-} \otimes \bigl(\n_-^{\psi} \backslash Y\bigr)^{N_-} \rightarrow \bigl(\n_-^{\psi} \backslash X\otimes_{\U(\gl_N)} Y\bigr)^{N_-}.
	\]
	\begin{Theorem}
		\label{thm::kw_exact_monoidal}
		The functor \smash{$\res_{\gl_N}^{\psi}$} is colimit-preserving and monoidal.
	\end{Theorem}
	
	Recall the monoidal functor $\free\colon \Rep(G) \rightarrow \HC(G)$.
	\begin{Theorem}
		\label{thm::kw_trivialization}
		The Kirillov projector defines a natural isomorphism
		\[
		(P_{\m_N} (\vec{u}))_V\colon\ V \otimes \Z_{\gl_N} \rightarrow \res_{\gl_N}^{\psi}(V \otimes \U(\gl_N)).
		\]
		of \emph{right} $\Z_{\gl_N}$-modules. Similarly, there is a natural isomorphism of \emph{left} modules
		\[
		\Z_{\gl_N} \otimes V \rightarrow \res_{\gl_N}^{\psi}(V \otimes \U(\gl_N))
		\]
	\end{Theorem}
	
	\begin{Remark}
		For a description of the bimodule structure on $\res_{\gl_N}^{\psi}(V \otimes \U(\gl_N))$, see \cite{BezrukavnikovFinkelberg}. We will only need the fact that it is free on both sides.
	\end{Remark}
	
	\begin{Theorem}
		\label{thm::kw_tensor_structure}
		There is a collection of \emph{constant} maps $F_{VW}(\vec{u})$ natural in $V,W\in \Rep(G)$ such that the diagram
		\[
		\xymatrix{
			(V \otimes \Z_{\gl_N}) \otimes_{\Z_{\gl_N}} (W \otimes \Z_{\gl_N}) \ar[d] \ar[r]^-{F_{VW}(\vec{u})} & V \otimes W \otimes \Z_{\gl_N} \ar[d] \\
			\res_{\gl_N}^{\psi}(V \otimes \U(\gl_N)) \otimes_{\Z_{\gl_N}} \res_{\gl_N}^{\psi}(W \otimes \U(\gl_N)) \ar[r] & \res_{\gl_N}^{\psi}(V \otimes W \otimes \U(\gl_N))
		}
		\]
		is commutative. The maps $F_{VW}(\vec{u})$ are constant in $\Z_{\gl_N}$ and induced from the maps in Theorem~{\rm\ref{thm::mirabolic_tensor_structure}}. In particular, the corresponding $R$-matrices $R_{VW}^{\CG}(\vec{u}) := J_{WV}^{-1}(\vec{u}) J_{VW}(\vec{u})$ quantize the family of rational Cremmer--Gervais solutions \eqref{eq::cg_family_equation}.
	\end{Theorem}
	
	\subsection{Example: dual of the vector representation}
	\label{subsect:vector_rep_constant}
	
	In this subsection, we will explicitly compute the resulting rational Cremmer--Gervais $R$-matrix in the case of the dual vector representation
	\[
	V := \bigl(\C^N\bigr)^* = \mathrm{span}(\phi_1,\dots,\phi_N), \qquad \ad_{E_{ij}}(\phi_k) = -\delta_{ik} \phi_j.
	\]
	It turns out that it is closely related to the action of the \emph{degenerate affine Hecke algebra}\footnote{I would like to thank Joel Kamnitzer for the suggestion!} $\mathfrak{H}_2$ on~${V \otimes V \otimes \U(\gl_N)}$, see \cite[Example 2.1]{Cherednik} for the general construction that we will adopt to our setting in what follows.
	
	Fix the vector of parameters $\vec{u}$ as in Proposition~\ref{prop::last_col_center_iso} and denote by
	\[
	P_{\m_N} := P_{\m_N} (N\hbar-\hbar,\dots,\hbar).
	\]
	Consider the map $X\colon V \otimes \U(\gl_N) \rightarrow V \otimes \U(\gl_N)$ defined by
	\[
	X(\phi\otimes z) \mapsto -\sum_{i,j=1}^N \ad_{E_{ij}}(\phi) \otimes E_{ji} z,\qquad \phi \in V, \quad z\in \U(\gl_N).
	\]
	By definition, it commutes with the right $\U(\gl_N)$-action, and one can easily check that it commutes with the left action as well. Therefore, it defines an endomorphism of Harish-Chandra bimodules $X \in \End_{\HC(\GL_N)} (V \otimes \U(\gl_N))$.
	
	\begin{Proposition}
		\label{prop:canonical_homomorphism}
		For any $1\leq k \leq N$, we have
		\[
		P_{\m_N} \phi_k P_{\m_N} = P_{\m_N}X^{k-1} (\phi_1) P_{\m_N}.
		\]
		\begin{proof}
			We will prove it by induction. The base $k=1$ is tautological. For the step, assume that
			\[
			P_{\m_N} \phi_k P_{\m_N}= P_{\m_N} X^{k-1} (\phi_1) P_{\m_N}.
			\]
			Since $X$ commutes with the left and the right actions of $\U(\gl_N)$, it descends to an endomorphism of \smash{$\res^{\psi}_{\gl_N}(V \otimes \U(\gl_N))$}, in particular, we can apply it to both sides of the equality to get
			\[
			P_{\m_N} X(\phi_k) P_{\m_N} = P_{\m_N}X^k(\phi_1) P_{\m_N}.
			\]
			Let us show that the left-hand side is equal to $P_{\m_N}\phi_{k+1} P_{\m_N}$. Indeed,
			\begin{align*}
				P_{\m_N}X(\phi_k) P_{\m_N} &= -\sum_{i,j=1}^N P_{\m_N} \ad_{E_{ij}}(\phi_k) E_{ji} P_{\m_N} = \sum_{j=1}^N P_{\m_N} \phi_j E_{jk} P_{\m_N} \\
				&=\sum_{j>k} P_{\m_N} \phi_j E_{jk} P_{\m_N} + \sum_{j <k} P_{\m_N} \phi_j E_{jk} P_{\m_N} + P_{\m_N} \phi_k E_{kk} P_{\m_N}.
			\end{align*}
			By the properties of the Kirillov projector of Theorem~\ref{thm::main_thm}, the second summand is zero (observe that $k<N$). Likewise, the third summand is $-(N-k)\hbar \phi_k$. As for the first summand, by pushing $E_{jk}$ to the left and using the right \smash{$\n_-^{\psi}$}-invariance of $P$, we obtain
			\begin{align*}
			\sum_{j>k} P_{\m_N} \phi_j E_{jk} P_{\m_N} &= \sum_{j>k} P_{\m_N} (E_{jk}\phi_j - \hbar \ad_{E_{jk}}(\phi_j)) P_{\m_N}\\
& = P_{\m_N}\phi_{k+1} P_{\m_N} + \hbar (N-k) P_{\m_N}\phi_k P_{\m_N}.
			\end{align*}
			Therefore, the total sum is $P_{\m_N}\phi_{k+1}P_{\m_N}$ as required.
		\end{proof}
	\end{Proposition}
	
	Recall the definition of the degenerate affine Hecke algebra introduced by Drinfeld \cite{DrinfeldDaH}.
	
	\begin{Definition}
		The \emph{degenerate affine Hecke algebra} $\mathfrak{H}_k$ is a $\C[\hbar]$-algebra generated by $\{s_i\}_{i=1}^{k-1}$ and $\{X_i\}_{i=1}^k$ subject to relations
		\[
		s_i^2 = 1, \qquad s_{i+1} s_i s_{i+1} = s_i s_{i+1} s_i, \qquad s_i s_j = s_j s_i,\qquad |i-j| \neq 1,
		\]
		and
$
		s_i X_i s_i = X_{i+1} - \hbar s_i$, $ 1 \leq i \leq k-1$.
	\end{Definition}
	
	Consider the Harish-Chandra bimodule $V \otimes V \otimes \U(\gl_N)$. Recall that there is an isomorphism
	\[
	V \otimes V \otimes \U(\gl_N) \cong (V \otimes \U(\gl_N)) \otimes_{\U(\gl_N)} (V \otimes \U(\gl_N)).
	\]
	Denote by $X_i$ for $i=1,2$ the endomorphism $X$ acting on the corresponding component of the tensor product. Also, denote by
	\[
	s \colon\ V \otimes V \otimes \U(\gl_N) \rightarrow V \otimes V \otimes \U(\gl_N), \qquad s(\phi \otimes \psi \otimes z) = \psi \otimes \phi \otimes z
	\]
	the map induced from the permutation on $V \otimes V$ that is also a homomorphism of Harish-Chandra bimodules.
	
	\begin{Proposition}
		The endomorphisms $\{s,X_1,X_2\}$ define an action of $\mathfrak{H}_2$ on $V \otimes V \otimes \U(\gl_N)$.
		\begin{proof}
			We only need to check the relation $s X_1 s = X_2 - \hbar s$. Indeed, the action of the left-hand side is given by
			\begin{align*}
				\phi_k \otimes \phi_l \otimes z &\mapsto \phi_l \otimes \phi_k \otimes z \mapsto -\bigg(\sum_{i,j} \ad_{E_{ij}}(\phi_l) \otimes E_{ji}\bigg) \otimes_{\U(\gl_N)} (\phi_k \otimes z) \\
				& = -\sum_{i,j} \ad_{E_{ij}}(\phi_l) \otimes \phi_k \otimes E_{ji}z - \hbar \sum_{i,j} \ad_{E_{ij}} (\phi_l) \otimes \ad_{E_{ji}} (\phi_k) \\
				&= -\sum_{i,j} \ad_{E_{ij}}(\phi_l) \otimes \phi_k \otimes E_{ji}z - \hbar \phi_k \otimes \phi_l \otimes z \\
				&\mapsto - \phi_k\otimes \sum_{i,j} \ad_{E_{ij}}(\phi_l) \otimes E_{ji}z - \hbar \phi_l \otimes \phi_k \otimes z,
			\end{align*}
			which coincides with the action of the right-hand side.
		\end{proof}
	\end{Proposition}
	
	In the notations of Theorem~\ref{thm::kw_tensor_structure}, it follows from Proposition~\ref{prop:canonical_homomorphism} that the tensor map $F_{VV}:=F_{VV}(N\hbar-\hbar,\dots,\hbar)$ is given by
	\[
	\phi_i \otimes \phi_j \mapsto X_1^{i-1} X_2^{j-1}(\phi_1 \otimes \phi_1)\in \res^{\psi}_{\gl_N}(V \otimes V \otimes \U(\gl_N)).
	\]
	In what follows, we identify $V[\hbar] \otimes_{\C[\hbar]} V[\hbar]$ with the subspace $\C[\hbar][X_1,X_2]_{<N}$ of polynomials of degree less than $N$ in each variable.
	
	Denote by $R^{\CG}_{VV}:= R^{\CG}_{VV}(N\hbar-\hbar,\dots,\hbar)$ the corresponding $R$-matrix given by $R^{\CG}_{VV}\! =\! \sigma_{VV}^{-1} \check{R}^{\CG}_{VV}$, where \smash{$\check{R}^{\CG}_{VV} = F_{VV}^{-1} \sigma_{VV} F_{VV}$}.
	
	\begin{Theorem}
		\label{thm:the_cremmer_gervais_solution}
		Under the identification $V[\hbar] \otimes_{\C[\hbar]} V[\hbar] \cong \C[\hbar][X_1,X_2]_{<N}$, the action of $R^{\CG}_{VV}$ is
		\[
		R^{\CG}_{VV}\colon\ f(X_1, X_2) \mapsto f(X_1,X_2) - \hbar \frac{f(X_1,X_2) - f(X_2,X_1)}{X_1-X_2}.
		\]
		\begin{proof}
			Under such identification, the action of the $R$-operator $\check{R}^{\CG}_{VV}$ is simply the action of $s$ on the elements
			\[
			\{f(X_1,X_2)(\phi_1 \otimes \phi_1) \in \res^{\psi}_{\gl_N} (V \otimes V \otimes \U(\gl_N)) \mid f(X_1,X_2) \in \C[\hbar][X_1,X_2]_{<N}\}.
			\]
			Observe that $s$ acts trivially on $\phi_1\otimes \phi_1$, therefore, its action can be computed using the polynomial representation $\C[\hbar][X_1,X_2] = \mathfrak{H}_2 \otimes_{\C[S_2]} \C$ of $\mathfrak{H}_2$. Due to \cite[Proposition 1.8]{KirillovHecke}, we have
			\[
			\check{R}^{\CG}_{VV} \colon\ f(X_1,X_2) \mapsto f(X_2,X_1) - \hbar \frac{f(X_1,X_2) - f(X_2,X_1)}{X_1 - X_2}.
			\]
			We conclude by noticing that the action of $\sigma_{VV}$ is given by $f(X_1,X_2) \mapsto f(X_2,X_1)$.
		\end{proof}
	\end{Theorem}
	
	This is exactly the formula of \cite[Theorem 1.1]{EndelmanHodges} (upon identification $\kappa=-\hbar$) which, by the results of \emph{loc.cit.}, is gauge equivalent to the rational degeneration of the Cremmer--Gervais $R$-matrix.
	
	\begin{Remark}
		As the example of \cite[Proposition 2.7.13]{Kalmykov} for $N=2$ suggests, there should be an $\sl_N$-version of the Kirillov projector involving its maximal parabolic subalgebra preserving the linear span of the last basis vector. We expect that it should reproduce the actual degeneration of the Cremmer--Gervais solution.
	\end{Remark}

	\section{Application: vertex-IRF transformation}
	\label{sect::vertex_irf}
	
	The main reference for this section is author's PhD thesis \cite{Kalmykov}. For reader's convenience, we repeat the main steps of the construction and refer the reader to loc.\ cit.\  for details.
	
	\subsection{Generalized gauge transformation}
	
	Let $R_{VW}(\lambda)$ be a quantum dynamical $R$-matrix as in Section~\ref{subsect::dynamica_r_mat}.
	\begin{Definition}
		Let $S_V(\lambda) \in \End(V)$ be a collection of invertible $\End(V)$-valued rational functions on $\h^* \times \mathbb{A}^1$ natural in $\Rep(G)$. A \emph{generalized gauge transformation} is
		\[
		R_{VW}	(\lambda) \mapsto R^S_{VW}(\lambda) = (S_V(\lambda) \otimes S_W(\lambda+\hbar h_V)) R_{VW}(\lambda) \bigl(S_V(\lambda +\hbar h_W)^{-1} \otimes S_W(\lambda)^{-1}\bigr).
		\]
		A generalized gauge transformation is called a \emph{gauge transformation} if $S_V(\lambda)$ respects the $\h$-module structure on $V$.
	\end{Definition}
	
	One can similarly define a gauge transformation for dynamical twists, see \cite{EtingofSchiffmann}.
	
	\begin{Proposition}
		Gauge transformations preserve the set of quantum dynamical $R$-matrices.
	\end{Proposition}
	
	In the setting of Section~\ref{subsect::dynamica_r_mat}, it has the following categorical interpretation:
	\begin{Proposition}[{\cite[Proposition 3.8]{KalmykovSafronov}}]
		The data of a gauge transformation between dynamical twists $J_1(\lambda)$ and $J_2(\lambda)$ is equivalent to the data of a natural monoidal isomorphism
		\[
		\xymatrix@C=3cm{
			\Rep(G) \rtwocell^{J_1(\lambda)}_{J_2(\lambda)} & \HC(H).
		}
		\]
	\end{Proposition}
	
	It is also possible to interpret generalized gauge transformations in this language. There is a forgetful functor $\HC(H) \rightarrow \RMod_{\U(\h)}$. While it is not monoidal, it can be upgraded to a~\emph{monoidal action} of $\HC(H)$ on $\RMod_{\U(\h)}$.
	\begin{Definition}
		\label{def::monoidal_action}
		Let $\cC$ be a monoidal category. A \emph{$\cC$-module category} $\cM$ is a category $\cM$ with a~functor $\cM \otimes \cC \rightarrow \cM$ which we denote by $(M,X) \mapsto X \otimes M$ with a~natural isomorphisms
		\[
		\Psi_{M,X,Y}\colon\ (M \otimes X) \otimes Y \rightarrow M \otimes (X \otimes Y),
		\]
		satisfying a pentagon axiom (see \cite{EGNO}). A \emph{functor} between $\cC$-module categories $F\colon \cM_1 \rightarrow \cM_2$ is a functor of plain categories with natural isomorphisms
		\[
		\alpha_{V,X} \colon\ F(V) \otimes X \rightarrow F(V \otimes X),
		\]
		satisfying the unit and the pentagon axioms.
	\end{Definition}
	
	The Harish-Chandra category $\HC(H)$ naturally acts on the category $\RMod_{\U(\h)}$ of right $\U(\h)$-modules by
	$
	(M,X) \mapsto M \otimes_{\U(\h)} X$.
	In particular, if we have a monoidal functor $\Rep(G) \rightarrow \HC(H)$, then there is an action of~$\Rep(G)$ on $\RMod_{\U(\h)}$. One can easily see that the structure morphisms of Definition~\ref{def::monoidal_action} are given by the dynamical twist maps $J_{VW}(\lambda)$ as in Section~\ref{subsect::dynamica_r_mat}.
	
	Let $A$ be an algebra over $\C[\hbar]$ with a map $A\rightarrow \U(\h)$. Assume that there is an action of~$\Rep(G)$ on $\RMod_{A}$ such that $(A,V) \mapsto A \otimes V$ from Definition~\ref{def::monoidal_action} is a free right $A$-module. Moreover, assume that the structure maps $\Psi_{A,V,W}\colon (A \otimes V) \otimes W \rightarrow A \otimes (V \otimes W)$ are given by $\id_A \otimes F_{VW}$ for some $\C[\hbar]$-\emph{constant} $F_{VW}$. One example is $A=\C[\hbar]$; another is the setup of Theorem~\ref{thm::kw_tensor_structure} with $A = \Z_{\gl_N}$.
	
	There is an extension of scalars functor $- \otimes_A \U(\h) \colon \RMod_{A} \rightarrow \RMod_{\U(\h)}$.
	\begin{Theorem}[{\cite[Proposition 2.4.3]{Kalmykov}}]
		\label{thm::vertex-irf}
		A generalized gauge transformation between $J_{VW}(\lambda)$ and $F_{VW}$ is equivalent to the data of a $\Rep(G)$-module structure on the functor $-\otimes_A \U(\h)$.
	\end{Theorem}
	
	\subsection{Categorical vertex-IRF transformation}
	
	Recall a monoidal functor $\Rep(\GL_N) \rightarrow \HC(H)^{\gen}$ from Section~\ref{subsect::parabolic_reduction} obtained from parabolic restriction whose monoidal structure morphisms are given by the standard dynamical twist of Definition~\ref{def::standard_dyn_twist}. In particular, there is an action of $\Rep(\GL_N)$ on the category $\RMod_{\U(\h)^{\gen}}$ of right $\U(\h)^{\gen}$-modules, whose structure morphisms are given by the standard dynamical twist. Likewise, there is an action of $\Rep(\GL_N)$ on the category $\RMod_{\Z_{\gl_N}}$ of right modules over the center $\Z_{\gl_N}$ whose structure morphisms are given by the rational Cremmer--Gervais twist of Definition~\ref{def::cremmer_gervais}. Also recall the Harish-Chandra homomorphism $\A_{\gl_{N}} \rightarrow \U(\h)$. To show that there is a vertex-IRF transformation between them, it is enough to prove by Theorem~\ref{thm::vertex-irf} (or, rather, generic version thereof) that the natural functor $- \otimes_{\Z_{\gl_N}} \U(\h)^{\gen}\colon \RMod_{\Z_{\gl_N}} \rightarrow \RMod_{\U(\h)^{\gen}}$ is a functor of module categories.
	
	Recall that the parabolic restriction $\res_{\gl_N}$ is a functor $\HC(\GL_N) \rightarrow \HC(H)^{\gen}$. The Harish-Chandra category $\HC(H)^{\gen}$ is a subcategory of $\BiMod{\U(\h)^{\gen}}{\U(\h)^{\gen}}$, and there is restriction functor from the latter to $\BiMod{\Z_{\gl_N}}{\U(\h)^{\gen}}$. In what follows, we consider $\res_{\gl_N}$ to take values in $\BiMod{\Z_{\gl_N}}{\U(\h)^{\gen}}$.
	\begin{Proposition}[{\cite[Corollary 2.7.7]{Kalmykov}}]
		\label{prop::two_res_iso}
		There is a natural isomorphism
		\[
		\res_{\gl_N} \rightarrow \res_{\gl_N}^{\psi} \otimes_{\Z_{\gl_N}} \U(\h)^{\gen}
		\]
		of functors $\HC(G) \rightarrow \BiMod{\Z_{\gl_N}}{\U(\h)^{\gen}}$.
		\begin{proof}
			Let $X\in \HC(\GL_N)$. Consider a sequence of maps
			\begin{equation}
				\label{eq::two_res_comp}
				\bigl(\n_-^{\psi} \backslash X\bigr)^{N_-} \otimes_{\Z_{\gl_N}} \U(\h) \xrightarrow{p_X \otimes \id_{\U(\h)}} \bigl(\n_-^{\psi} \backslash X / \n\bigr) \otimes_{\Z_{\gl_N}} \U(\h) \xrightarrow{\act} \n_-^{\psi} \backslash X / \n,
			\end{equation}
			where $p_X$ is the projection and $\act$ is the right action of $\U(\h)$ on $\n_-^{\psi} \backslash X / \n$. By \cite[Lemma~6.2.1]{GinzburgKazhdan}, it is an isomorphism. On the other hand, there is a map
			\begin{equation}
				\label{eq::two_res_invert}
				(X/\n)^{N} \rightarrow \n_-^{\psi} \backslash X / \n.
			\end{equation}
			We claim that after extension of scalars to $\U(\h)^{\gen}$, it is an isomorphism. Indeed, since $\HC(\GL_N)$ is generated by free Harish-Chandra bimodules by Proposition~\ref{prop::hc_generators} and both functors in question are colimit-preserving by Theorems~\ref{thm::res_exact_monoidal} and~\ref{thm::kw_exact_monoidal}, it is enough to check the statement for free Harish-Chandra bimodules $X = \U(\g) \otimes V$ (here, it is convenient to take \emph{left} free modules). There is an $\n$-stable filtration on $V$ with one-dimensional quotients; in particular, we can choose a (weight) basis $\{v_{\mu}\}$ and a partial order on it such that $\ad_{\n}(v_{\mu}) > v_{\mu}$.
			
			By the PBW theorem, we have a right $\U(\h)$-module isomorphism
$
			\n_-^{\psi} \backslash X /\n \cong \U(\h) \otimes V$,
			given by the generators $\{ 1\otimes v_{\mu}\}$. Likewise, the extremal projector gives an isomorphism
			\[
			P\colon\ \U(\h)^{\gen} \otimes V \rightarrow (X/\n)^{N}
			\]
			as in Theorem~\ref{thm::dynamical_trivialization} (upon identifying $V \otimes \U(\h)^{\gen} \cong \U(\h)^{\gen} \otimes V$). It follows from definition of $P$ that
			\[
			P(1\otimes v_{\mu}) = 1\otimes v_{\mu} + \sum_{\lambda > \mu} f_{\lambda} \cdot v_{\lambda} \in (X/\n) \otimes_{\U(\h)} \U(\h)^{\gen} \cong \U(\b)^{\gen} \otimes V,
			\]
			where $f_{\lambda} \in \U(\b)$ are some elements. In particular, its class in the left quotient by $\n_-^{\psi}$ is given by some upper-unitriangular transformation in the basis $\{1\otimes v_{\mu}\}$ with coefficients in $\U(\h)^{\gen}$, hence is invertible.
			
			Therefore, composing \eqref{eq::two_res_comp} with the inverse of \eqref{eq::two_res_invert}, we conclude.
		\end{proof}
	\end{Proposition}

	\begin{Proposition}
		The functor $\RMod_{\Z_{\gl_N}} \rightarrow \RMod_{\U(\h)^{\gen}},\ M \mapsto M \otimes_{\Z_{\gl_N}} \U(\h)^{\gen}$ is a functor of $\HC(\GL_N)$-module categories, in particular of $\Rep(\GL_N)$-module categories.
		\begin{proof}
			Since the category $\RMod_{\Z_{\gl_N}}$ is generated by free modules, it is enough to check compatibility of actions on $\Z_{\gl_N}$. By Proposition~\ref{prop::two_res_iso}, we have an isomorphism
			\[ \res^{\psi}_{\gl_N} (X) \otimes_{\Z_{\gl_N}} \U(\h) \rightarrow \res_{\gl_N}(X) \]
			(for brevity, we omit the index ``gen''). Therefore, we only need to check the pentagon axiom:
			\[
			\xymatrix@C=-1cm{
				& \res^{\psi}_{\gl_N}(X) \otimes_{\Z_{\gl_N}} \res^{\psi}_{\gl_N}(Y) \otimes_{\Z_{\gl_N}} \U(\h) \ar[dr] \ar[dl] & \\
				\res^{\psi}_{\gl_N}(X \otimes_{\U(\g)} Y) \otimes_{\Z_{\gl_N}} \U(\h) \ar[dr] & & (X/\n)^{N} \otimes _{\U\h} (Y/\n)^{N} \ar[dl] \\
				& (X\otimes_{\U\g} Y/\n)^{N}. &
			}
			\]
			In terms of elements, let \smash{${}_{\n^{\psi}_-} [y] \in \bigl(\n_-^{\psi} \backslash Y\bigr)^{N_-}$} and \smash{${}_{\n^{\psi}_-} [x] \in \bigl(\n_-^{\psi} \backslash X\bigr)^{N_-}$}. We denote their image in the double quotient \smash{$\n_-^{\psi} \backslash Y / \n$} by \smash{${}_{\n^{\psi}_-} [y]_{\n}$} (same for $X$). Let $F$ be the inverse of \eqref{eq::two_res_invert}. Then we need to show that
			\[ F({}_{\n^{\psi}_-} [x]_{\n}) \otimes_{\U\g} F({}_{\n^{\psi}_-} [y]_{\n}) = F({}_{\n^{\psi}_-} [x\otimes_{\U(\g)} y]_{\n}).\]
			We show that the projections of both sides to \smash{$\n_-^{\psi} \backslash (X \otimes_{\U(\g)} Y) /\n$} coincide. Indeed, since $F$ is an isomorphism, we have \smash{$ F({}_{\n^{\psi}_-} [y]_{\n}) \in [y]_{\n} + \n_-^{\psi} Y/\n$}.
			The element \smash{${}_{\n_-^{\psi}}[x]$ is $\n_-^{\psi}$}-invariant, that is, \smash{$x\cdot \n_-^{\psi} \in \n_-^{\psi} X$}, therefore, the left-hand side belongs to
			\[[x\otimes_{\U(\g)} y]_{\n} + \n_-^{\psi} (X \otimes_{\U\g} Y)/\n,\]
			and so the projections to both sides coincide.
		\end{proof}		\label{prop::kostantwhittaker_parabolic_monoidal_equiv}
	\end{Proposition}
	
	In particular, it implies the following.
	\begin{Theorem}
		\label{thm::the_vertex-irf}
		There is generalized gauge transformation between the standard quantum dynamical $R$-matrix $R^{\dyn}(\lambda)$ from Definition~{\rm\ref{def::standard_dyn_twist}} and the constant quantum rational Cremmer--Gervais $R$-matrix $R^{\mathrm{CG}}(\vec{u})$ from Definition~{\rm\ref{def::cremmer_gervais}}.
	\end{Theorem}
	
	In the next subsection, we will compute the categorical construction in the particular case of the dual vector representation $V = \bigl(\C^N\bigr)^*$ and the vector of parameters $\vec{u} = (N\hbar - \hbar,\dots,\hbar)$.
	
	\subsection{Vertex-IRF transformation}
	
	In what follows, we will use the notations of Example~\ref{ex:vector_rep_dynamical} and Section~\ref{subsect:vector_rep_constant}; in particular, we denote $\lambda_i = E_{ii} + \hbar(N-i)$. Recall that the canonical homomorphism $X$ commutes with both $\U(\gl_N)$-actions on $V \otimes \U(\gl_N)$, in particular, it induces an action on the parabolic reduction~$\res_{\gl_N}(V \otimes \U(\gl_N))$.
	
	\begin{Lemma}
		\label{lm:canonical_homomorphism_parabolic}
		The action of $X$ on $\res_{\gl_N}(V \otimes \U(\gl_N)) \cong V \otimes \U(\h)^{\gen}$ is given by
		\[
		P \phi_i P \mapsto P\phi_i P \cdot \lambda_i,
		\]
		where $P$ is the extremal projector of Section~{\rm\ref{subsect::extremal_proj}}.
		\begin{proof}
			We have
			\begin{align*}
				PX(\phi_i)P &= -\sum\limits_{a,b} P \ad_{E_{ab}}(\phi_i) E_{ba} P = \sum\limits_b P \phi_b E_{bi} P \\
&= \sum\limits_{b>i} P \phi_b E_{bi} P + P\phi_i P E_{ii} + \sum\limits_{b<i} P \phi_b E_{bi} P \\
				&= \sum\limits_{b>i} P E_{bi} \phi_b P +\hbar (N-i) P\phi_i P + P\phi_i P E_{ii} + 0 = P \phi_i P (E_{ii}+N\hbar - i\hbar) \\
				&= P\phi_i P \cdot \lambda_i,
			\end{align*}
			where we used the properties of the extremal projector from Theorem~\ref{thm::extremal_proj}.
		\end{proof}
	\end{Lemma}
	
	Denote by
	\[
	\phi_i^{\psi} := \phi_i P_{\m_N}(N\hbar-\hbar,\dots,\hbar) \in \res^{\psi}_{\gl_N}(V \otimes \U(\gl_N)), \qquad \bar{\phi}_i := P \phi_i \in \res_{\gl_N}(V \otimes \U(\gl_N))
	\]
	the basis elements of the corresponding reductions. Recall that by construction of Proposition~\ref{prop::two_res_iso}, we need to explicitly compute the isomorphism
	\begin{align}
		\res^{\psi}_{\gl_N}(V \otimes \U(\gl_N)) \otimes_{\Z_{\gl_N}} \U(\h)^{\gen}& \rightarrow \n_-^{\psi} \backslash V \otimes \U(\gl_N)^{\gen} / \n\nonumber \\
&\rightarrow \res_{\gl_N}(V \otimes \U(\gl_N)).\label{eq::sequence_of_isos}
	\end{align}
	We identify the middle module with $V \otimes \U(\h)^{\gen}$ using PBW theorem.
	
	\begin{Lemma}
		Under the projection map $\res_{\gl_N}(V \otimes \U(\gl_N)) \rightarrow \n_-^{\psi} \backslash V \otimes \U(\gl_N) / \n$, we have
		\[
		\bar{\phi}_i \mapsto \sum\limits_{j=i}^{N} \phi_j \prod_{k=i+1}^j (\lambda_i-\lambda_k)^{-1}.
		\]
		\begin{proof}
			Recall the formula \eqref{eq::n_invariant_vectors}. It is clear that in order to be nonzero in the left quotient by \smash{$\n_-^{\psi}$}, the indices $\{i_1 < \dots < i_s\}$ should be $\{i+1 < \dots < j\}$ for some $j \leq N$. The result follows.
		\end{proof}
	\end{Lemma}
	
	In what follows, we identify the elements $\bigl\{\phi_i^{\psi}\bigr\}$ as well as $\bigl\{\bar{\phi}_i\bigr\}$ with their images in the double quotient \smash{$\n_-^{\psi} \backslash V \otimes \U(\gl_N) /\n$}.
	
	\begin{Lemma}
		We have
		\[
		\phi_1^{\psi} = \sum\limits_{i=1}^N \bar{\phi}_i \prod_{j<i} (\lambda_i - \lambda_j)^{-1}.
		\]
		\begin{proof}
			We need to show that the coefficient in front of each $\phi_k$ for $k=2,\dots,N$ on the right-hand side is zero. It is equal to \smash{$\sum_{i=1}^k \prod_{j\neq i, j\leq k} (\lambda_i-\lambda_j)^{-1}$}. But this is the sum of all residues with respect to variable $u$ of the rational function $\prod_{j=1}^k (u-\lambda_j)^{-1}$ which is zero.
		\end{proof}
	\end{Lemma}
	
	\begin{Corollary}
		For any $k\leq N$, we have
		\[
		\phi_k^{\psi} = \sum\limits_{i=1}^N \bar{\phi}_i \cdot \lambda_i^{k-1} \prod_{j<i} (\lambda_i - \lambda_j)^{-1}.
		\]
		\begin{proof}
			Observe that the sequence \eqref{eq::sequence_of_isos} commutes with the canonical homomorphism $X$. The corollary follows from Proposition~\ref{prop:canonical_homomorphism} and Lemma~\ref{lm:canonical_homomorphism_parabolic}.
		\end{proof}
	\end{Corollary}
	
	In particular, we see that the transition matrix from the Kostant--Whittaker reduction of~$V \otimes \U(\gl_N)$ to its parabolic reduction is essentially given by the Vandermonde matrix
	\begin{equation}
		\label{eq::vertex_irf_transformation}
		S(\lambda) = D(\lambda) \cdot V(\lambda), \qquad V(\lambda) = \bigl(\lambda_i^{j-1}\bigr)_{ij}, \qquad D(\lambda) = \operatorname{diag}\Big(\prod_{j<i} (\lambda_i-\lambda_j)^{-1} \Big)_{i=1}^N,
	\end{equation}
	which agrees (up to a diagonal factor) with a rational analog of the original construction \cite{CremmerGervais} (see also \cite{BuffenoirRocheTerras}) studied in \cite{BalogDabrowskiFeher}.
	
	\begin{Theorem}
		\label{thm:the_vertex_irf_transformation}
		The generalized gauge transformation between the functors of the Kostant--Whittaker reduction and parabolic restriction for the Harish-Chandra bimodules $\bigl(\C^N\bigr)^* \otimes \U(\gl_N)$ gives the vertex-IRF transformation $S(\lambda)$ as in \eqref{eq::vertex_irf_transformation} between the standard solution $R^{\dyn}(\lambda)$ to the quantum dynamical Yang--Baxter equation and the rational degeneration $\R^{\CG}$ of the Cremmer--Gervais solution as in Theorem~{\rm\ref{thm:the_cremmer_gervais_solution}}
		\[
		R^{\CG} = \bigl(S(\lambda)^{-1} \otimes S(\lambda+\hbar h_{(1)})^{-1} \bigr) R^{\dyn}(\lambda) (S(\lambda+\hbar h_{(2)}) \otimes S(\lambda)).
		\]
	\end{Theorem}
	
\subsection*{Acknowledgements} The author would like to thank Roman Bezrukavnikov, Pavel Etingof, Boris Feigin, Michael Finkelberg, Joel Kamnitzer, Vasily Krylov, and Leonid Rybnikov for helpful discussions and explanations as well as the anonymous referees for the comments. The author would also like to thank the contributors of \cite{SageMath}; this project would not be possible without their libraries.
	
This work was initially supported by the Ministry of Science and Higher Education of the Russian Federation (agreement no. 075–15–2022–287); the major part of it was accomplished at the Massachusetts Institute of Technology. The author is very grateful to the Department of Mathematics of MIT for hospitality and for the opportunity to avoid (a form of) politically motivated persecution in Russia.
	
\pdfbookmark[1]{References}{ref}
\LastPageEnding

\end{document}